\documentclass[preprint,12pt,3p]{elsarticle}
\setcitestyle{authoryear,round}
\usepackage{amssymb}
\usepackage{amsthm}
\usepackage{lineno,amsthm,amsmath,amssymb,mathrsfs,makecell,bm}
\usepackage{algorithm,algorithmicx}
\usepackage{algpseudocode}
\usepackage{subfigure}
\usepackage{booktabs,array,multirow,threeparttable,caption}
\usepackage{float}
\usepackage{subfigure}
\usepackage{graphicx}
\usepackage{epstopdf}
\usepackage{color}

\usepackage[colorlinks,
linkcolor=blue,
anchorcolor=blue,
citecolor=blue]{hyperref}

\begin{document}

\begin{frontmatter}

\title{Persistent homology-based explicit topological control for 
		2D topology optimization with MMA}

\author[zju]{Gengchen Li\fnref{equal}}
\author[zju]{Depeng Gao\fnref{equal}}
\author[zju]{Wenliang Yin}
\author[zju]{Hongwei Lin\corref{cor}}
\ead{hwlin@zju.edu.cn}
\cortext[cor]{Corresponding author.}
\address[zju]{School of Mathematical Sciences, Zhejiang University, Hangzhou, 310058, China}
\fntext[equal]{These authors contribute equally to this study.}

\begin{abstract}
Controlling structural complexity, particularly the number of holes, 
	remains a fundamental challenge in topology optimization, 
	with significant implications 
	for both theoretical analysis and manufacturability. 
Most existing approaches rely on indirect strategies, 
	such as filtering techniques, minimum length-scale control, 
	or specific level-set initializations, 
	which influence topology only implicitly and 
	do not allow precise regulation of topological features. 
In this work, we propose an explicit and differentiable topology-control framework 
	by integrating persistent homology 
	into the classical minimum-compliance topology optimization problem. 
The design domain and density field 
	are represented using non-uniform rational B-splines (NURBS), 
	while persistence diagrams are employed to rigorously and quantitatively 
	characterize topological features. 
Given a prescribed number of holes, 
	a differentiable topology-aware objective is 
	constructed from the persistence pairs and 
	incorporated into the compliance objective, 
	leading to a unified optimization formulation. 
The resulting problem is efficiently solved using the method of moving asymptotes (MMA).
Numerical experiments demonstrate that the proposed approach enables explicit control 
	over structural connectivity and the number of holes, 
	thereby providing a systematic and mathematically grounded strategy 
	for topology regulation.  
\end{abstract}

\end{frontmatter}

\section{Introduction}	
\label{sec:intro}
Topology optimization \citep{bendsoe1988generating}, 
	first introduced by Bends$\varnothing$e and Kikuchi, 
	has evolved into a powerful design paradigm for determining 
	optimal material distributions within prescribed domains 
	to achieve superior structural performance. 
It has been widely adopted across numerous engineering applications. 
Despite these advances, several challenges remain, 
	among which the control of structural complexity is of particular importance \citep{zuo2022explicit}. 

Structural complexity is particularly critical from both theoretical 
	and practical perspectives. 
From a theoretical standpoint, continuum topology optimization is known to be ill-posed, 
	as the global optimum may contain infinitely many holes \citep{sigmund1998numerical}. 
Imposing constraints on the number of holes regularizes the problem 
	and contributes to restoring well-posedness. 
From a manufacturing perspective, excessively complex designs are often difficult 
	or even infeasible to fabricate using conventional processes. 
Although additive manufacturing alleviates some of these limitations, 
	certain technologies, such as selective laser sintering 
	and fused deposition modeling, 
	still prohibit the fabrication of enclosed internal cavities. 
Consequently, controlling the number of holes is essential for 
	ensuring manufacturability 
	and maintaining structural reliability.

Existing approaches for controlling structural complexity can generally be classified into 
	implicit or explicit strategies.
Implicit methods do not directly regulate the number of holes, 
	but instead influence topology indirectly through auxiliary mechanisms.
For instance, sensitivity and density filters \citep{bourdin2001filters,sigmund1994design} 
	reduce numerous small holes as a by-product while suppressing checkerboard patterns.
Geometric constraints \citep{zhou2015minimum} restrict the minimum feature size 
	rather than explicitly controlling the number of holes.
The classical level-set method typically reduces the number of holes monotonically, 
	making the final topology strongly dependent on the initial configuration \citep{wang2003level}.
Virtual temperature approaches \citep{liu2015identification,yamada2022topology,li2018topology} 
	can enforce structural connectivity but lack an explicit relationship between 
	the control variables and the resulting number of holes. 
Consequently, these methods offer limited precision in topology control.

In contrast, explicit approaches aim to establish a direct relationship between 
	the number of holes and the design variables.
The moving morphable component method (MMC) \citep{zhang2017structural} enables  
	direct manipulation of holes through geometric components, 
	while skeleton-based level-set formulations \citep{zhang2017explicit} 
	assign an independent level-set function to each hole.  
Discrete counting formulas, such as the digital Gauss–Bonnet relation \citep{han2021topological}, 
	as well as graph-theory-inspired approaches \citep{zhao2020direct}, 
	have also been incorporated into BESO-type frameworks.

In this paper, we develop a persistent-homology-based framework that enables explicit 
	and differentiable control of structural topology within a density-based optimization setting.
Persistent homology provides a rigorous mathematical tool for quantifying topological features, 
	allowing holes and connectivity to be measured directly 
	and incorporated into the optimization process. 
By constructing differentiable objectives from persistence diagrams, 
	the proposed method can be seamlessly integrated into 
	gradient-based solvers, such as MMA. \citep{svanberg1987method}.
To summarize, the main contributions of this work are as follows:
\begin{itemize}
	\item An explicit topology-control strategy based on persistent homology 
		that accurately captures and regulates the number of holes and connectivity.   
	\item A differentiable persistence-based objective that integrates 
		naturally with compliance minimization and enables efficient gradient-based optimization.	
	\item A scalable framework that is readily extendable to three-dimensional topology optimization.

\end{itemize}

The remainder of this paper is organized as follows. 
In Section \ref{sec:pre},
	NURBS representation and the preliminaries on persistent homology 
	are introduced. 
The topology control method based on persistent 
	homology is presented in Section \ref{sec:method}. 
In Section \ref{sec:example}, 
	the experimental results are presented to validate our proposed method.
Finally, Section \ref{sec:conclud} concludes this study.

\section{Preliminaries}
\label{sec:pre}

\subsection{NURBS representation}
In isogeometric analysis (IGA), the same basis functions are employed for 
	both geometric modeling and numerical analysis, 
	enabling seamless integration between computer-aided design (CAD) 
	and computer-aided engineering (CAE). 
Non-uniform rational B-splines (NURBS) are among the most widely used representations 
	in CAD due to their ability to exactly describe common engineering geometries 
	such as conic sections and free-form surfaces.

Let $\left\{\xi_{1},\cdots,\xi_{n+p+1}\right\}$ be a non-decreasing sequence of real numbers, 
   referred to as a knot vector. 
Here, $p$ denotes the degree of the B-spline basis functions, 
   and $n$ represents the number of basis functions. 
A knot vector is called uniform if its knots are equally spaced.
An open uniform knot vector is a uniform knot vector with $p+1$ repeated knots at each end.
In the following analysis and experiments, 
	open uniform knot vectors are employed throughout.

A univariate B-spline basis function is defined recursively. 
For $p=0$:
\begin{align}
	\nonumber
	\begin{aligned}
		N_{i,0}(\xi)&=
		\begin{cases}
			1,\quad \text{if }\xi_{i}\leq\xi<\xi_{i+1}, \\
			0,\quad \text{otherwise}.
		\end{cases}
	\end{aligned}
	\quad
\end{align}
For $p>0$:
\begin{equation}
	N_{i,p}(\xi)=\frac{\xi-\xi_{i}}{\xi_{i+p}-\xi_{i}}N_{i,p-1}(\xi)+\frac{\xi_{i+p+1}-\xi}{\xi_{i+p+1}-\xi_{i+1}}N_{i+1,p-1}(\xi).
\end{equation}
The NURBS basis function of degree $p$ is defined as
\begin{equation}
	\label{a1}
	R_{i,p}(\xi)=\frac{N_{i,p}(\xi)\omega_{i}}{\sum_{i'=1}^{n}N_{i',p}(\xi)\omega_{i'}}\nonumber,
\end{equation}
where $\omega_{i}$ is the weight associated with the B-spline basis function $N_{i,p}(\xi)$. 

A NURBS curve of degree $p$ is defined as 
\begin{equation}
	\label{eq:nurbs curve}
	C(u)=\sum_{i=0}^{n}\frac{\omega_{i}N_{i,p}(u)P_{i}}{\sum_{i'=0}^{n}\omega_{i'}N_{i',p}(u)},	
\end{equation}
where $P_{i}$, $i\in\left\{0,1,\cdots,n\right\}$, 
	are the control points, 
	and $\omega_{i}$, $i\in\left\{0,1,\cdots,n\right\}$ are the 
	corresponding weights\citep{piegl2012nurbs}.

Higher-dimensional extensions of these spline concepts are constructed via tensor products.
For example, a bivariate NURBS surface of degree $(p,q)$ is defined as: 
\begin{equation}
	\label{eq:b-spline surface}
	S(u,v)=\sum_{i=0}^{m}\sum_{j=0}^{n}\frac{\omega_{i,j}N_{i,p}(u)N_{j,q}(v)P_{i,j}}{\sum_{i'=0}^{m}\sum_{j'=0}^{n}\omega_{i',j'}N_{i',p}(u)N_{j',q}(v)},	
\end{equation}
	where $P_{i,j}$, $i\in\left\{0,1,\cdots,m\right\}$, 
	$j\in\left\{0,1,\cdots,n\right\}$ are the control points 
	in the $u$, $v$ directions, respectively;  
	$N_{i,p}(u)$, $N_{j,q}(v)$ are the B-spline 
	basis functions of  
	degree $p$ and $q$ in the corresponding parametric directions \citep{piegl2012nurbs}.
In this paper, the structural geometry is represented by a
	bivariate tensor-product NURBS surface.



\subsection{Persistent homology of cubical complex}
Persistent homology is a fundamental tool in Topological Data
	Analysis (TDA) for extracting topological features of datasets
	filtered by real-valued functions. 
To characterize the topology of the structure under consideration, 
	the continuous density field is discretized, 
	thereby inducing a well-defined topological representation.
In the two-dimensional setting, the discretized density field 
	naturally results in image-like data.
For such image-based representations, the cubical complex provides 
	a particularly suitable computational framework.

\begin{figure}[!htb]
	\centering
	\subfigure[0-cube]{
		\includegraphics[width=0.15\linewidth]{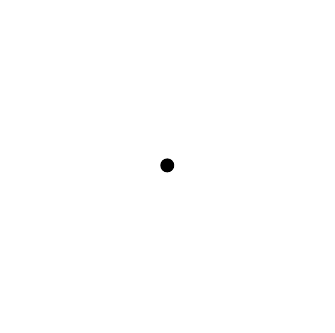}
	}
	\subfigure[1-cube]{
		\includegraphics[width=0.15\linewidth]{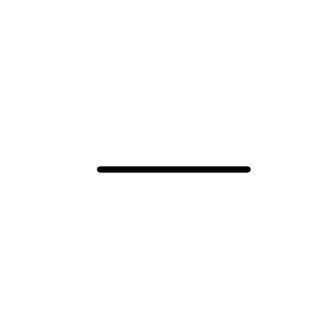}
	}
	\subfigure[2-cube]{
		\includegraphics[width=0.15\linewidth]{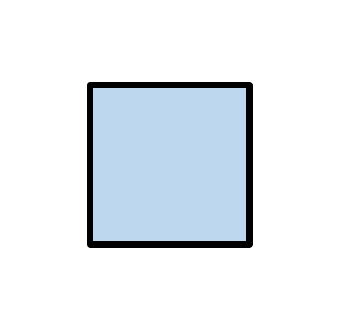}
	}
	\subfigure[3-cube]{
		\includegraphics[width=0.16\linewidth]{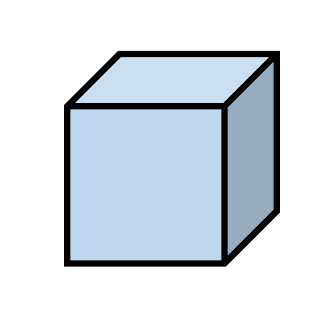}
	}\\
	\subfigure[Cubical complex]{
		\includegraphics[width=0.2\linewidth]{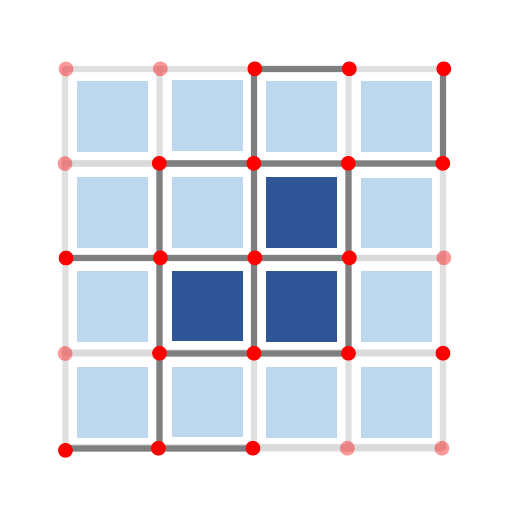}
	}
	\subfigure[Underlying space]{
		\includegraphics[width=0.2\linewidth]{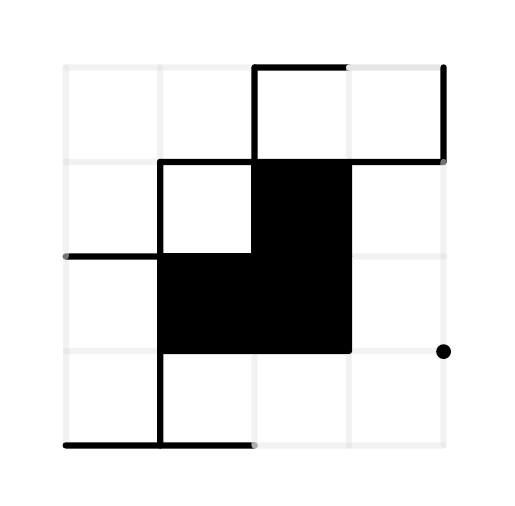}
	}
	\caption{Cubical complex. (a)-(d) Examples of $k$-dimensional cubes. 
	(e) An example of cubical complex represented by dark part, 
	which is a collection of cubes. 
	(f) The underlying space constructed from (e) by merging all 
	cubes as a topological space on the plane.}
	\label{fig:cubical_complex}
\end{figure}

In Euclidean space, an elementary interval is defined as 
	$I=[l,l+1]$, where $l\in\mathbb{Z}$.
A $k$-cube $\sigma^{k}$ is defined as the 
	Cartesian product of $k$ elementary intervals, 
	such as points ($0$-cubes), segments ($1$-cubes),
	and squares ($2$-cubes). 
A face of a cube is a subset $F\subset I^{k}$ of the form 
	$F=\prod_{i=1}^{k}J_{i}$ where $J_{i}$ is either $\{l\}$, $\{l+1\}$, 
	or $\left[l,l+1\right]$ for all $1\leq i \leq k$. 
We denote $\tau \leq \sigma^{k}$ if $\tau$ is a face of $\sigma^{k}$.
A cubical complex $\mathcal{K}$ is a collection of $k$-cubes,
	such that each face of each cube in $\mathcal{K}$ is also contained in $\mathcal{K}$.
The union of all cubes in $\mathcal{K}$ forms its underlying space $\mathcal{S}$. 
Fig. \ref{fig:cubical_complex} shows examples of cubical complexes.

In general, the topological features of a cubical complex can be classified 
	according to their dimension, namely connected components, 
	one-dimensional loops (holes) and voids corresponding to 
	$dim=0,1,2$, respectively. 
For two-dimensional topology optimization problems, 
	the structural topology is fully characterized by 
	connected components and one-dimensional loops.
More precisely, zero-dimensional topological features correspond to 
	the connected components of the structure, 
	whereas one-dimensional features represent the holes within the structure.
From the perspective of algebraic topology, 
	these features are formally described by homology groups.
For a cubical complex $(\mathcal{K})$, its $d$th homology group $H_{d}(\mathcal{K})$ 
	encodes the $d$-dimensional topological features mentioned above. 
The rank of $H_{d}(\mathcal{K})$, denoted by $\beta_{d}(\mathcal{K})$, 
	is called the $d$-th Betti number and represents the number of 
	independent $d$-dimensional homological features in $(\mathcal{K})$.
For example, $\beta_{0}(\mathcal{K})$ equals the number of connected 
	components of $\mathcal{K}$ while  
	$\beta_{1}(\mathcal{K})$ corresponds to the number of independent 
	one-dimensional loops (i.e., holes).

\begin{figure}[!htb]
	\centering
	\subfigure[]{
		\includegraphics[width=0.22\linewidth]{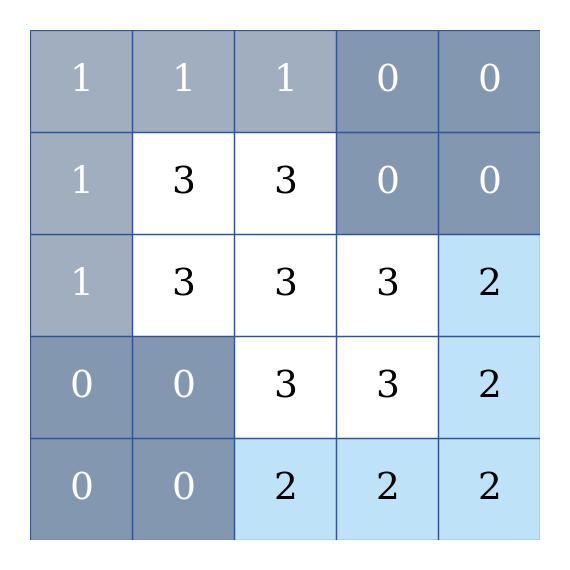}
	}
	\subfigure[]{
		\includegraphics[width=0.3\linewidth]{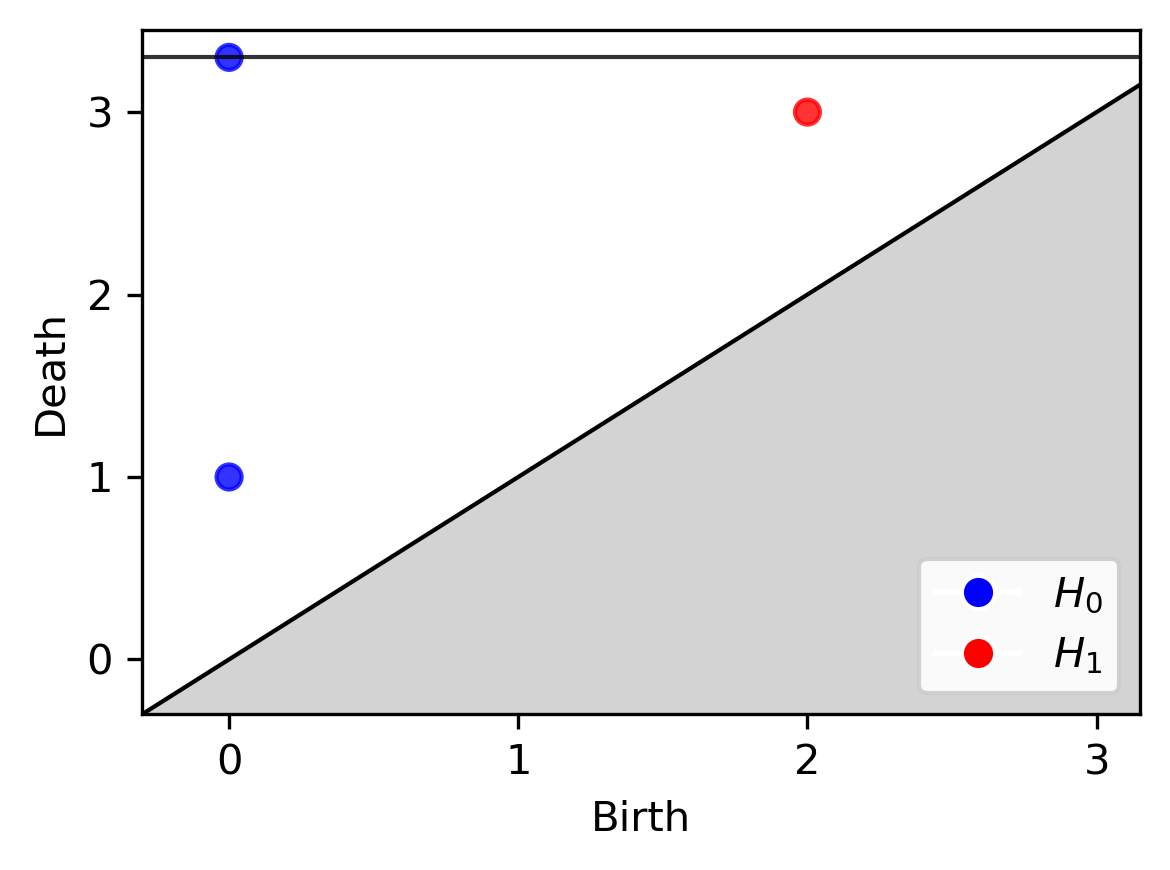}
	}
	\subfigure[]{
		\includegraphics[width=0.33\linewidth]{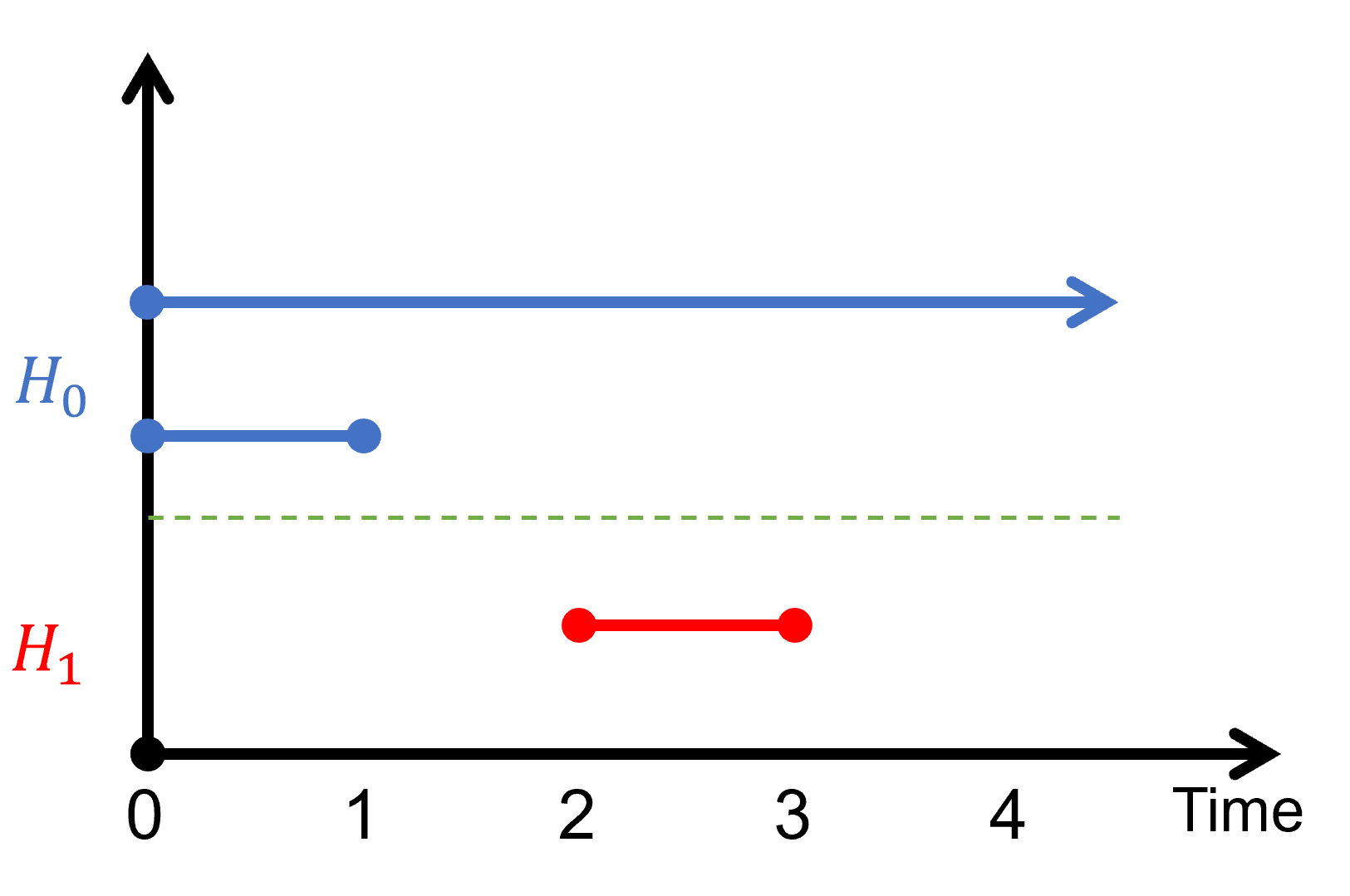}
	}\\
	\subfigure[]{
		\includegraphics[width=0.20\linewidth]{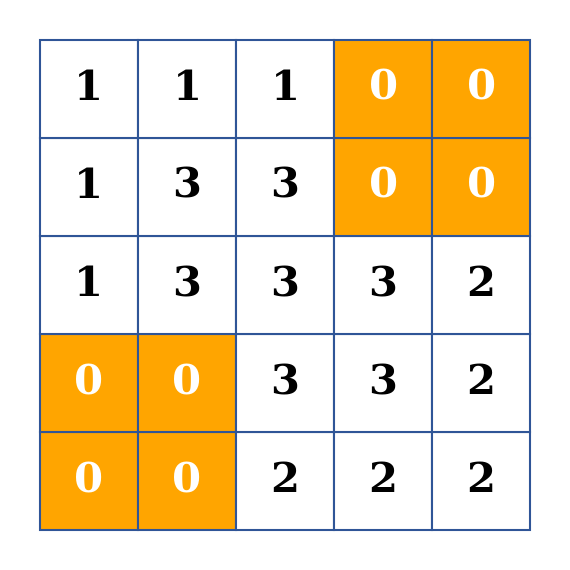}
	}
	\subfigure[]{
		\includegraphics[width=0.20\linewidth]{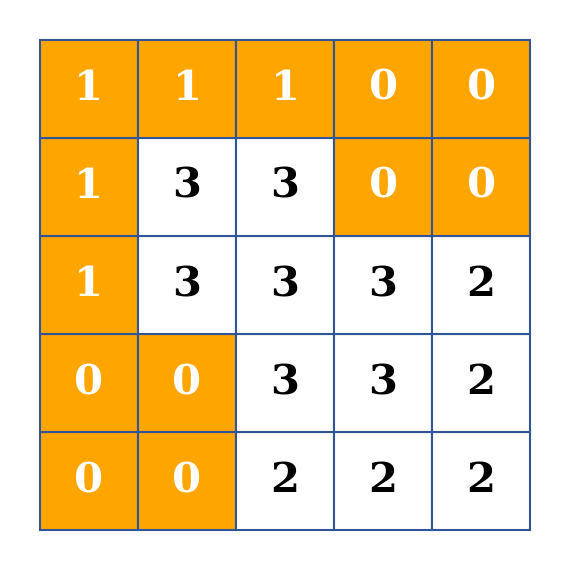}
	}
	\subfigure[]{
		\includegraphics[width=0.20\linewidth]{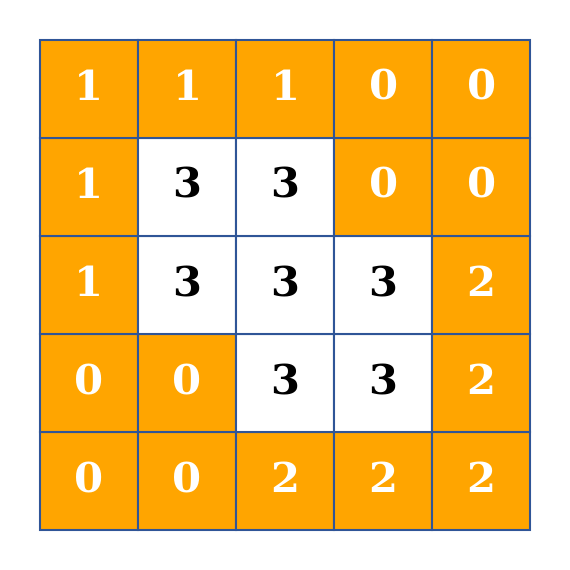}
	}
	\subfigure[]{
		\includegraphics[width=0.20\linewidth]{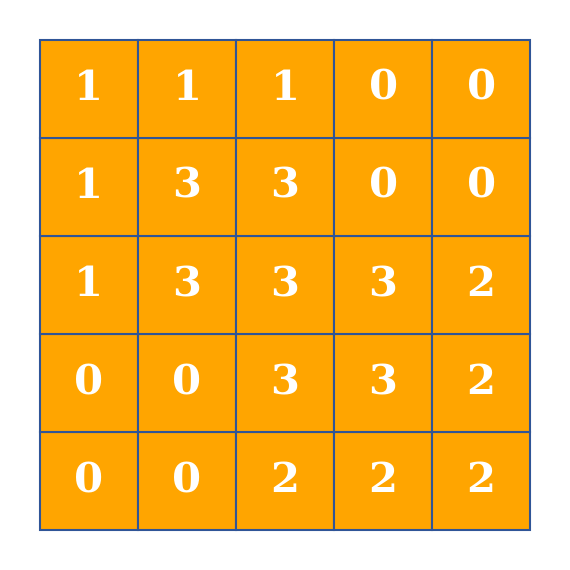}
	}
	\caption{Filtration on a cubical complex. (a) is a cubical data with 
	height function values on its grid cells. 
	(d)$-$(g) show the corresponding filtration of (a)
	and each item of this filtration is a cubical complex 
	which represents the sublevel set of original cubical data.
	The topological features can be extracted from
	those expanding cubical complexes with threshold value increasing. 
	Using persistent homology, one can obtain the persistence diagram (b) and persistence barcode (c)
	of filtration (d)–(g).}
	\label{fig:filtration}
\end{figure}

Topological features defined on discrete structures cannot be 
	directly optimized using conventional discrete metrics. 
To address this limitation, persistent homology is introduced, 
	enabling the tracking of the birth and death of topological 
	features along a sequence of nested complexes, known as a filtration. 
A filtration is constructed from a function defined on a cubical complex
	$f:\mathcal{K}\rightarrow \mathbb{R}$, 
	by considering its sublevel sets.
Let $f:\mathcal{K}\rightarrow \mathbb{R}$ be a function defined on all elementary cubes 
	of $\mathcal{K}$ that satisfies: $\forall \sigma \in \mathcal{K}$, 
	$\forall \tau \leq \sigma$, $f(\tau)\leq f(\sigma)$. 
Then, for any $a\in \mathbb{R}$, the sublevel set of $\mathcal{K}$ with respect to $f$ 
	is defined as 
	$\mathcal{K}(f,a)=f^{-1}((-\infty,a])$. 
For a series of $a_{0}\leq a_{1}\leq \cdots \leq a_{t}$, we have 
	\begin{equation}
		\nonumber
		\mathcal{K}(f,a_{0})\subseteq\mathcal{K}(f,a_{1})\subseteq \cdots \subseteq\mathcal{K}(f,a_{t})
	\end{equation}
The sublevel filtration $\mathcal{F}$ is defined as a nested 
	sequence of cubical complexes 
	$\mathcal{K}(f,a_{0})\subseteq \mathcal{K}(f,a_{1})\subseteq 
	\cdots \subseteq\mathcal{K}(f,a_{t})$, 
	where the threshold parameter satisfies $a_{0}\leq a_{1}\leq \cdots \leq a_{t}$. 
As the parameter $a_{t}$ increases, the emergence (birth) and disappearance (death) of 
	topological features in different dimensions can be systematically tracked.
The filtration therefore reveals the multi-scale  
	topological structure underlying the data.
Persistent homology computes not only the homology group of
	each complex in the filtration but also records the birth and
	death times of topological features according 
	to well-defined algebraic rules \citep{edelsbrunner2010computational}.
Through this framework, topological characteristics 
	are encoded as continuous variables instead of discrete Betti
	numbers, thereby enabling gradient-based optimization of topology.
The lifespan of each topological feature can be
	visualized as an interval along the filtration parameter axis, 
	forming what is known as a persistence barcode.
Specifically, if a topological feature appears 
	in $\mathcal{K}(f,b)$ 
	and disappears in $\mathcal{K}(f,d)$, 
	it is represented as a persistent pair $(b,d)$.
These persistent pairs $(b_{i},d_{i})$ can equivalently be mapped to points 
	$(b_{i},d_{i})$ in the Cartesian plane, 
	yielding a persistence diagram (PD), 
	which is one of the most widely used representations in TDA.
Fig. \ref{fig:filtration} illustrates the persistent pairs $(b_{i},d_{i})$ 
	represent these topological features and can be plotted 
	into a persistence diagram or a persistence barcode.

\section{Methodology}
\label{sec:method}
\subsection{Topology control via persistent homology}
PD provides a compact and rigorous representation of the 
	topological characteristics of a structure, 
	thereby enabling modifications of the underlying topology
	through its optimization. 
Within this framework, the PD is optimized 
	via gradient-descent methods \citep{poulenard2018topological}. 
Subsequent theoretical studies have rigorously established the correctness 
	and convergence properties of persistence-based optimization schemes 
	and have demonstrated that a broad class of persistence-based 
	functions can be optimized using stochastic 
	subgradient descent algorithms \citep{carriere2021optimizing}. 
Moreover, persistence-based optimization of PDs through gradient descent 
	has shown significant effectiveness across a wide range 
	of applications \citep{clough2020topological}.

\subsubsection{Filtration induced by density field}
\begin{figure}[!htb]
	\centering
	\includegraphics[scale=0.8]{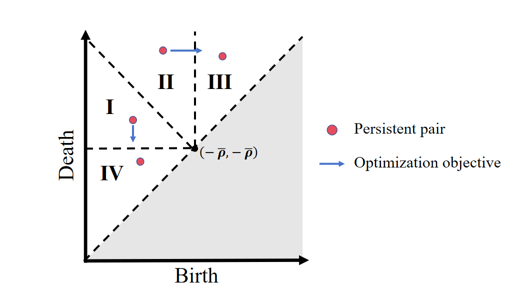}
	\caption{A two-dimensional persistence diagram. Regions I, II, III, and IV are segmented
	by lines $y=-x$, $x=-\bar{\rho}$ and $y=-\bar{\rho}$. 
	The persistence pairs in regions I and II correspond
	to the isolated holes within the binarized structure. 
	The blue arrows indicate the
	direction of the optimization objective. }
	\label{fig:density_filter}
\end{figure}

Let $\rho$ be a density field represented by NURBS: 
\begin{equation}
	\rho(u,v)=\sum_{i=0}^{m}\sum_{j=0}^{n}\frac{\omega_{i,j}N_{i,p}(u)N_{j,q}(v)P_{i,j}}{\sum_{i'=0}^{m}\sum_{j'=0}^{n}\omega_{i',j'}N_{i',p}(u)N_{j',q}(v)},	
	\label{eq:density}
\end{equation}
For the binarized structure 
\begin{equation}
	\tilde{\rho}=\left\{
	\begin{aligned}
		&1, \quad \rho\geq \bar{\rho}\\
		&0, \quad \rho< \bar{\rho}
	\end{aligned}
	\right.
\end{equation}
	where 
	$\bar{\rho}$ is the threshold, 
	and the sublevel filtration can be induced 
	by employing the equivalent structural representation $-\rho\leq -\bar{\rho}$.
Therefore, by substituting the filtration function $f$ with $-\rho$,  
	computing the cubical complexes of the discretized density field, 
	and applying the sublevel filtration mentioned above, 
	the topological features of the structure can be obtained. 
Suppose the resulting $k$-dimensional
	persistence diagrams consist of the following persistence pairs:
	\begin{equation}
	 P_{k}=\{b_{i,k},d_{i,k}\}_{i=0}^{N_{k}-1},	
	 \label{eq:per_pair}
	\end{equation}
	where $b_{i,k}$ and $d_{i,k}$ are the birth and 
	death times of the $i$-th persistence pair in 
	dimension $k$.
As shown in Fig. \ref{fig:density_filter}, 
	only those persistence pairs born before $-\bar{\rho}$
	and die after $-\bar{\rho}$ (i.e., persistence pairs 
	located in regions I and II) 
	correspond to topological features present in the structure.
Based on this correspondence between persistence pairs and topological 
	features, 
	the topology of the structure can be controlled by optimizing the PD.
The key to persistence-based optimization is to find a correspondence between 
	underlying complex and persistence pairs, that is, $(b_{i},d_{i})$.
Let $f:\mathcal{K}\rightarrow \mathbb{R}$ be a function 
	on all elementary cubes in $\mathcal{K}$ 
	such that it generates a filtration $\mathcal{F}$. 
Every homology generator gives birth at an elementary cube $\tau_{i}$ 
	and dies at another elementary cube $\sigma_{i}$. 	
Thus, persistence pair $(b_{i},d_{i})$ can be 
	written as $(f_{\tau_{i}},f_{\sigma_{i}})$. 
Suppose $\{c_{j}\}$ are the parameters to be optimized, 
	we can replace gradient $\frac{\partial b_{i}}{\partial c_{j}}$ 
	and $\frac{\partial d_{i}}{\partial c_{j}}$ with 
	$\frac{\partial f_{\tau_{i}}}{\partial c_{j}}$ 
	and $\frac{\partial f_{\sigma_{i}}}{\partial c_{j}}$ 
	due to the invariance of critical elementary cubes within a local
	parameter neighborhood.



\subsubsection{Zero-dimensional topology objective function}
The presence of extra connected components may affect manufacturability.
Therefore, we should design a zero-dimensional topology objective function 
	to ensure the structural connectivity.
A two-dimensional PD as shown in Fig. \ref{fig:density_filter} 
	is utilized to illustrate the persistence-based optimization. 
Since the death times of the topological features are always larger
	than their birth times, 
	all persistence pairs in the PD are located above the line $y=x$. 
The diagram is segmented by the lines $y=x$, $x=-\bar{\rho}$, $y=-\bar{\rho}$ and $y=-x$
	into four regions. 
In 0-dimensional situations, the set $\{-\rho \leq +\infty\}$ is
	always one connected component, which corresponds to a persistent
	pair $(b,+\infty), b\in \mathbb{R}$ in the zero-dimensional persistence diagram. 
To ensure no additional connected components exist in the structure, 
	our optimization objective is to eliminate persistent pairs 
	in regions I and II of the zero-dimensional PD.
Excess persistence pairs located in regions I and II should 
	be removed and translocated to regions III and IV. 
Specifically, persistence pairs in region I are translocated to region IV, 
	and persistence pairs
	in region II are translocated to region III to avoid a significant
	structure change caused by a considerable movement of persistence pairs. 

The 0-dimensional persistence pairs $P_{0}$ defined in Eq. (\ref{eq:per_pair}) 
	are categorized into distinct sets based on regions.
\begin{align}
	\left\{
		\begin{aligned}
			P_{0}^{I}&=\{(b_{i,0}^{I},d_{i,0}^{I})\}_{i=0}^{N_{0}^{I}-1}\\
			P_{0}^{II}&=\{(b_{i,0}^{II},d_{i,0}^{II})\}_{i=0}^{N_{0}^{II}-1},
		\end{aligned}
		\right.
\end{align}
where $N_{0}^{I}$ and $N_{0}^{II}$ denote the counts of persistence pairs within regions
	I and II, respectively. 
The repositioning of 0-dimensional persistence pairs from region I to IV and region II to III,
	respectively, is achieved by defining the objective function as
	\begin{equation}
		C_{top}^{0}=\sum_{i=0}^{N_{0}^{I-1}}d_{i,0}^{I}-\sum_{i=0}^{N_{0}^{II-1}}b_{i,0}^{II}.
		\label{eq:zero_optimization}
	\end{equation}

\subsubsection{One-dimensional topology objective function}

To control the structural complexity, 
	we impose a constraint on the maximum number of holes $\bar{N}_{1}$ in the structure.
We design the objective function from a dual perspective: by treating the 
	material as voids and voids as material. 
That is, inducing filtration on $\rho < \bar{\rho}$ to capture the holes.
Increasing the threshold $\bar{\rho}$ naturally induces a 
	corresponding filtration on the void phase.
	and then the holes in the structure can be captured by computing 
	the 0-dimensional PD. 

However, not all 0-dimensional persistence pairs 
	represent structural holes.
Only those that do not intersect the boundary of the design domain 
	are regarded as holes.
Therefore, after inducing the filtration, 
	we employ topological inverse mapping to determine the birth locations 
	of the 0-dimensional topological features in the parametric domain.
In the discretized density field, a breadth-first search (BFS) is then 
	performed from each birth location to check whether the corresponding component 
	intersects the boundary of the design domain, 
	thereby filtering out the actual structural holes.
	
We control the maximum number of holes in the structure. 
The number of holes $N_{1}$ is constrained not to exceed $\bar{N}_{1}$.
We can induce filtration based on 
	$\rho < \bar{\rho}$ and obtain the corresponding 
	0-dimensional persistence pairs of holes through BFS: 
	\begin{equation}
		\tilde{P_{0}}=\{(\tilde{b}_{i,0},\tilde{d}_{i,0})\}_{i=0}^{N_{1}-1}.
	\label{eq:flip_p}
	\end{equation}
Suppose $\{(\tilde{b}_{i,0},\tilde{d}_{i,0})\}_{i=0}^{N_{1}-1}$ are sorted in ascending order of area. 
The area of each hole is approximated by the number of discrete elements it contains.
When $N_{1}>\bar{N}_{1}$, 
	we drive the first $N_{1}-\bar{N}_{1}$ persistence pairs 
	out of regions I and II.
Similar to the zero-dimensional case, 
	we categorize these $N_{1}-\bar{N}_{1}$ persistence pairs 
	into distinct sets based on their respective regions.
\begin{align}
	\left\{
		\begin{aligned}
			\tilde{P_{0}}^{I}&=\{(\tilde{b}_{i,0}^{I},\tilde{d}_{i,0}^{I})\}_{i=0}^{(N_{1}-\bar{N}_{1})^{I}-1}\\
			\tilde{P_{0}}^{II}&=\{(\tilde{b}_{i,0}^{II},\tilde{d}_{i,0}^{II})\}_{i=0}^{(N_{1}-\bar{N}_{1})^{II}-1},
		\end{aligned}
		\right.
\end{align}
where $(N_{1}-\bar{N}_{1})^{I}$ and $(N_{1}-\bar{N}_{1})^{II}$ denote the number of 
	persistence pairs needed to be removed from regions
	I and II, respectively. 
Thus, one-dimensional topology objective function can be defined as:
\begin{equation}
	C_{top}^{1}=\sum_{i=0}^{(N_{1}-\bar{N}_{1})^{I-1}}\tilde{d}_{i,0}^{I}-
	\sum_{i=0}^{(N_{1}-\bar{N}_{1})^{II-1}}\tilde{b}_{i,0}^{II}.
	\label{eq:one_optimization}
\end{equation}

\subsection{Problem and formulation}
The minimum compliance optimization problem with 
	zero-dimensional and one-dimensional topological constraints 
	can be summarized as follows:

\begin{equation}
	\begin{aligned}
		\label{eq:optimization_problem}
		&\text{Find: }\rho=\sum_{i=0}^{m}\sum_{j=0}^{n}\omega_{i,j}N_{i,p}(u)N_{j,q}(v)\rho_{i,j}\\
		&\text{Minimize: }C(\rho)=\frac{1}{2}\mathbf{F}^{T}\mathbf{U}
							+C_{top}(\rho)\\
		&\begin{array}{ll}
		\text{Subject to:} &V(\rho)\leq V_{0}\\
								&\mathbf{KU=F}\\
								&\rho_{min}\leq \rho \leq 1,\\
		\end{array}
	\end{aligned}
\end{equation}
where $\rho_{i,j}$ are the design variables, 
	denoting the control coefficients of the density field, 
	$V(\rho)$ and $V_{0}$ represent the volume 
	fraction of the optimized structure and a 
	predefined threshold, respectively,  
	and $\rho_{min}=0.1$ is the lower bound of density to 
	prevent singular stiffness matrices.
Based on Eq. (\ref{eq:zero_optimization}) and (\ref{eq:one_optimization}), 
	topological objective function $C_{top}(\rho)$ is defined as follows:
\begin{equation}
	\centering	
	C_{top}(\rho)=\mu_{0}C_{top}^{0}+\mu_{1}C_{top}^{1}
	\label{eq:topo_opt}
\end{equation}
where $\mu_{0}$ and $\mu_{1}$ denote the weights of 0-dimensional and 1-dimensional topology 
	objective functions, respectively.
 
\subsection{Sensitivity analysis}
In this section, we perform a sensitivity analysis 
	of the objective function and the constraints 
	with respect to the design variables.

The sensitivity of $C_{top}^{0}$ and $C_{top}^{1}$ with respect to 
	$\rho_{i,j}$ are computed utilizing the topological inverse mapping \citep{poulenard2018topological,bruel2020topology} 
\begin{equation}
	\pi_{\bm{\rho}}(b_{i},d_{i})=(\pi_{b}(b_{i}),\pi_{d}(d_{i})), 
\end{equation}
where $(b_{i},d_{i})$ denotes a persistence pair, 
	$\pi_{b}(b_{i})$ and $\pi_{d}(d_{i})$ denote the locations in the 
	parameter domain that correspond to the birth and death of 
	the topological feature, respectively.
Based on the inverse mapping $\pi_{\bm{\rho}}$ and the differential 
	chain rule, the partial derivatives of $i$-dimensional
	topological objective function 
	with respect to the control coefficients of the density field 
	can be derived:
\begin{equation}
	\begin{split}
		\frac{\partial C_{top}^{i}}{\partial \rho_{i,j}}
		&=\frac{\partial C_{top}^{i}}{\partial b}\frac{\partial b}{\partial \rho_{i,j}}+
		\frac{\partial C_{top}^{i}}{\partial d}\frac{\partial d}{\partial \rho_{i,j}}\\
		&=-\sum_{(b_{i},d_{i})\in I}\frac{\partial (-\rho)}{\partial \rho_{i,j}}(\pi_{b}(b_{i}))
			+\sum_{(b_{i},d_{i})\in II}\frac{\partial (-\rho)}{\partial \rho_{i,j}}(\pi_{d}(d_{i}))\\
		&=\sum_{(b_{i},d_{i})\in I}\frac{\partial \rho}{\partial \rho_{i,j}}(\pi_{b}(b_{i}))
			-\sum_{(b_{i},d_{i})\in II}\frac{\partial \rho}{\partial \rho_{i,j}}(\pi_{d}(d_{i}))
	\end{split}
\end{equation}
Then the sensitivity of the whole objective function is:
\begin{equation}
	\begin{split}
		\frac{\partial C(\rho)}{\partial \rho_{i,j}}
		&=\frac{\partial \bm{U}^{T}\bm{KU}}{2\partial \rho_{i,j}}
			+\mu_{0}\frac{\partial C_{top}^{0}}{\partial \rho_{i,j}}
			+\mu_{1}\frac{\partial C_{top}^{1}}{\partial \rho_{i,j}}\\
		&=-\frac{1}{2}p\rho^{p-1}\bm{U}^{T}\bm{KU}\frac{\partial \rho}{\partial \rho_{i,j}}\\
		&+
		\sum_{i=0}^{1}\mu_{i}\left[\left(\sum_{(b_{i},d_{i})\in I}\frac{\partial \rho}{\partial \rho_{i,j}}(\pi_{b_{i}}(b_{i}))
		-\sum_{(b_{i},d_{i})\in II}\frac{\partial \rho}{\partial \rho_{i,j}}(\pi_{d_{i}}(d_{i}))\right)\right].
	\end{split}
	\label{eq:sens_top}
\end{equation}

The sensitivity of the volume constraint is given by:
\begin{equation}
	\begin{split}
		\frac{\partial V(\rho)}{\partial \rho_{i,j}}
		&=\frac{\partial \int_{\Omega}\rho d\Omega}{\partial \rho_{i,j}}\\
		&=\frac{\partial \int_{\Omega}\sum_{i=0}^{m}\sum_{j=0}^{n}\frac{\omega_{i,j}N_{i,p}(u)N_{j,q}(v)P_{i,j}}{\sum_{i'=0}^{m}\sum_{j'=0}^{n}\omega_{i',j'}N_{i',p}(u)N_{j',q}(v)}d\Omega}
		{\partial \rho_{i,j}}\\
		&=\int_{\Omega}\frac{\omega_{i,j}N_{i,p}(u)N_{j,q}(v)P_{i,j}}{\sum_{i'=0}^{m}\sum_{j'=0}^{n}\omega_{i',j'}N_{i',p}(u)N_{j',q}(v)}d\Omega.
	\end{split}
	\label{eq:sens_vol}
\end{equation}
\subsection{Optimization algorithm}
In summary, the proposed topology control method 
	first employs the density field $\rho$ to induce a sublevel filtration, 
	yielding the corresponding zero-dimensional persistence pairs.
To ensure the optimization stability, 
	we control the number of holes also by optimizing the 
	associated zero-dimensional topological features. 
By adopting a dual-phase strategy, 
	$\rho$ is used to induce sublevel filtration for hole detection.
Topological inverse mapping combined with BFS is then employed to identify the 
	holes in the structure, 
    resulting in zero-dimensional persistence pairs corresponding to the void phase.
When the prescribed topological constraints are violated, 
	the persistence pairs associated with redundant connected components 
	or excessive holes are optimized. 
Specifically, by minimizing their death times 
	or maximizing their birth times, 
	these undesired topological features are gradually eliminated from the structure.

Owing to the differentiability of the persistence-based objective function, 
	the proposed framework can be seamlessly integrated into 
	classical topology optimization procedures, 
	enabling rigorous and continuous control of structural complexity. 
The overall optimization procedure is summarized in Algorithm~\ref{alg:opt}: 
\begin{algorithm}[!htb]
	\caption{Optimization algorithm with topology control} 
	\label{alg:opt}
	\hspace*{0.02in} {\bf Input:} 
	Design domain represented by NURBS; \\
	\hspace*{0.56in}Boundary conditions;\\
	\hspace*{0.56in}Maximum number of holes $\bar{N}_{1}$;\\
	\hspace*{0.56in}Poisson's ratio, Young's modulus;\\
	\hspace*{0.56in}Volume fraction $V_{0}$;\\
	\hspace*{0.56in}Maximum number of iterations;\\
	\hspace*{0.02in} {\bf Output:} 
	Optimized density field represented by NURBS (\ref{eq:density}).
	\begin{algorithmic}[1]
		\While{the maximum number of iterations is not reached}
		\State Perform finite element analysis, and compute the compliance and volume sensitivities;
		\State Induce filtration on $-\rho<-\bar{\rho}$ and compute the number of connected components $N_{0}$;
		\State Induce filtration on $\rho<\bar{\rho}$ and perform BFS in the parameter domain;
        \State Compute the number of holes $N_{1}$ and sort them in ascending order of area;
		\If{$N_{0}>1$ \textbf{or} $N_{1}>\bar{N}_{1}$}
		\State Incorporate the topological objective function into the total objective function;
		\State Compute the sensitivity of the augmented objective function; 
		\EndIf	
		\If{ $N_{0}=1$ \textbf{and} $N_{1}\leq \bar{N}_{1}$ \textbf{and} excess material exists}
		\State Freeze the density of the existing holes;
		\State Compute the sensitivities to optimize the compliance of the remaining structure;
		\EndIf
		\State Update the density control coefficients using MMA;
		\EndWhile
	\end{algorithmic}
\end{algorithm}

\section{Numerical examples}
\label{sec:example}
In this section, several benchmark examples are presented to validate the effectiveness 
	of the proposed method. 
All problems are solved using MMA. 
The density field is represented by NURBS, 
	and persistent homology is evaluated at each iteration 
	to extract topological features of the structure.
The results demonstrate that the proposed framework can explicitly ensure the 
	structural connectivity and accurately control the number of holes.

Suppose the solid material has a Young's modulus of $1.9\times 10^{11}$ and 
	Poisson's ratio of $0.3$.
The filtering radius of the sensitivity filter is $1.5h$, where $h$ represents the element size.
In this study, the threshold $\bar{\rho}$ is $0.4$.
The maximum iteration number is set to $400$.
Due to structural topology instability during the initial iterations,
	we implemented topology control starting from the $20-30$th iteration in our experiments.
Moreover, the density field discretization resolution for calculating persistent homology 
	is set to $200 \times 200$, 
	and the corresponding parameter domain is $\left[0,1\right]^{2}$. 
This resolution ensures that the structural topology information derived 
	from the discrete density field remains consistent with 
	the structure defined by the continuous density field.
	
\subsection{Short beam example}
\begin{figure}[!htb]
	\centering
	\setcounter{subfigure}{0}
	\subfigure[]{
		\label{fig:model_beam}
		\includegraphics[width=0.3\linewidth]{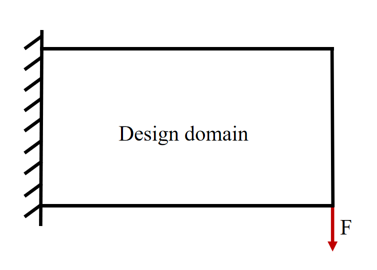}
	}
	\subfigure[]{
		\label{fig:model_beam_no_opt}
		\includegraphics[width=0.32\linewidth]{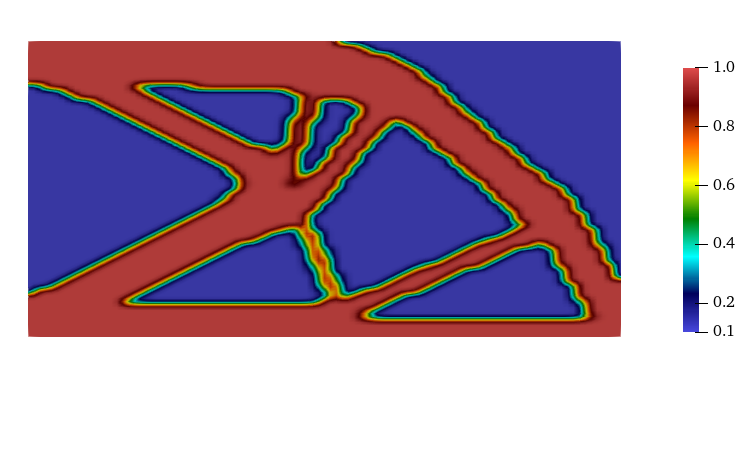}
	}
	\subfigure[]{
		\label{fig:model_beam_no_opt}
		\includegraphics[width=0.3\linewidth]{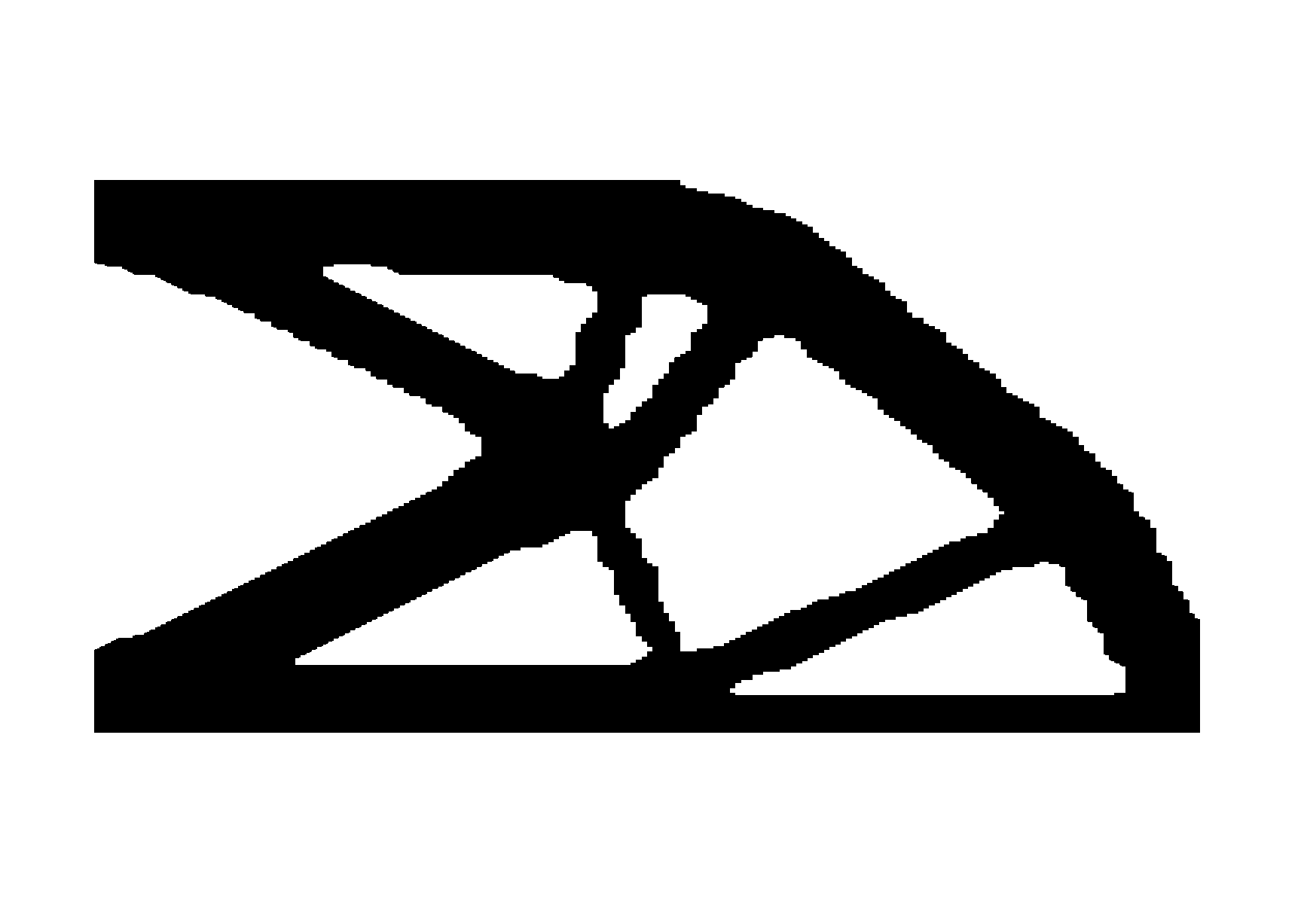}
	}
	\caption{(a) The short beam model. (b) Optimization result without topology control. 
	(c) Binary structure corresponding to the density field in (b).}
	\label{fig:short_beam_no_opt}
\end{figure}

In this example, the classical short beam problem shown in 
	Fig. \ref{fig:short_beam_no_opt} is used to investigate the capability of the 
	proposed method for topological control.
The length of the beam is set to two times its width. 
The design domain is represented by a $61\times 61$ bicubic NURBS surface.
The left boundary of the design domain is fixed, and a downward force of $F=10^{5}$ 
	is applied at its lower-right corner. 
The classical minimum compliance problem 
	without topology constraint is solved with $V_{0}=0.5$.
As shown in Fig. \ref{fig:model_beam_no_opt}, there exist $5$ holes in 
	the structure.

\begin{figure}[!htb]
	\centering
	\setcounter{subfigure}{0}
	\subfigure[$\bar{N}_{1}=1$]{
		\label{fig:short_1}
		\includegraphics[width=0.3\linewidth]{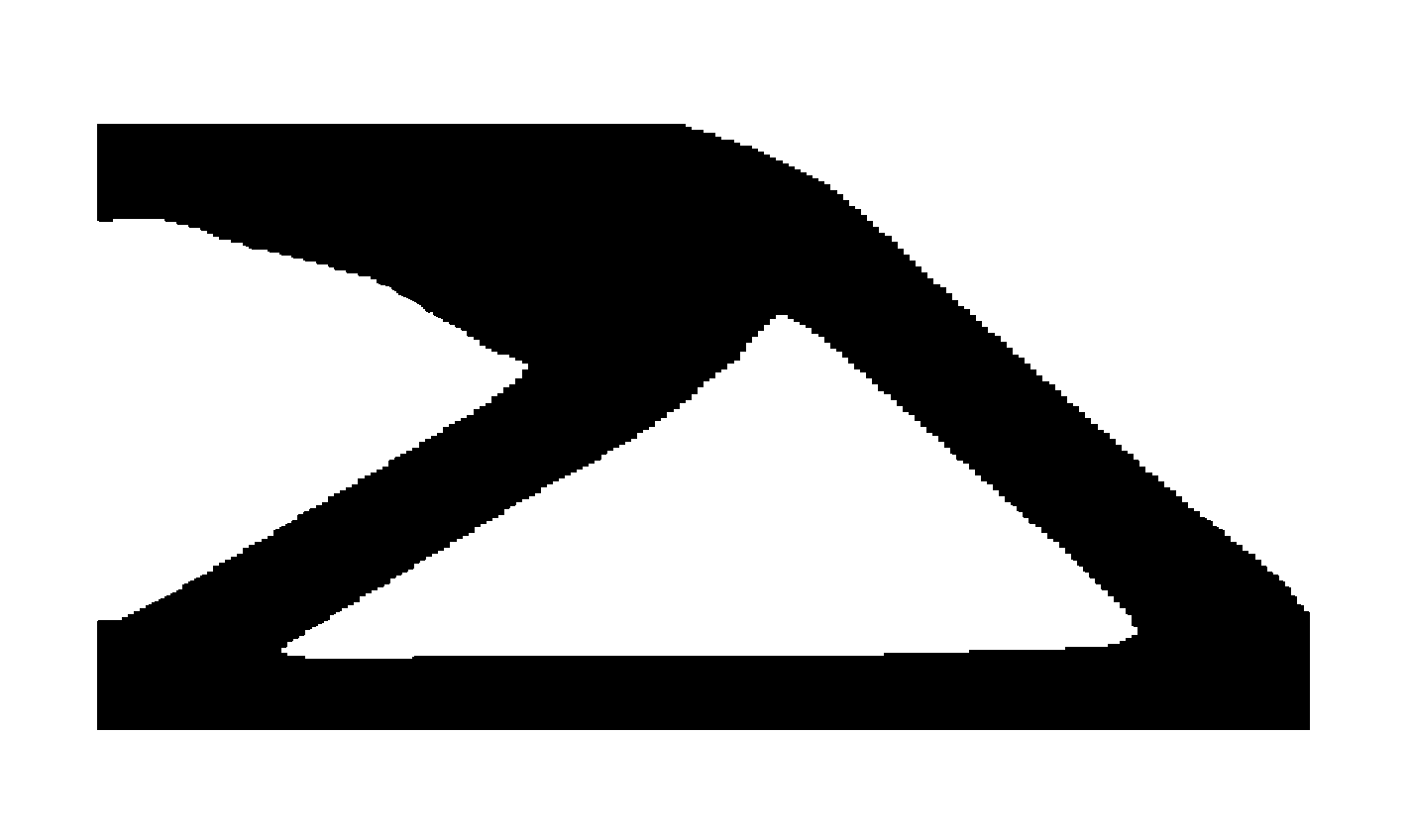}
	}
	\subfigure[$\bar{N}_{1}=2$]{
		\label{fig:short_2}
		\includegraphics[width=0.3\linewidth]{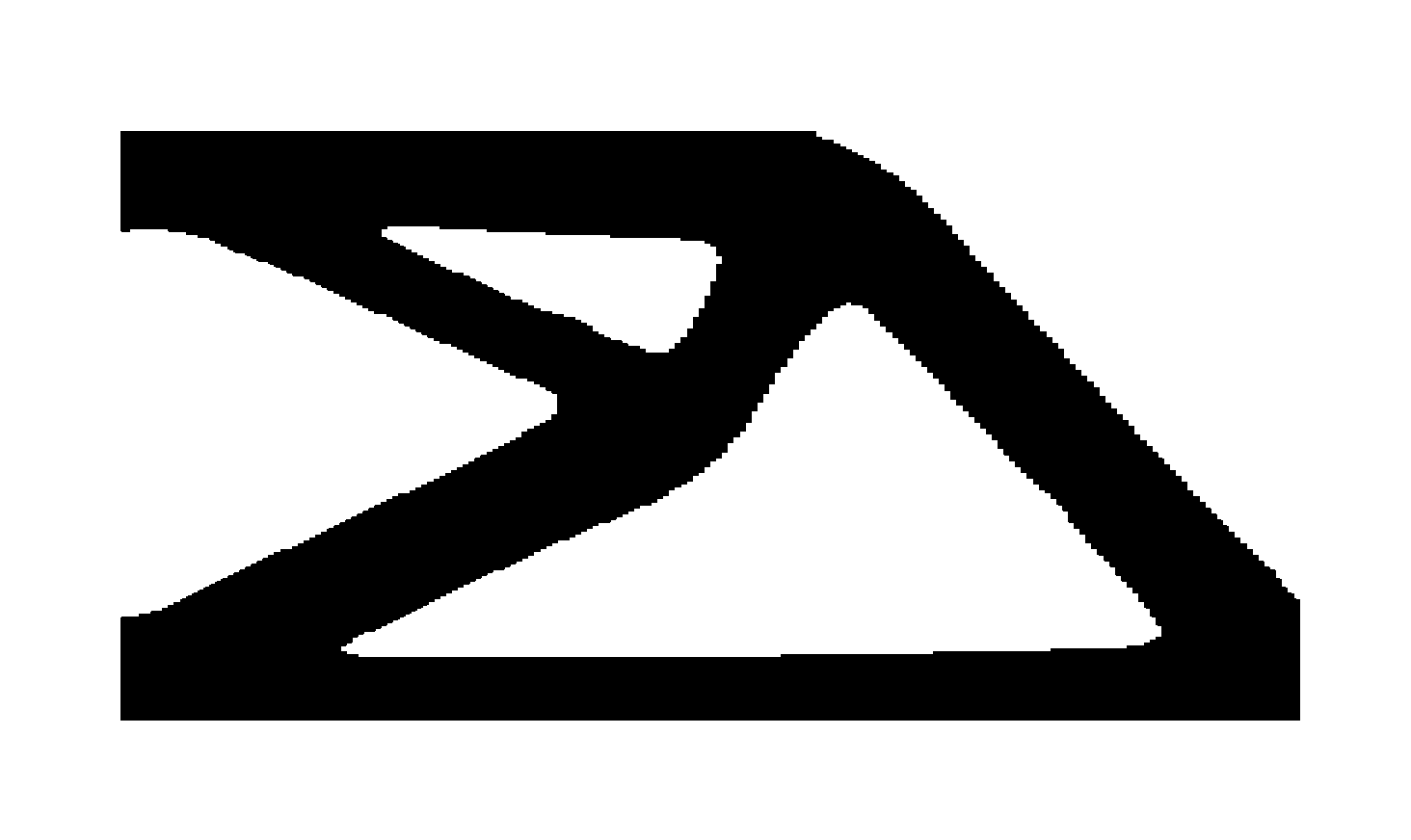}
	}
	\subfigure[$\bar{N}_{1}=3$]{
		\label{fig:short_3}
		\includegraphics[width=0.3\linewidth]{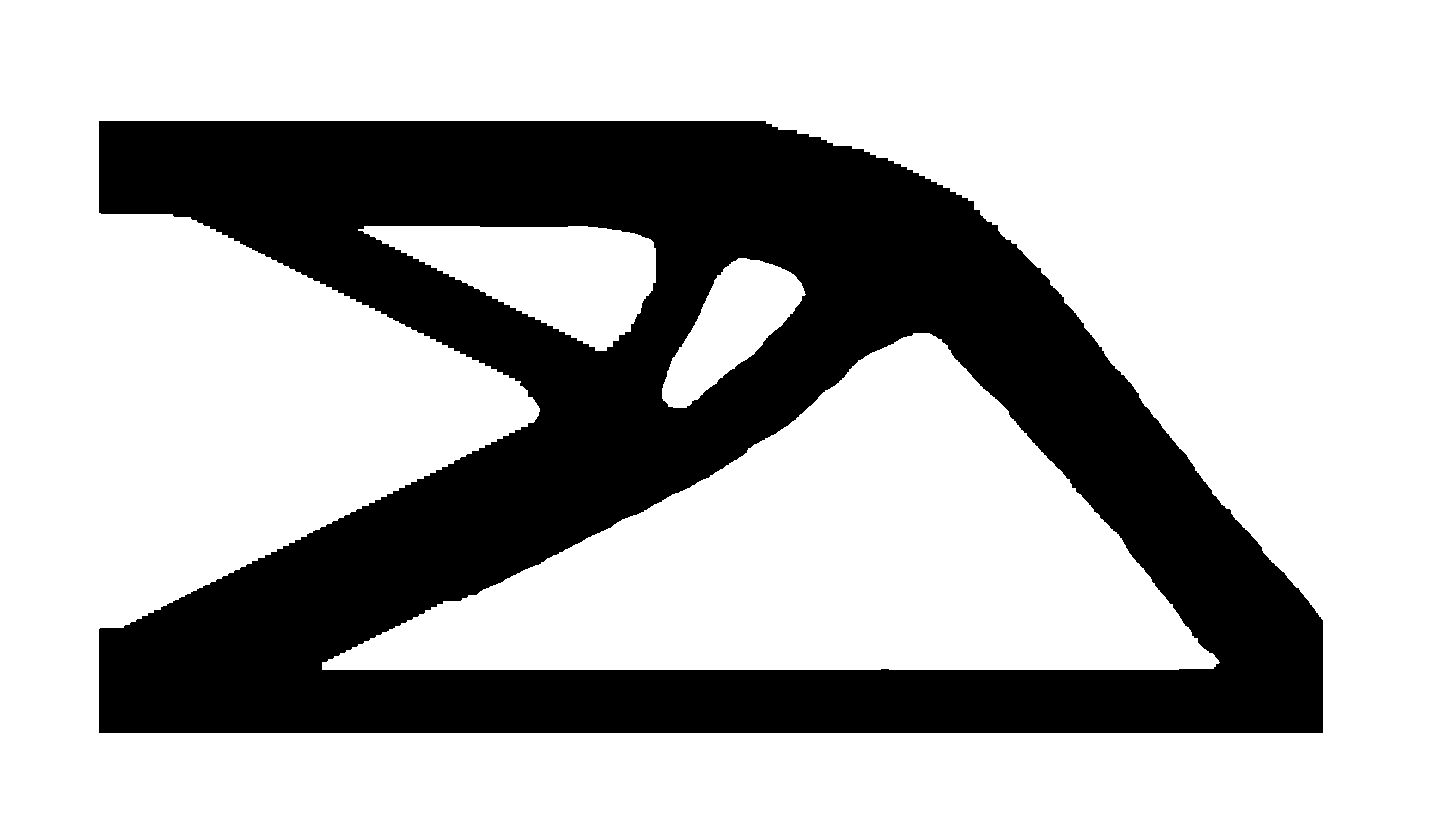}
	}\\
	\subfigure[$\bar{N}_{1}=4$]{
		\label{fig:short_4}
		\includegraphics[width=0.3\linewidth]{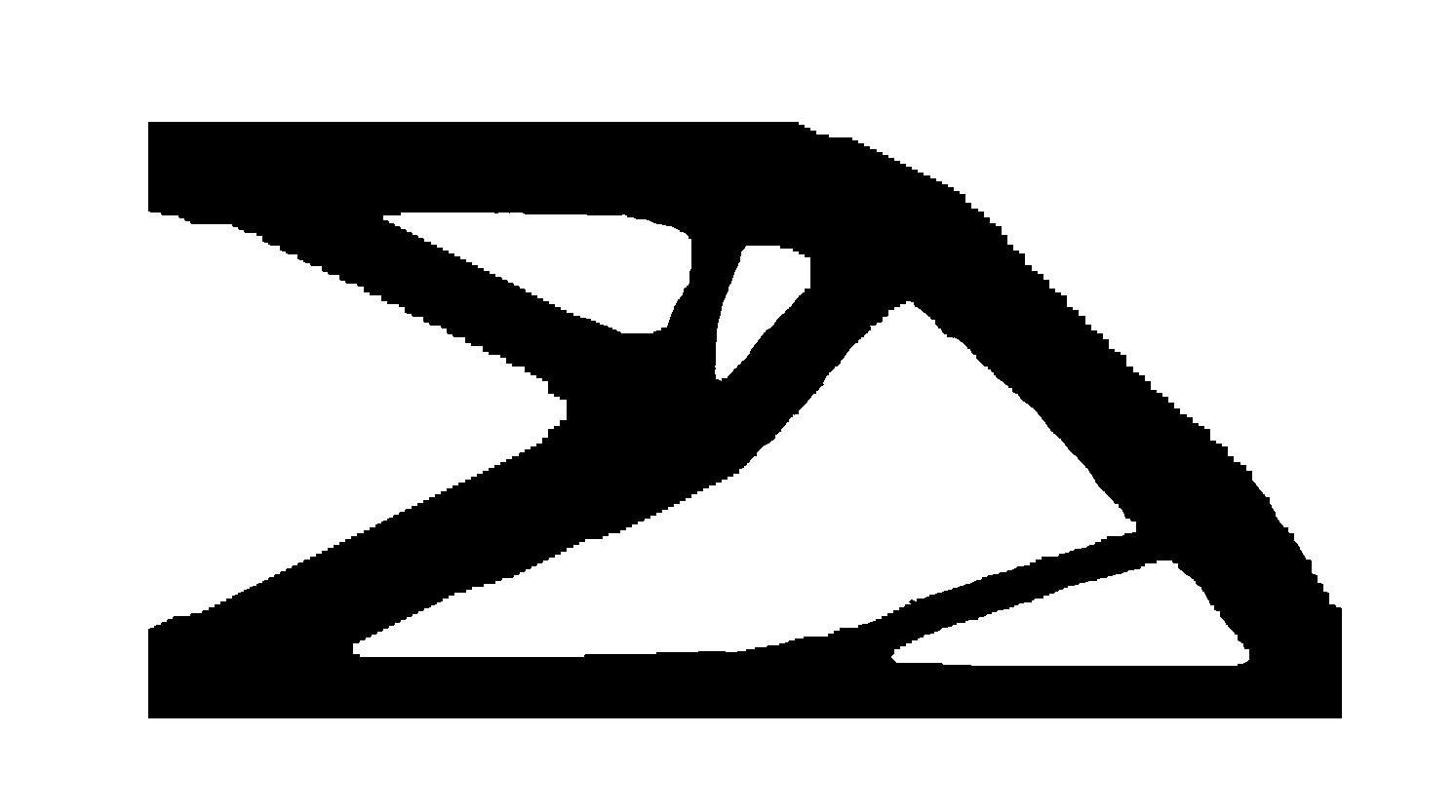}
	}
	\subfigure[$\bar{N}_{1}=5$]{
		\label{fig:short_5}
		\includegraphics[width=0.3\linewidth]{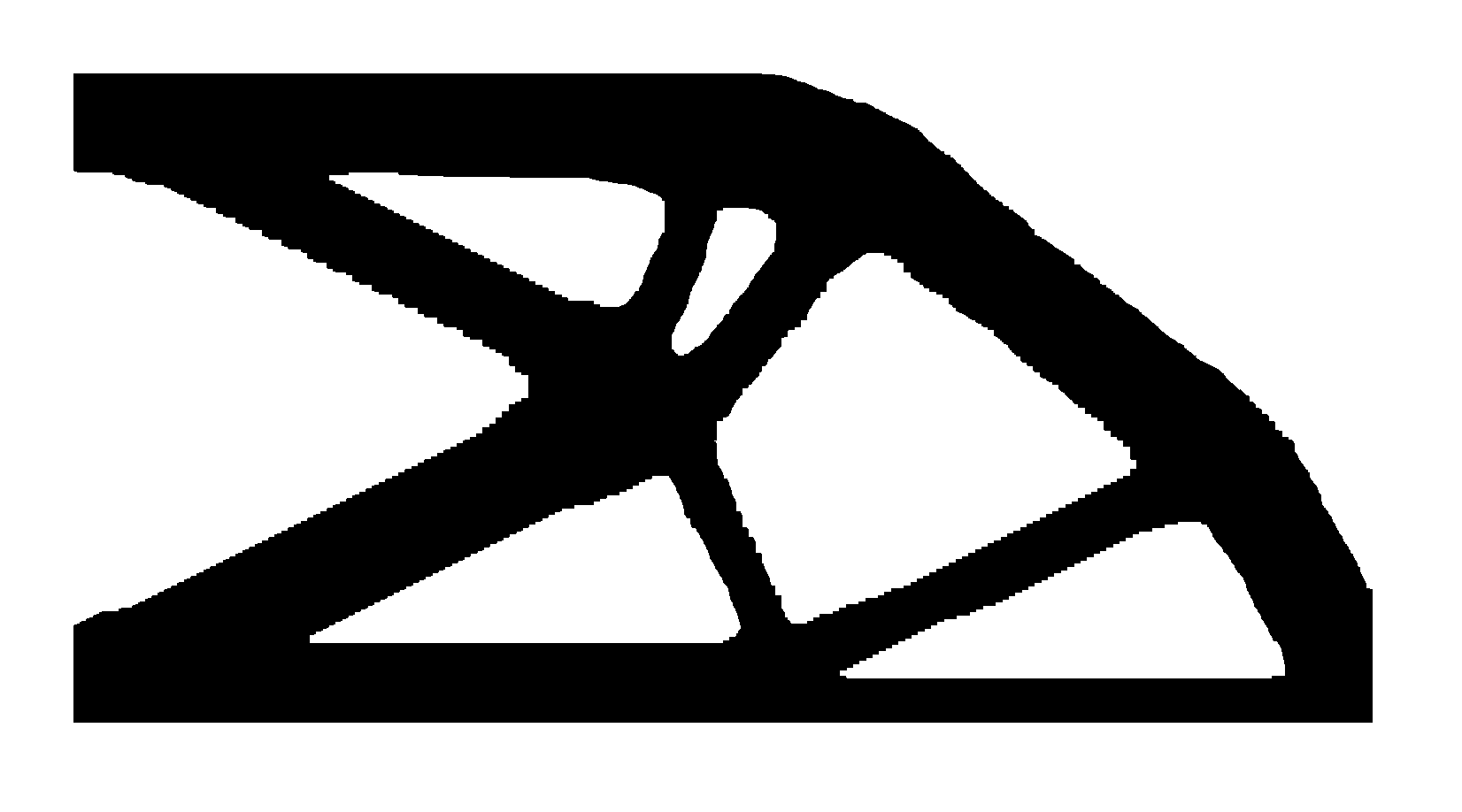}
	}
	\subfigure[$\bar{N}_{1}=6$]{
		\label{fig:short_6}
		\includegraphics[width=0.3\linewidth]{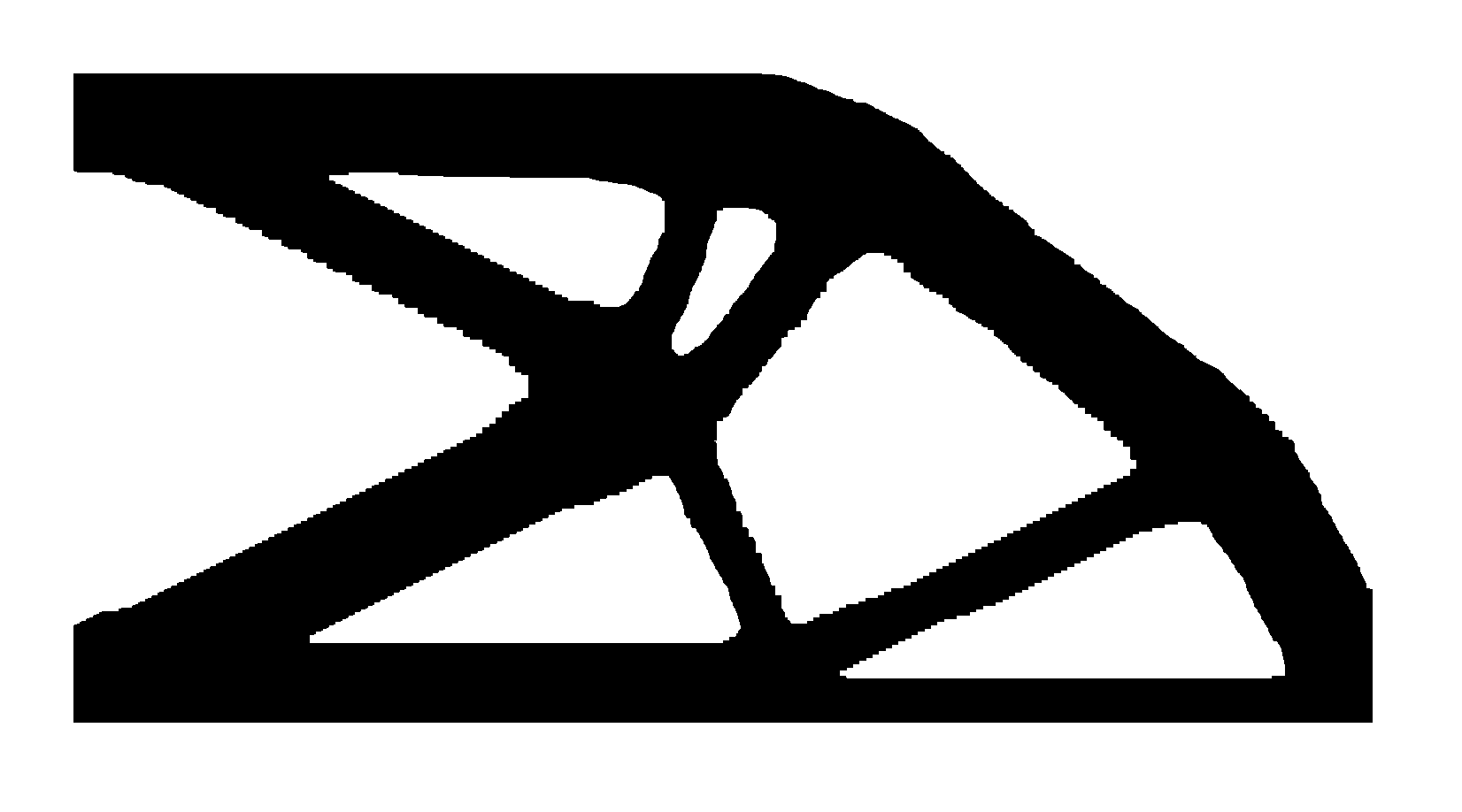}
	}	
	\caption{Optimization results of the short beam for $\bar{N}_{1}=1,2,3,4,5,6$, respectively.}
	\label{fig:short_beam_opt}
\end{figure}

A topological constraint is imposed such that the number of holes in the optimized structure 
	is no more than $\bar{N}_{1}$, 
	while guarantee its connectivity at the same time. 
Fig. \ref{fig:short_beam_opt} illustrates the optimized results for 
	$\bar{N}_{1} = 1, 2, 3, 4, 5,$ and $6$. 
It can be observed that our method successfully achieves effective control 
	of the structural topology. 
Table \ref{tbl:obj_short} summarizes the compliance values, 
	topology-control weights, 
	and computational time for different prescribed maximum number of holes $\bar{N}_{1}$.
It can be observed that the compliance increases gradually 
	as the allowable number of holes decreases. 
This trend is physically intuitive: restricting the number of holes reduces structural design freedom 
	and limits the ability of the optimizer to distribute material in a compliance-optimal manner.
Nevertheless, the compliance variation remains moderate, 
	indicating that the proposed topology-control strategy achieves 
	structural regularization with only a limited performance penalty.

Moreover, stricter topological constraints require larger weights 
	for the topology objective function.
When the target topology deviates significantly from the unconstrained solution, 
	a stronger topological driving force is necessary to overcome the 
	compliance-minimizing tendency of the optimization process.
This observation confirms that the proposed persistence-based objective 
	provides a controllable and tunable mechanism 
	for balancing structural performance and topological complexity.	

In addition, the optimized structures demonstrate smooth geometric transitions 
	without artificial discontinuities.
The structural connectivity is consistently preserved across all cases, 
	confirming the effectiveness of the zero-dimensional topology control.
When $\bar{N}_{1} = 5$ or $6$, the optimized structures are nearly 
	identical to the unconstrained design, 
	further validating that the proposed method does not introduce 
	unnecessary distortions when the prescribed topology is already satisfied.
These results demonstrate that the proposed framework enables explicit 
	and reliable regulation of structural topology 
	while maintaining competitive mechanical 
	performance and numerical stability.


\begin{table}[!htb]
	\caption{Comparison of compliance, number of holes, and time costs 
for the short beam benchmark under different topology-control settings.}
	\begin{tabular}{lllllll}
		\toprule
		        &$\bar{N}_{1}=1$ & $\bar{N}_{1}=2$  & $\bar{N}_{1}=3$ & $\bar{N}_{1}=4$ & $\bar{N}_{1}=5$ & $\bar{N}_{1}=6$\\
		\midrule
		weights $(\mu_{0},\mu_{1})$    &(0.8,0.6)  &(0.1,0.1)  &(0.8,0.6)  &(0.1,0.1)  &(0.1,0.1)   & (0.1,0.1)   \\
		Compliance  &5.01  &4.95  &4.91  &4.91  & 4.90  &4.89    \\
		Time (s)  &445.3  &440.2  &423.2  &422.1  & 418.6  &412.3    \\
		\bottomrule
	\end{tabular}
	\label{tbl:obj_short}
\end{table}

\begin{figure}[!htb]
	\centering
	\setcounter{subfigure}{0}
	\subfigure[Compliance vs. iteration]{
	\label{fig:ex1_iter}
		\includegraphics[width=0.3\linewidth]{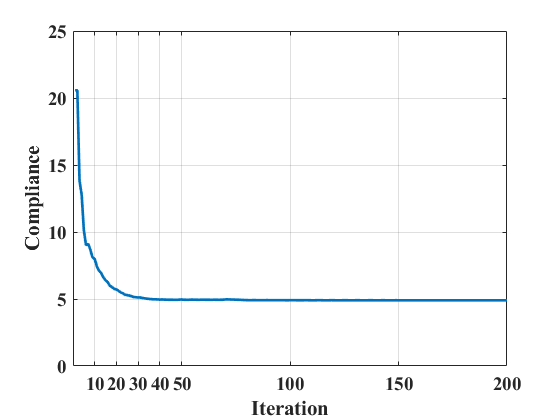}
	}
	\subfigure[iter=9: $N_{1}=2$]{
		\label{fig:ex1_iter_9}
		\includegraphics[width=0.31\linewidth]{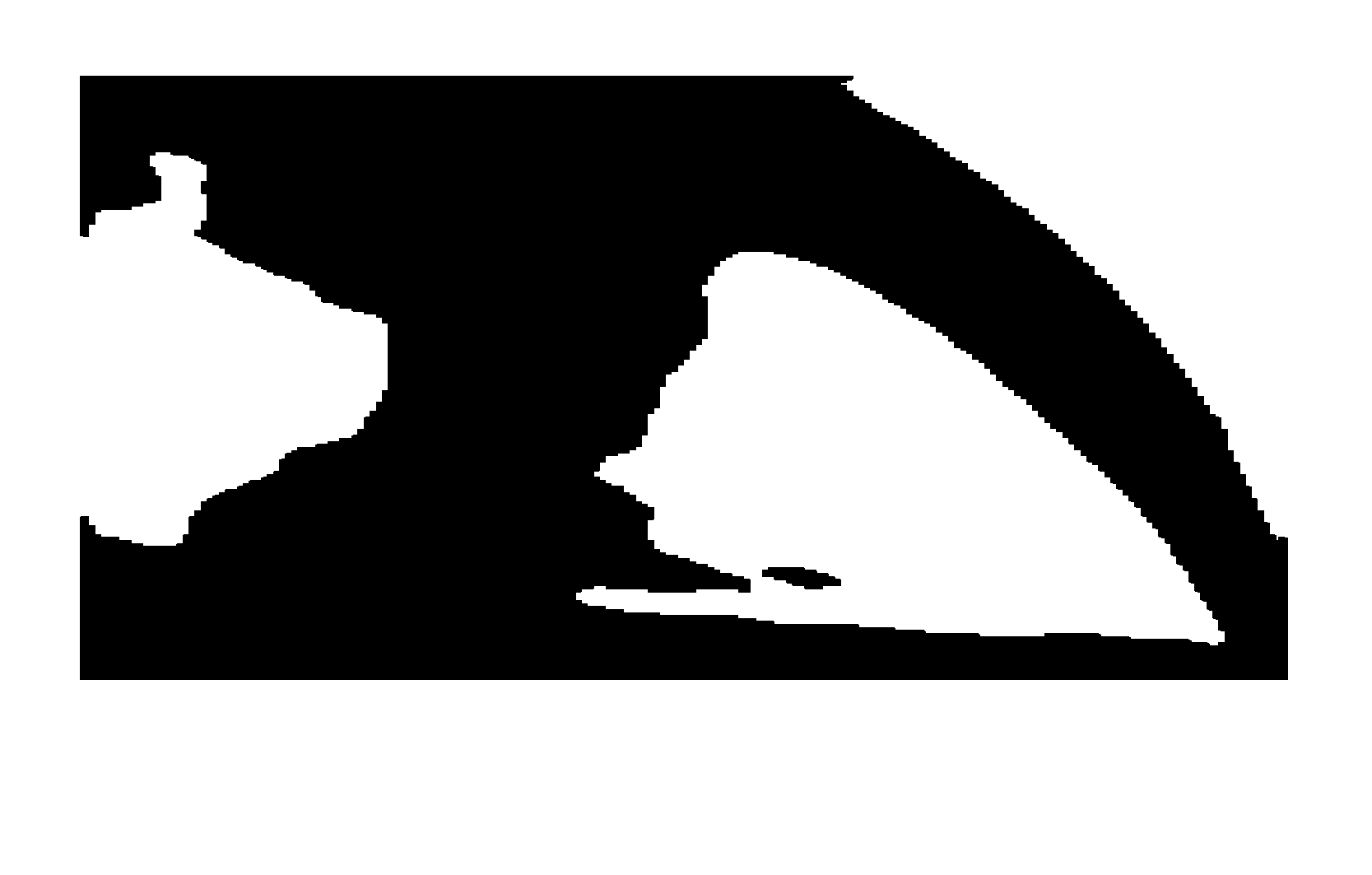}
	}
	\subfigure[iter=17: $N_{1}=7$]{
		\label{fig:ex1_iter_17}
		\includegraphics[width=0.31\linewidth]{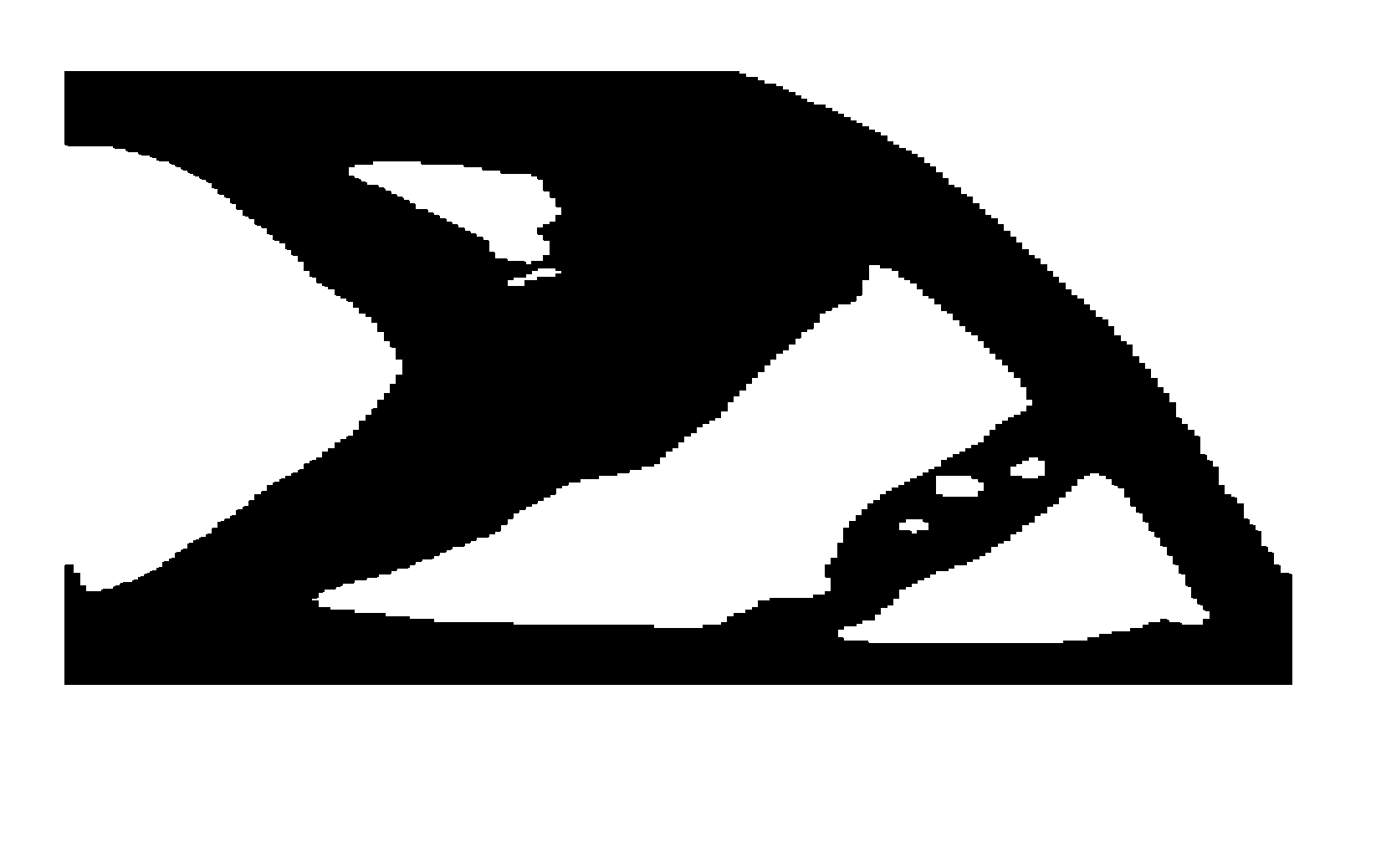}
	}\\
	\subfigure[iter=20: $N_{1}=4$]{
		\label{fig:ex1_iter_20}
		\includegraphics[width=0.31\linewidth]{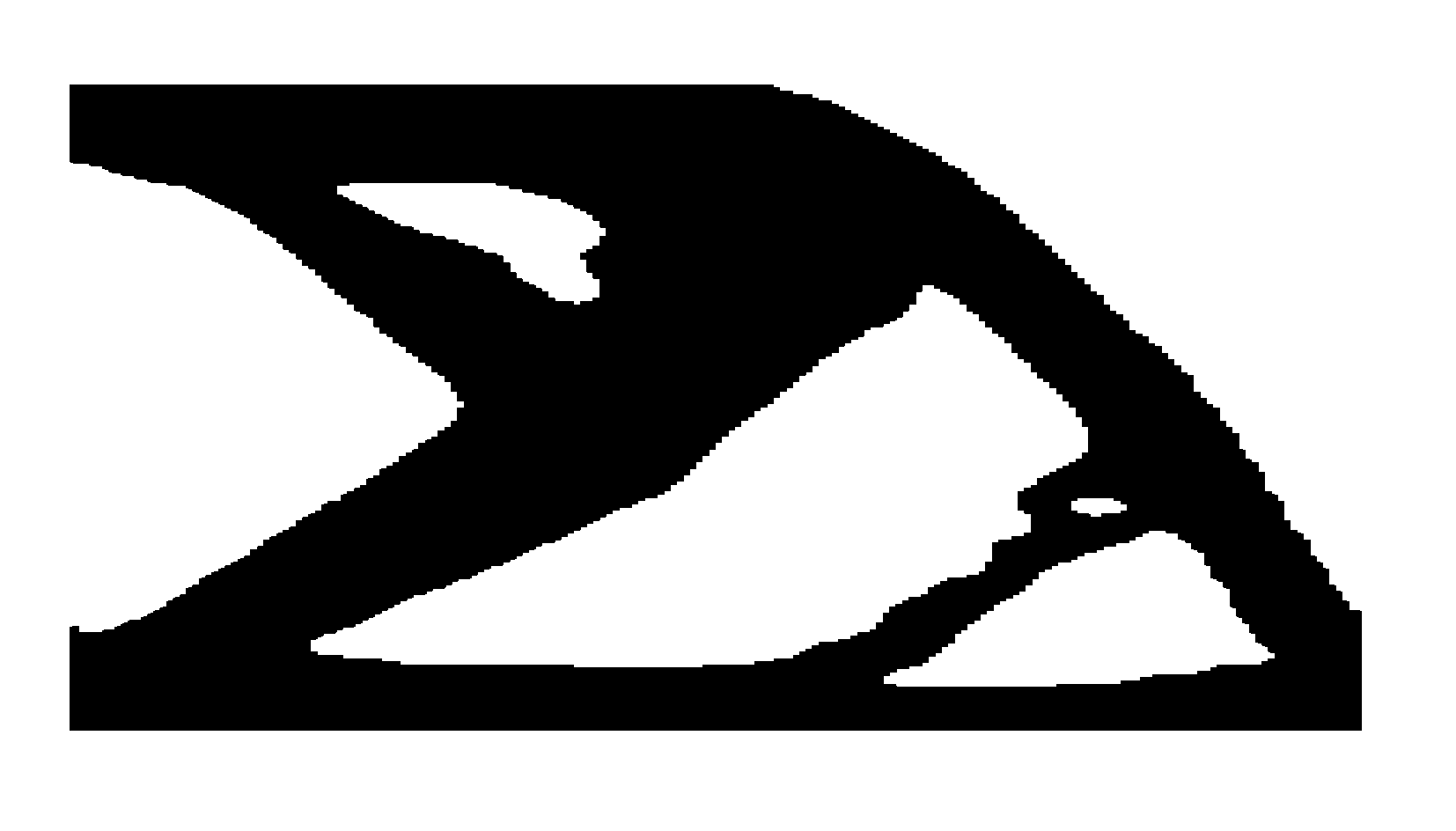}
	}
	\subfigure[iter=30: $N_{1}=3$]{
		\label{fig:ex1_iter_30}
		\includegraphics[width=0.31\linewidth]{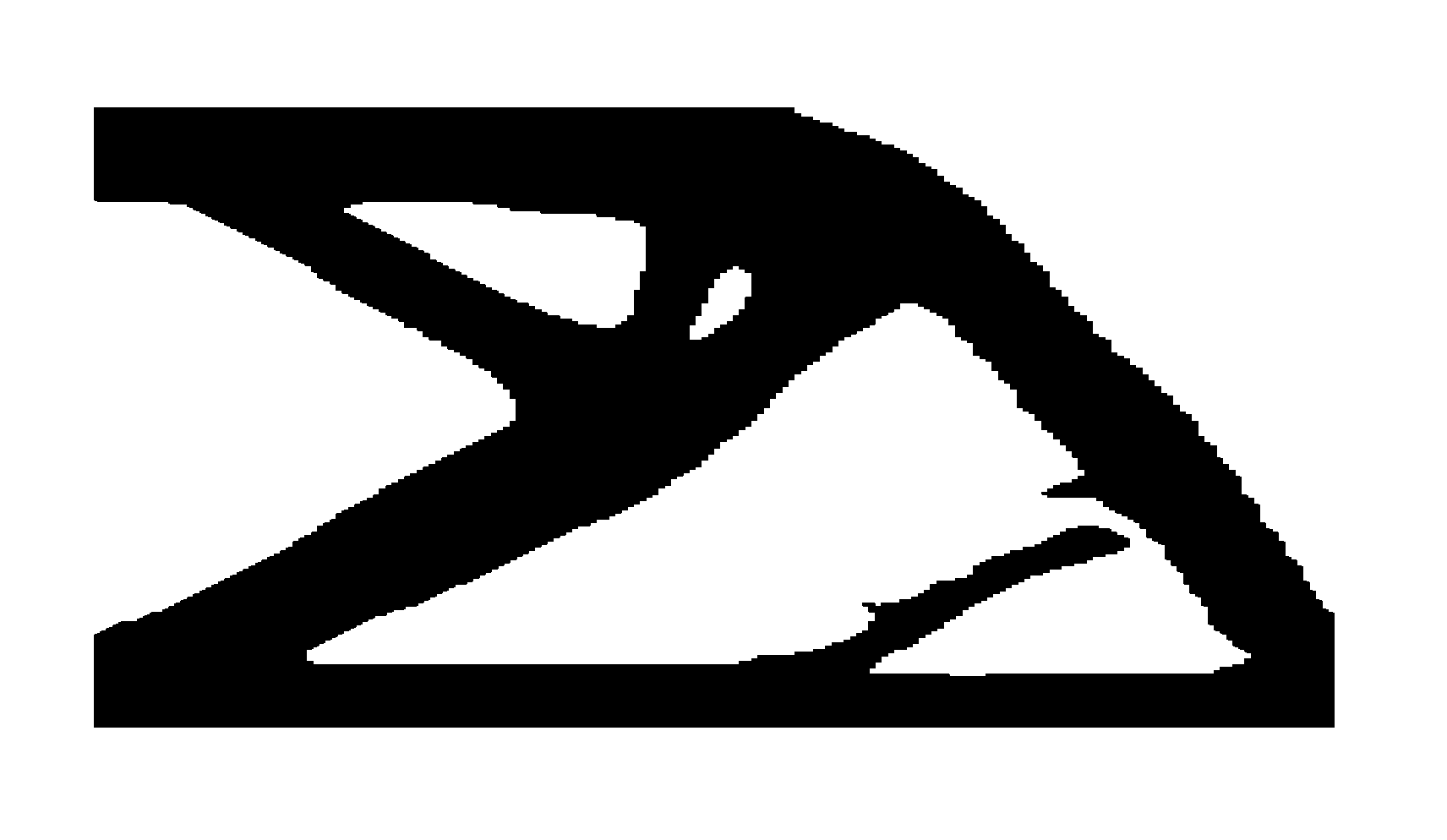}
	}
	\subfigure[iter=75: $N_{1}=3$]{
		\label{fig:ex1_iter_75}
		\includegraphics[width=0.31\linewidth]{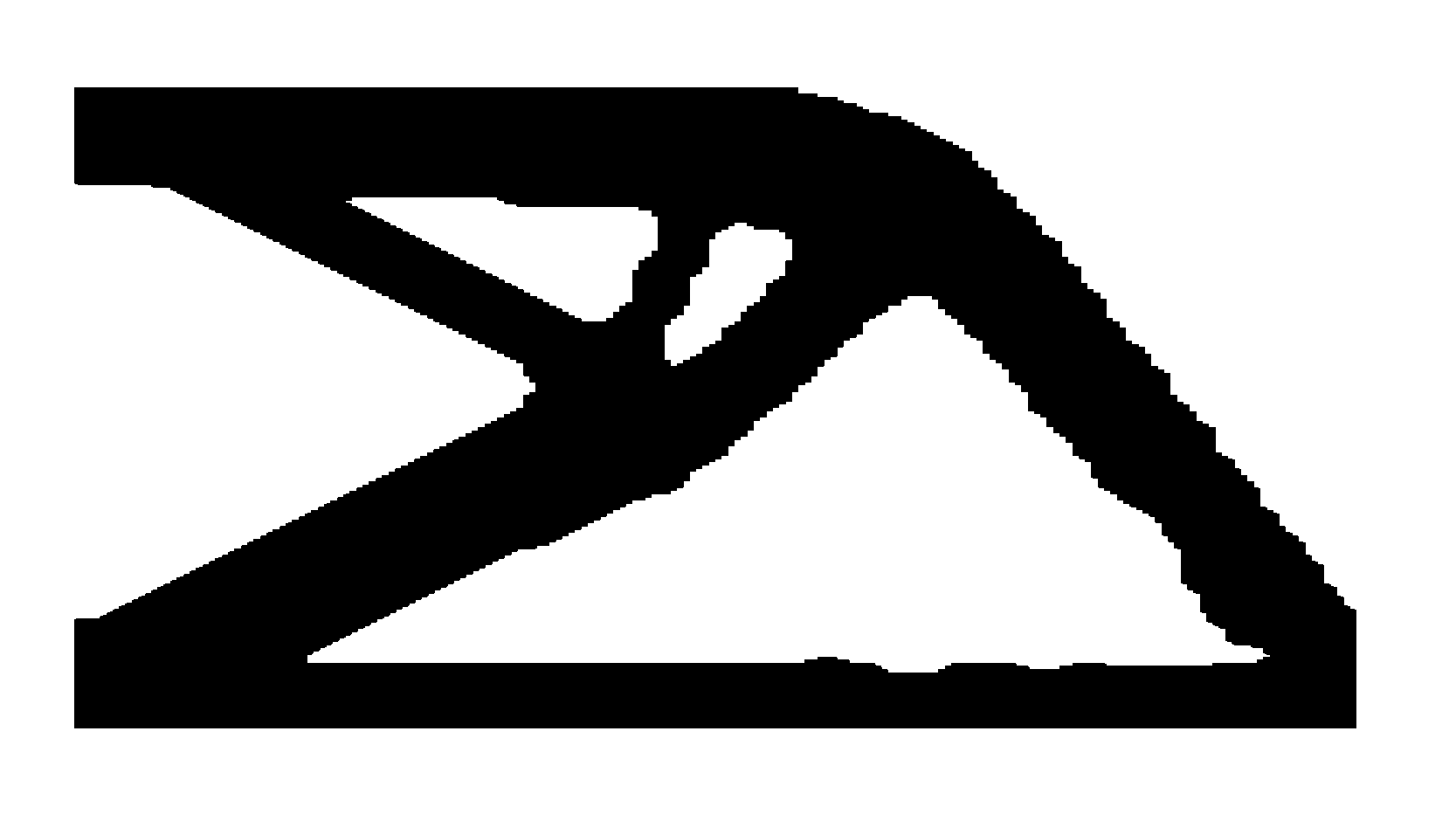}
	}
	\caption{Evolution of the optimized topology for the short beam when $\bar{N}_{1}=3$.
Binary density fields at representative iterations are shown.
The proposed topology-control strategy progressively eliminates redundant holes
and enforces the prescribed connectivity constraint.}
	\label{fig:short_process}
\end{figure}

Fig. \ref{fig:short_process} illustrates the evolution of the compliance 
	and the binarized density structures during the optimization process 
	under the constraint $\bar{N}_{1}=3$. 
As shown in Fig. \ref{fig:ex1_iter}, 
	the compliance decreases rapidly within the first $30$ 
	iterations and then gradually stabilizes. 
This behavior indicates that achieving a stable structure requires 
	approximately $30$ iterations without topological control.

Figs. \ref{fig:ex1_iter_9} - \ref{fig:ex1_iter_75} depict the structural evolution 
	and demonstrate the effectiveness of the proposed topological objective function.
Fig. \ref{fig:ex1_iter_17} shows that at the $17$th iteration, 
	the number of holes in the structure reaches 7, 
	exceeding the prescribed maximum and thereby violating the topological constraint. 
After activating topology control, the number of holes is reduced to 4 by the 20th iteration (Fig.~\ref{fig:ex1_iter_20}), 
	confirming that the persistence-based objective effectively suppresses redundant topological features.
Compared with the configuration at the 17th iteration, 
	several small holes within the structure are removed, 
	indicating that the optimization preferentially eliminates topologically insignificant holes.	
By the $30$th iteration, the number of holes further decreases to 3, 
	satisfying the prescribed constraint.
This reduction is achieved through structural reconnection and merging of adjacent void regions, 
	which alters the homological structure while preserving overall mechanical performance.
The structural transitions occur smoothly without introducing numerical oscillations 
	or instability in the compliance curve.
Finally, at the $75$th iteration, excess material is completely eliminated 
	and the structure reaches a stable topology.
The compliance remains nearly constant during the late stage of the optimization, 
	demonstrating that the proposed topology-control mechanism 
	does not compromise convergence stability.

\subsection{L-shaped beam example}
\begin{figure}[!htb]
	\centering
	\setcounter{subfigure}{0}
	\subfigure[]{
		\label{fig:ex2_model}
		\includegraphics[width=0.28\linewidth]{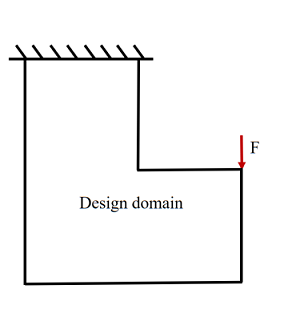}
	}
	\subfigure[]{
	\label{fig:ex2_density}
		\includegraphics[width=0.30\linewidth]{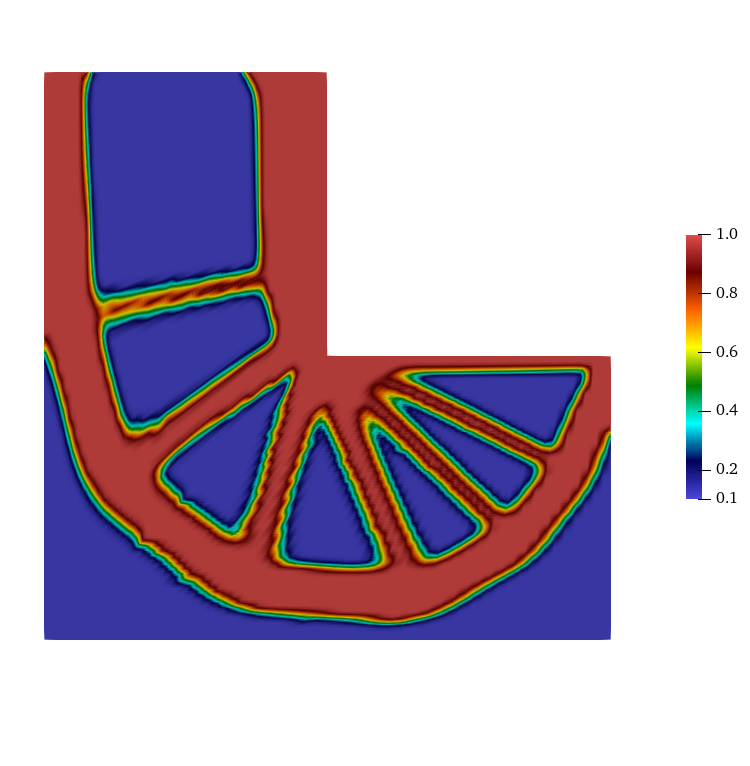}
	}
	\subfigure[]{
	\label{fig:ex2_binar}
		\includegraphics[width=0.32\linewidth]{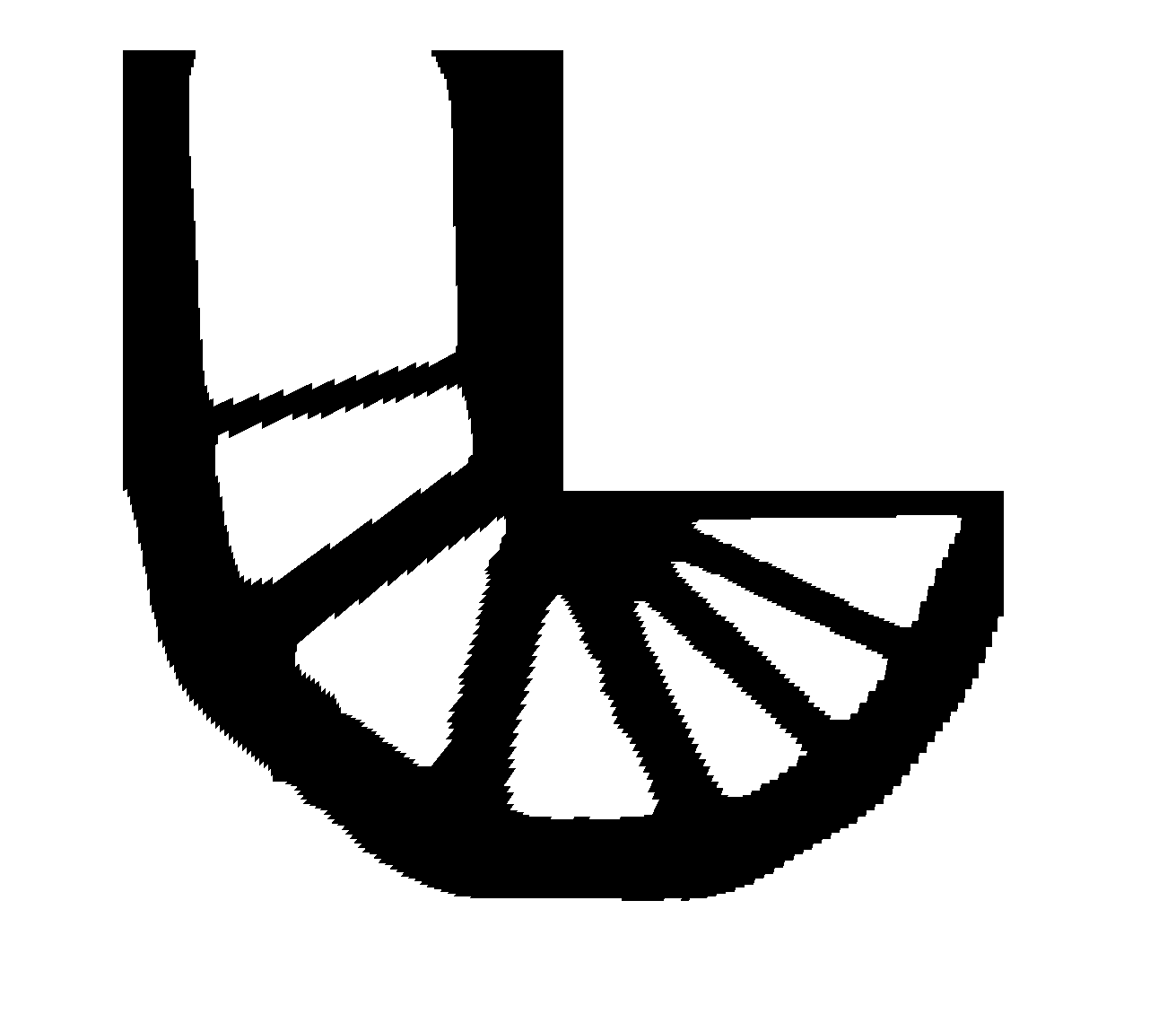}
	}
	\caption{(a) L-shaped beam model. (b) Optimization result without topology control. 
	(c) Binary structure corresponding to the density field in (b).}
	\label{fig:L_beam_no_opt}
\end{figure}

Fig. \ref{fig:L_beam_no_opt} presents the L-shaped beam benchmark.
The design domain is represented by a $103\times 52$ biquadratic NURBS surface.
The top boundary of the domain is fully fixed, and a downward point load of $F=10^{5}$ 
		is applied, as illustrated in Fig.~\ref{fig:ex2_model}. 
The classical minimum compliance problem 
	without topology constraint is solved with $V_{0}=0.5$.
As shown in Figs.~\ref{fig:ex2_density} and~\ref{fig:ex2_binar}, 
	there exist $6$ holes in the structure without 
	topological control.

\begin{figure}[!htb]
	\centering
	\setcounter{subfigure}{0}
	\subfigure[$\bar{N}_{1}=1$]{
	    \includegraphics[width=0.22\linewidth]{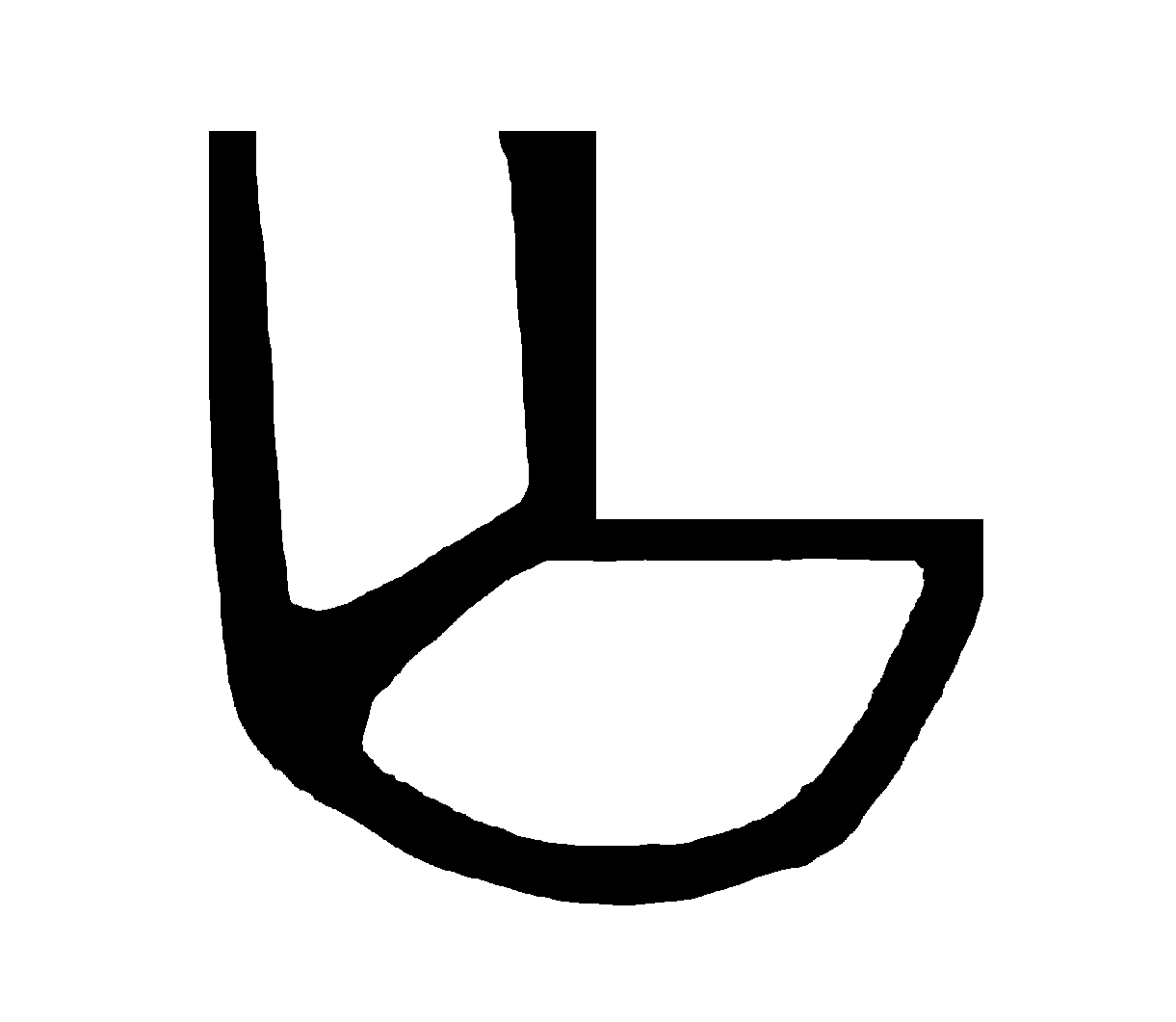}
	}
	\subfigure[$\bar{N}_{1}=2$]{
	     \includegraphics[width=0.22\linewidth]{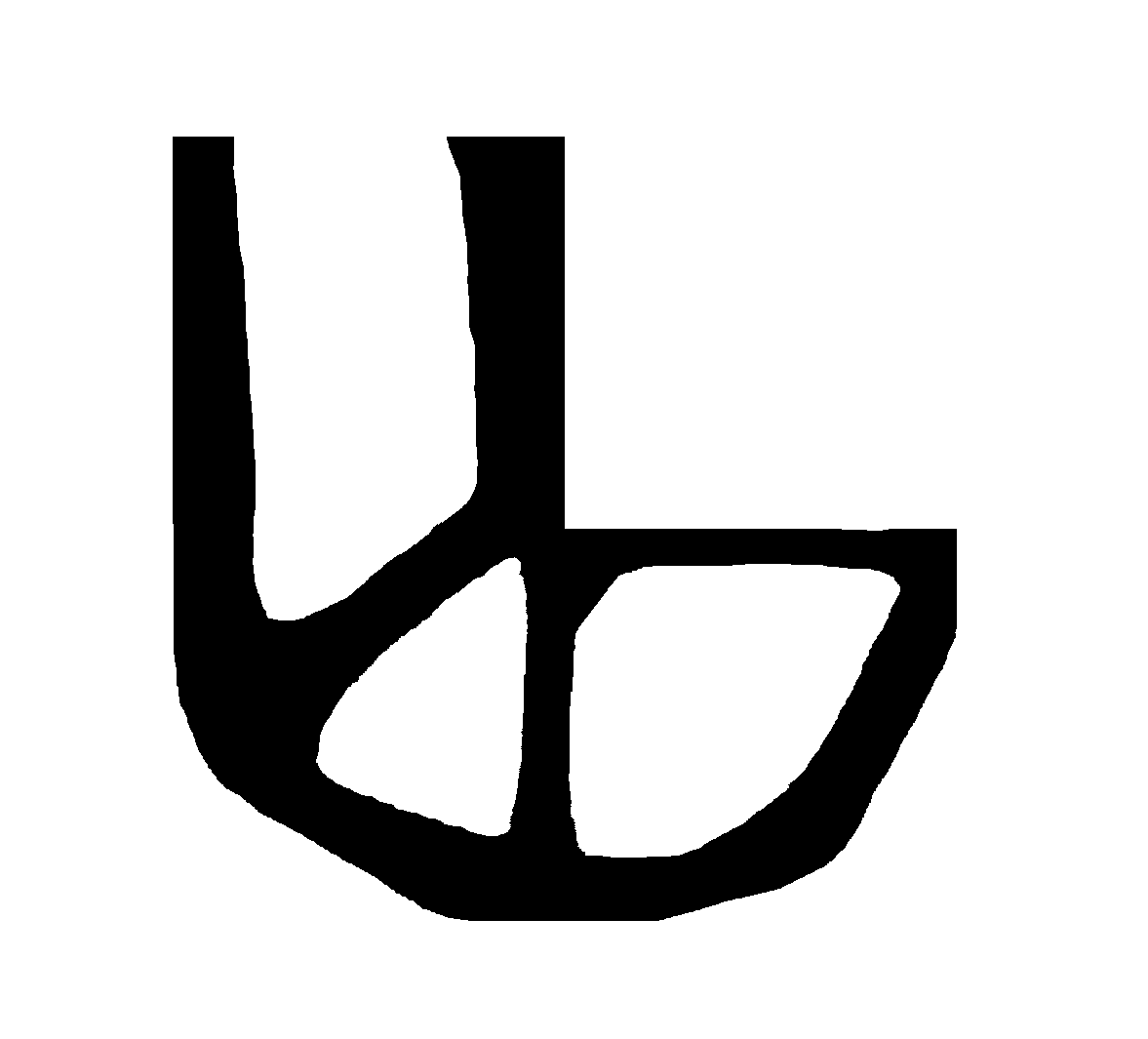}
	}
	\subfigure[$\bar{N}_{1}=3$]{
	     \includegraphics[width=0.24\linewidth]{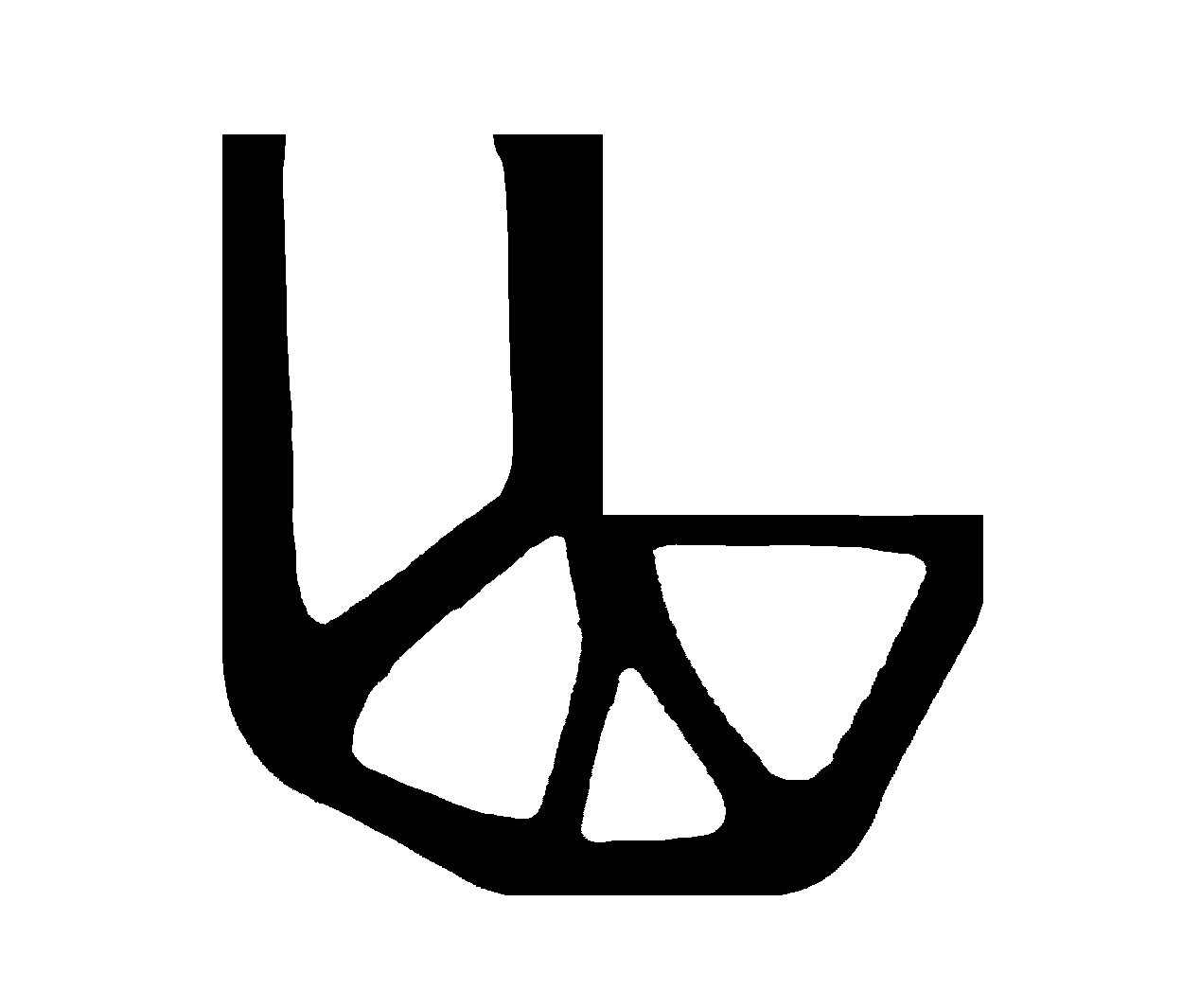}
	}
	\subfigure[$\bar{N}_{1}=4$]{
		\includegraphics[width=0.21\linewidth]{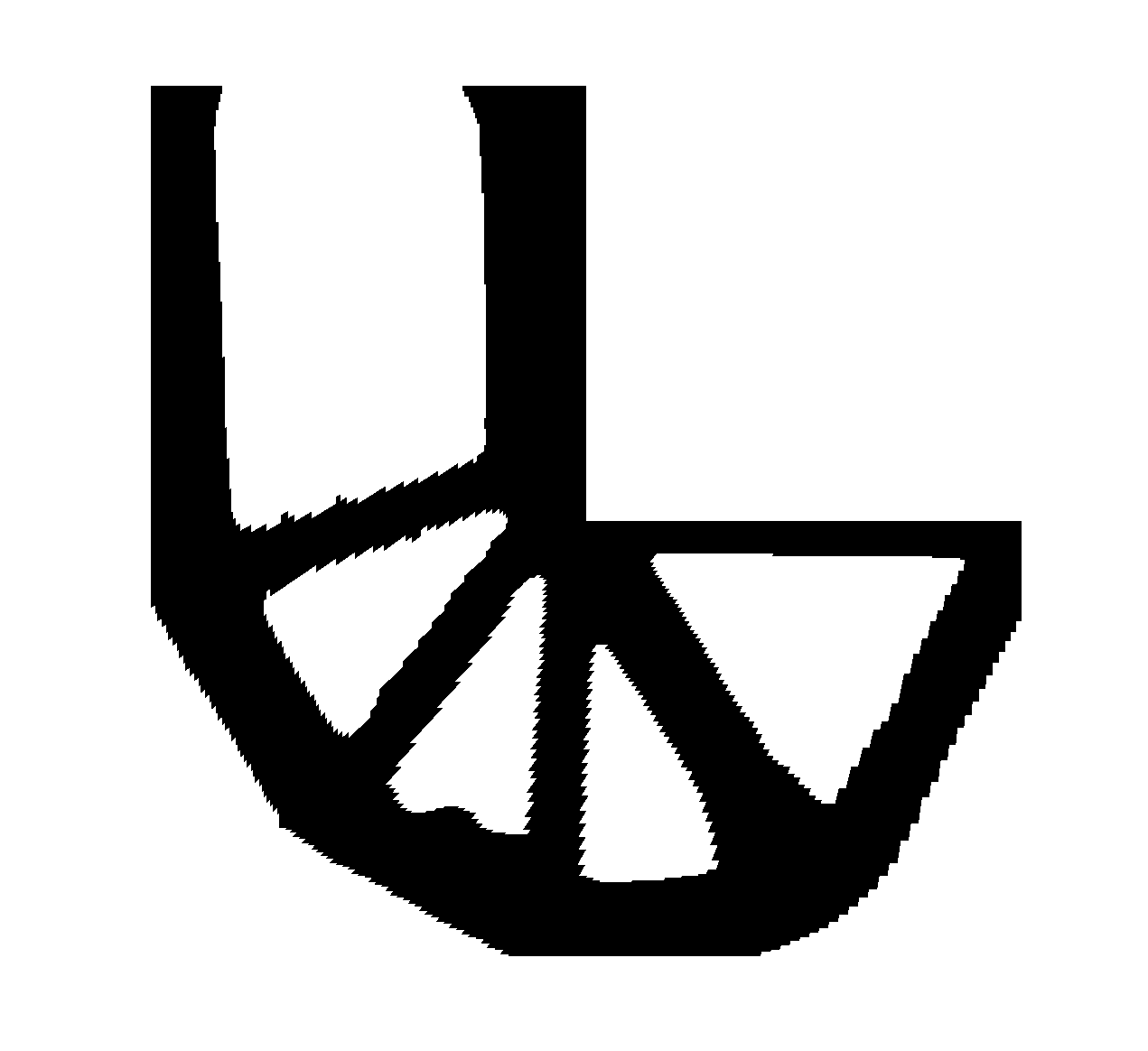}
	}\\
	\subfigure[$\bar{N}_{1}=5$]{
		\includegraphics[width=0.22\linewidth]{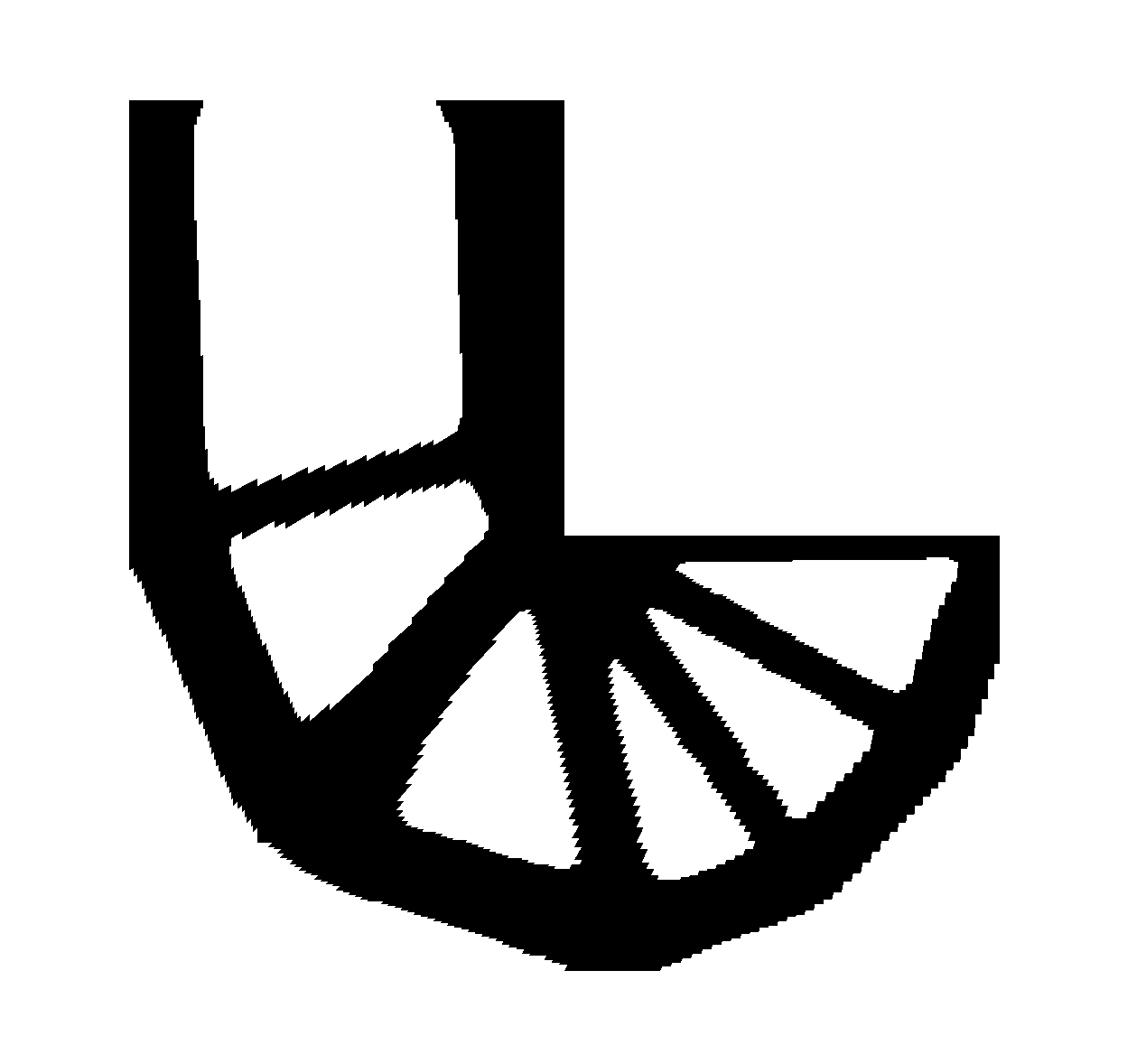}
	}
	\subfigure[$\bar{N}_{1}=6$]{
		\includegraphics[width=0.22\linewidth]{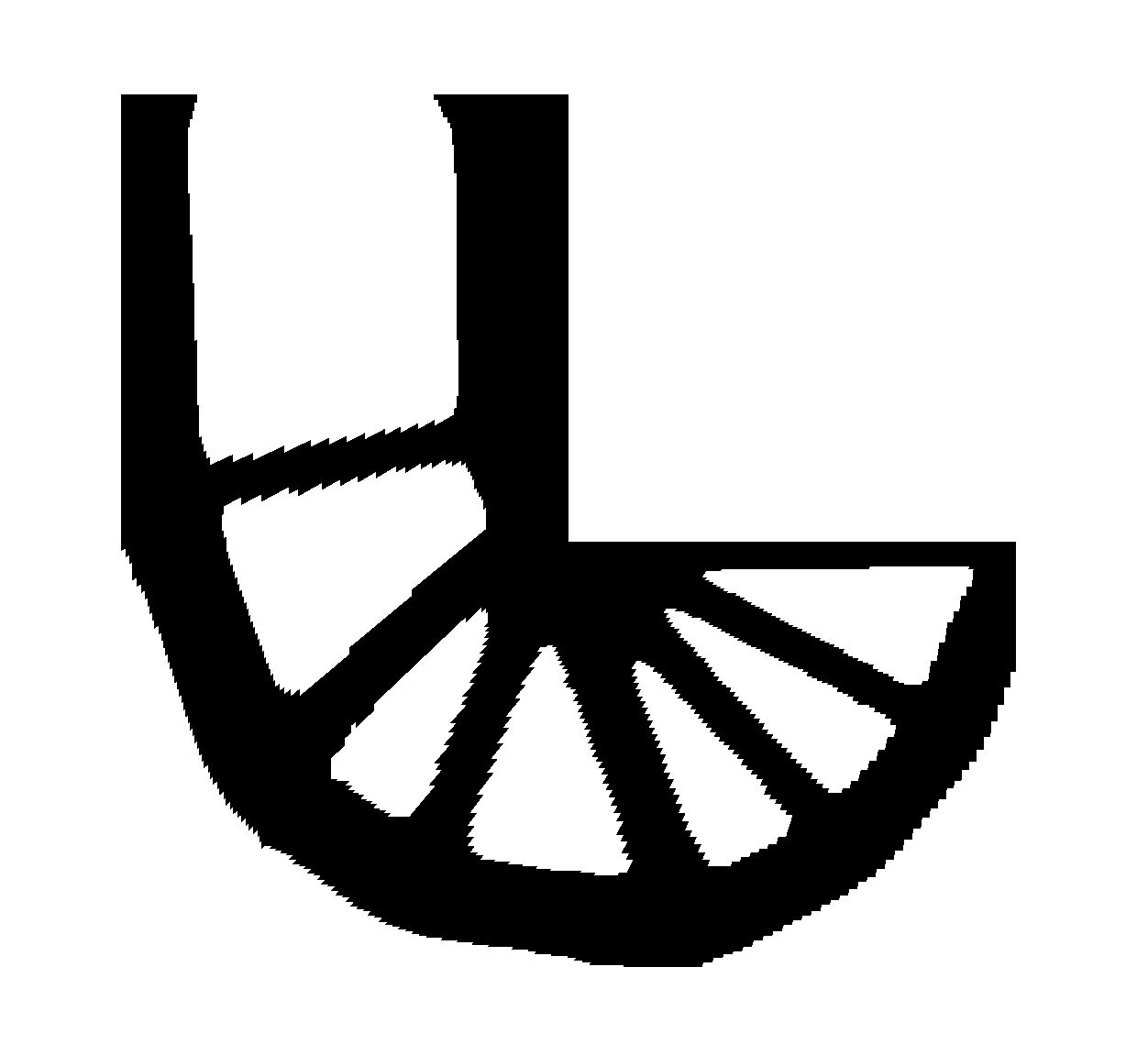}
	}
	\subfigure[$\bar{N}_{1}=7$]{
		\includegraphics[width=0.22\linewidth]{example_fig/L_beam/snapshot700.png}
	}
	\caption{Optimization results of the L-shaped beam for $\bar{N}_{1}=1,2,3,4,5,6,7$, respectively.}
	\label{fig:L_beam_opt}
\end{figure}

We impose topological constraints on the maximum number of holes 
	in the L-shaped beam and guarantee the structural connectivity.
Fig. \ref{fig:L_beam_opt} shows the optimization results for 
	$\bar{N}_{1} = 1, 2, 3, 4, 5, 6, 7$.
Table \ref{tbl:obj_l} reports the compliance values,
	topology-objective weights 
	and computational time for $\bar{N}_{1} = 1, 2, 3, 4, 5, 6, 7$. 
The results exhibit a clear performance–complexity trade-off: 
	as $\bar{N}{1}$ decreases, the compliance increases moderately, 
	reflecting the reduced design freedom induced 
	by stricter topological constraints.
Meanwhile, tighter hole limits generally require larger topology weights 
	$(\mu{0},\mu_{1})$, indicating that a stronger topological driving 
	force is needed when the prescribed topology deviates 
	substantially from the compliance-optimal unconstrained solution.
When $\bar{N}_{1}\ge 6$, the resulting structures are close to the unconstrained design, 
	confirming that the proposed method does not introduce unnecessary distortions 
	once the target topology is already satisfied.

\begin{table}[!htb]
	\centering
	\caption{Comparison of compliance, number of holes, and time costs 
for the L-shaped beam benchmark under different topology-control settings.}
	\begin{tabular}{llllllll}
		\toprule
		        &$\bar{N}_{1}=1$ & $\bar{N}_{1}=2$  & $\bar{N}_{1}=3$ &$\bar{N}_{1}=4$ & $\bar{N}_{1}=5$  & $\bar{N}_{1}=6$ & $\bar{N}_{1}=7$\\
		\midrule
		weights $(\mu_{0},\mu_{1})$    &(0.8,0.6)  &(0.8,0.6)  &(0.4,0.3) &(0.4,0.3)  &(0.4,0.3)  &(0.2,0.1)  &(0.1,0.1) \\
		Compliance  &5.46  &5.44  &5.42 &5.41 &5.41  &5.39  &5.38 \\
		Time (s)  &538.4 &527.9  &501.3 &503.4  &478.8  &467.8  &468.2 \\
		\bottomrule
	\end{tabular}
	\label{tbl:obj_l}
\end{table}

\begin{figure}[!htb]
	\centering
	\setcounter{subfigure}{0}
	\subfigure[Compliance vs. iteration]{
		\label{fig:ex2_iter}
		\includegraphics[width=0.33\linewidth]{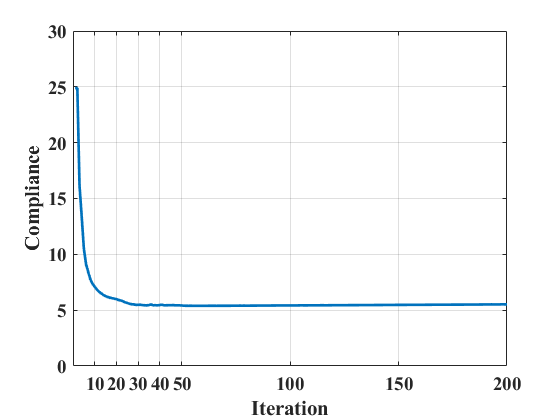}
	}
	\subfigure[iter=11: $N_{1}=3$]{
		\label{fig:ex2_iter_11}
		\includegraphics[width=0.29\linewidth]{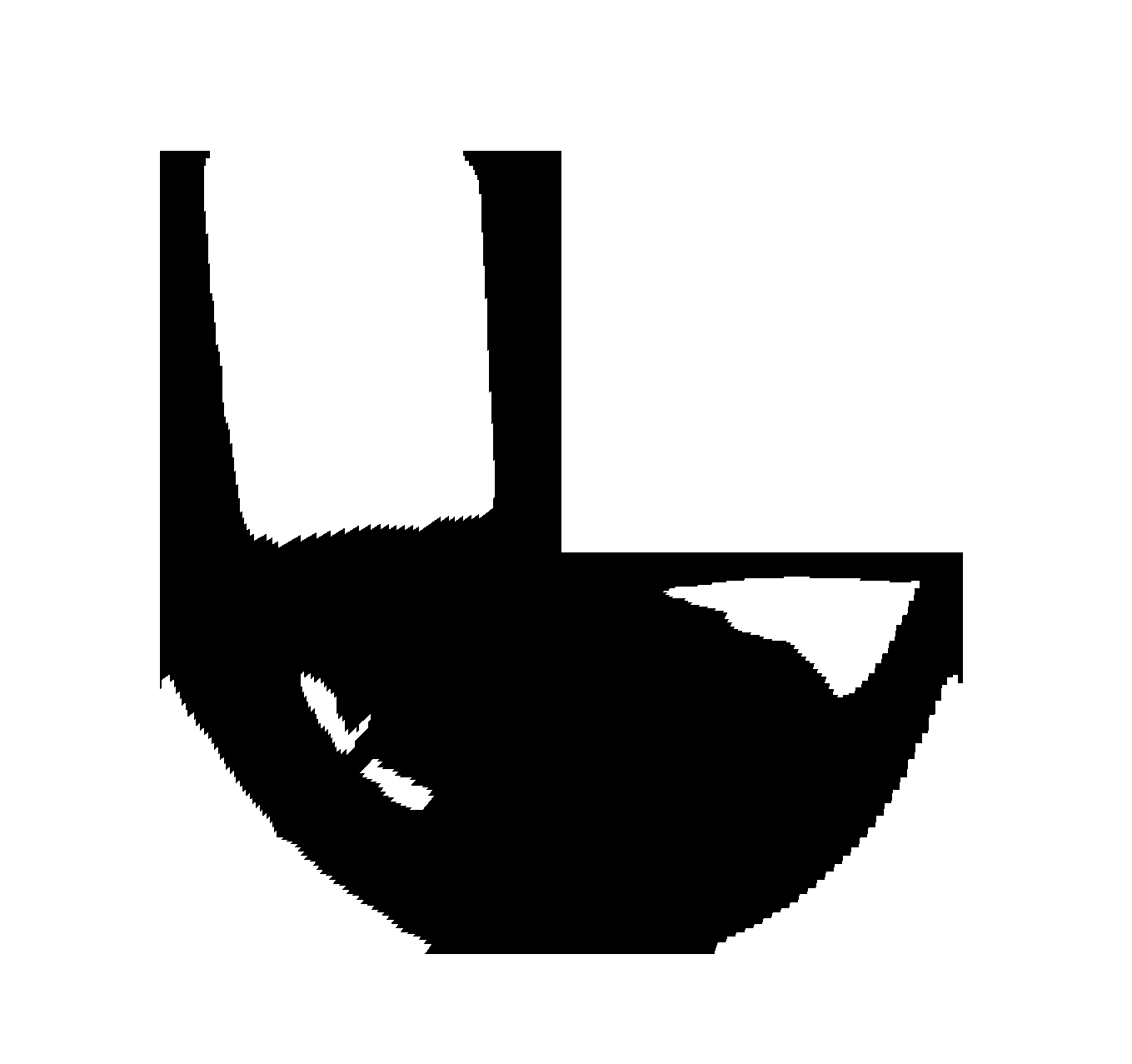}
	}
	\subfigure[iter=19: $N_{1}=9$]{
		\label{fig:ex2_iter_19}
		\includegraphics[width=0.29\linewidth]{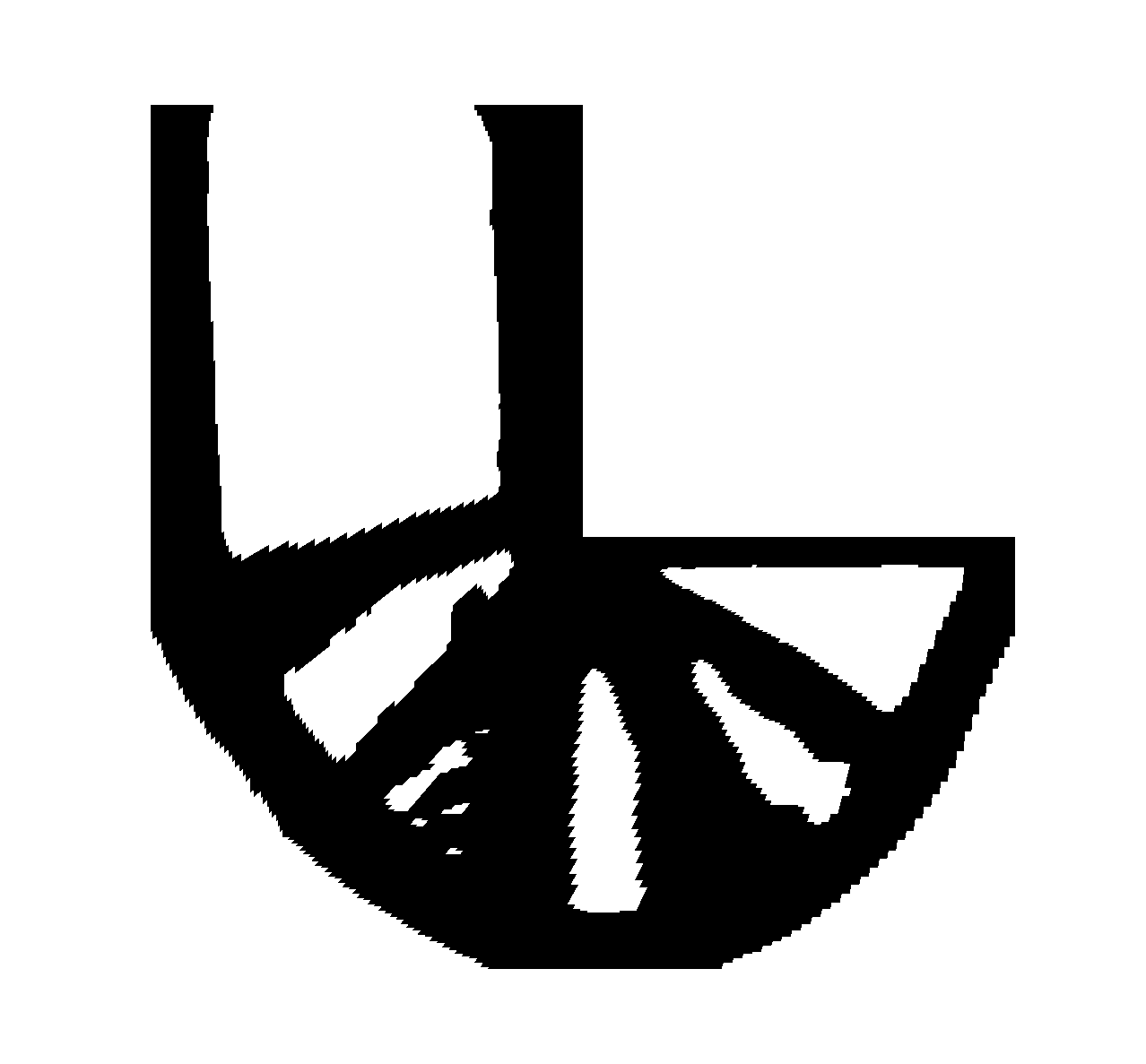}
	}\\
	\subfigure[iter=20: $N_{1}=6$]{
		\label{fig:ex2_iter_20}
		\includegraphics[width=0.29\linewidth]{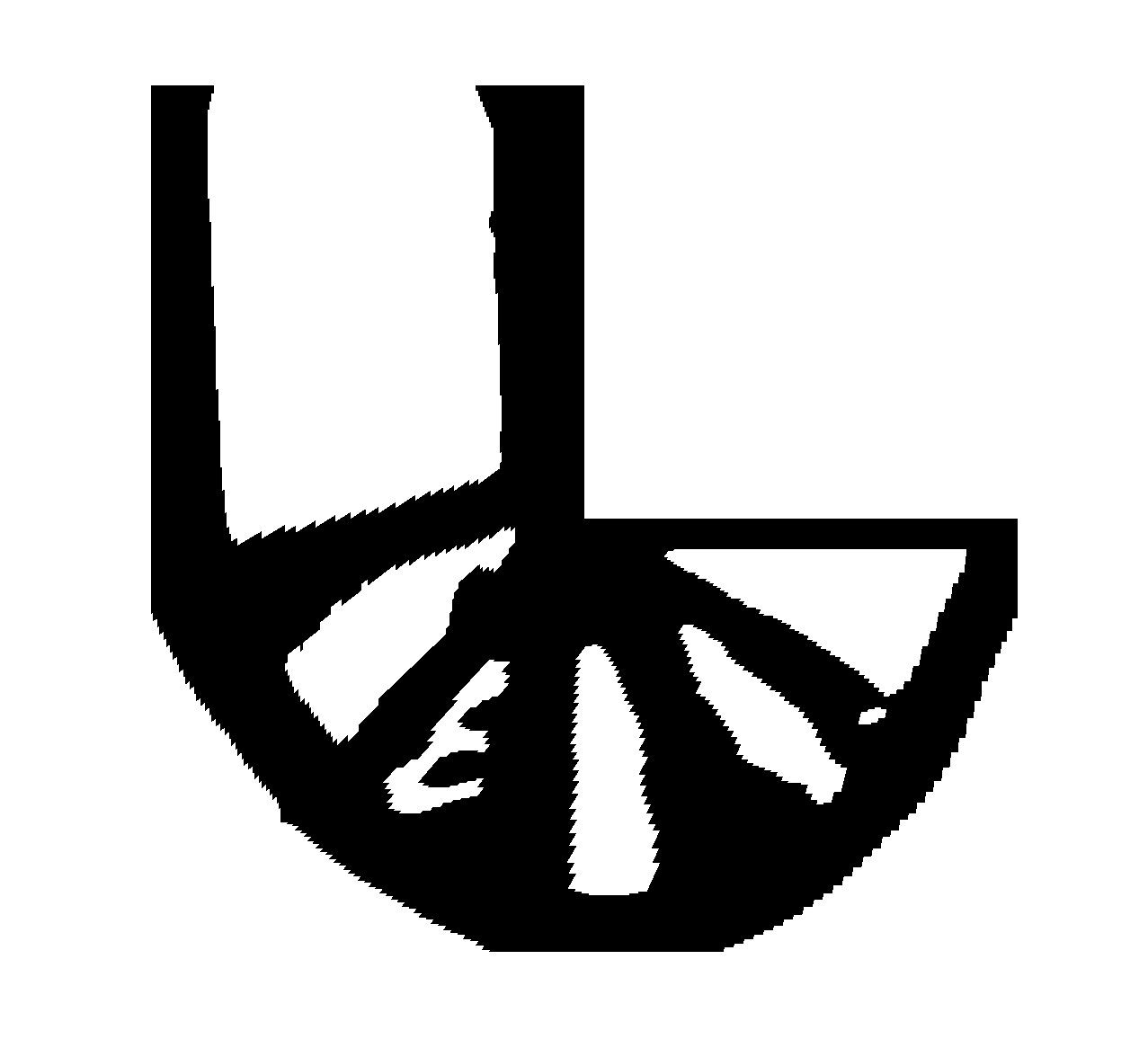}
	}
	\subfigure[iter=30: $N_{1}=4$]{
		\label{fig:ex2_iter_31}
		\includegraphics[width=0.29\linewidth]{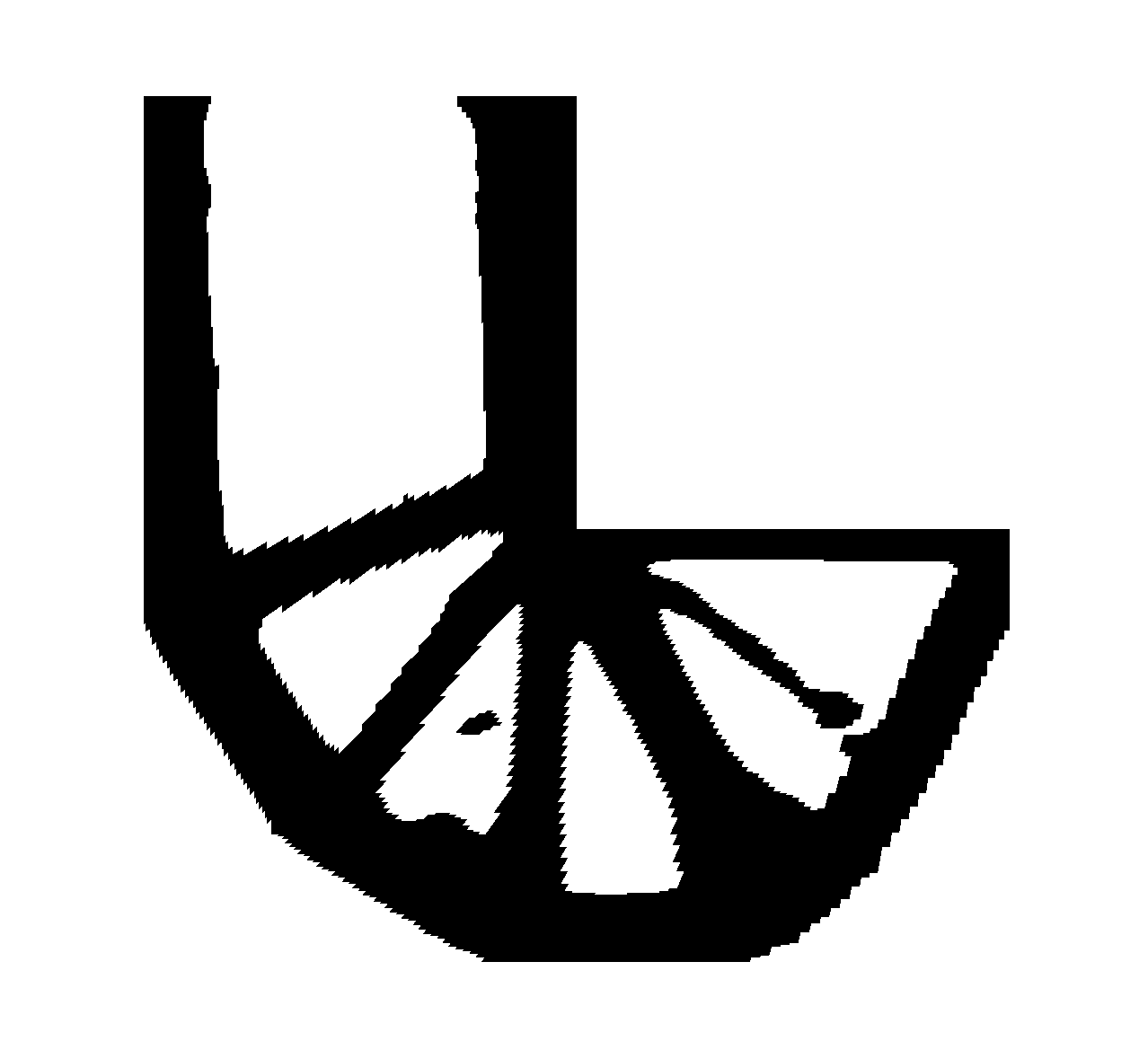}
	}
	\subfigure[iter=43: $N_{1}=4$]{
		\label{fig:ex2_iter_44}
		\includegraphics[width=0.29\linewidth]{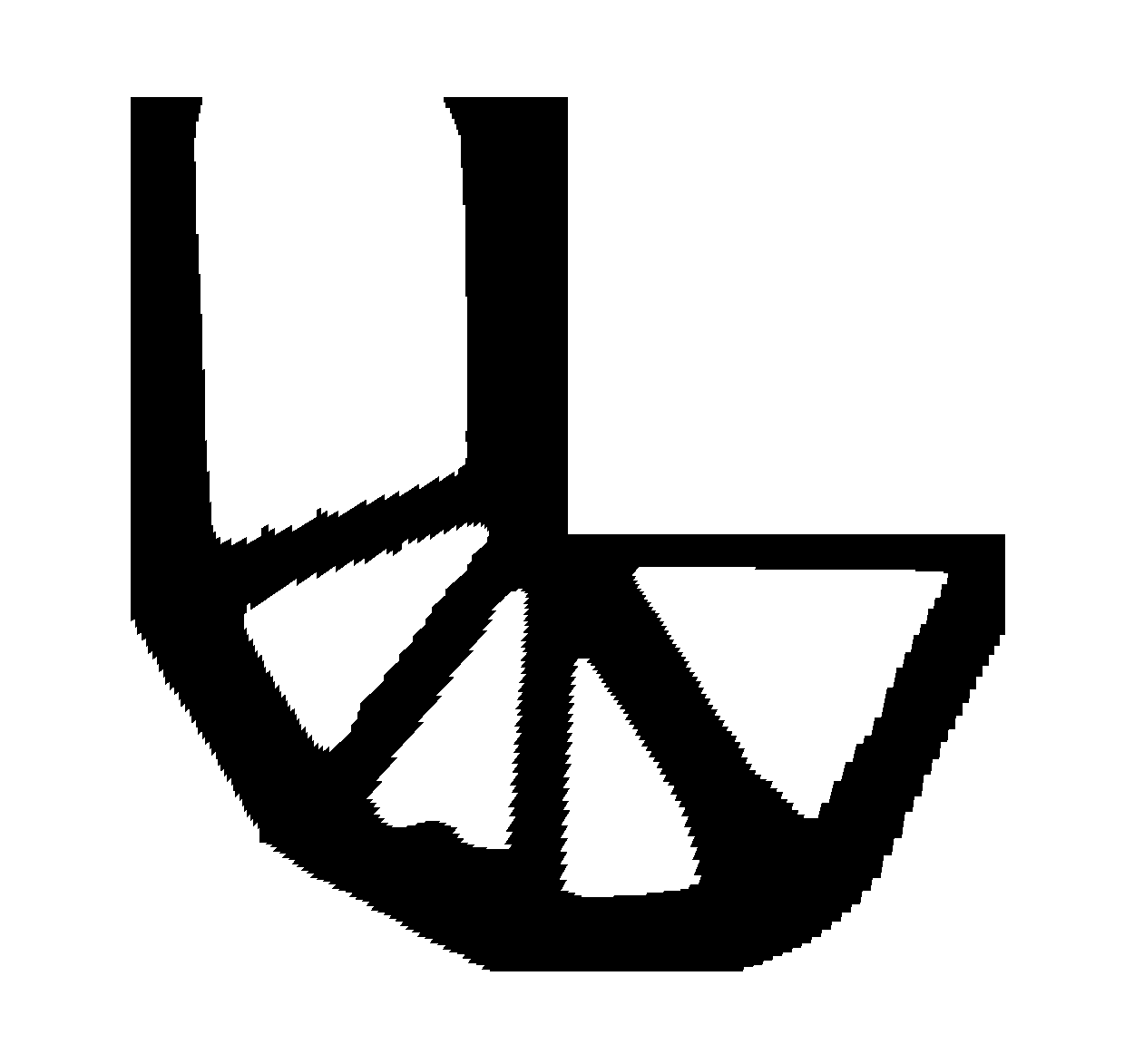}
	}
	\caption{Evolution of the optimized topology for the L-shaped beam when $\bar{N}_{1}=4$.
Density fields at representative iterations are shown.}
	\label{fig:L_process}
\end{figure}

Fig. \ref{fig:L_process} illustrates the convergence behavior 
	and topology evolution for the case $\bar{N}_{1} = 4$.
As shown in Fig.~\ref{fig:ex2_iter}, the compliance decreases rapidly 
	during the early iterations as the main load-carrying skeleton forms, 
	and then remains stable.
Importantly, although the compliance stabilizes, 
	the topology continues to evolve under the influence of the topological objective, 
	demonstrating that the proposed method can adjust topological features 
	without compromising convergence stability.
From Figs. \ref{fig:ex2_iter_11}-\ref{fig:ex2_iter_44}, 
	the number of holes increases to $9$ by iteration $19$, 
	which violates the prescribed constraint. 
Subsequently, redundant small holes are eliminated, 
	and by iteration 30, adjacent void regions merge through local structural reconfiguration, 
	reducing the hole count to satisfy $N_{1}\le 4$. 
A transient isolated component is observed at iteration 30, 
	but it is removed by iteration 43, 
	confirming the effectiveness of the 
	zero-dimensional objective in enforcing global connectivity.


\subsection{Cantilever beam example}
\begin{figure}[!htb]
	\centering
	\setcounter{subfigure}{0}
	\subfigure[]{
		\label{fig:mid_beam}
		\includegraphics[width=0.31\linewidth]{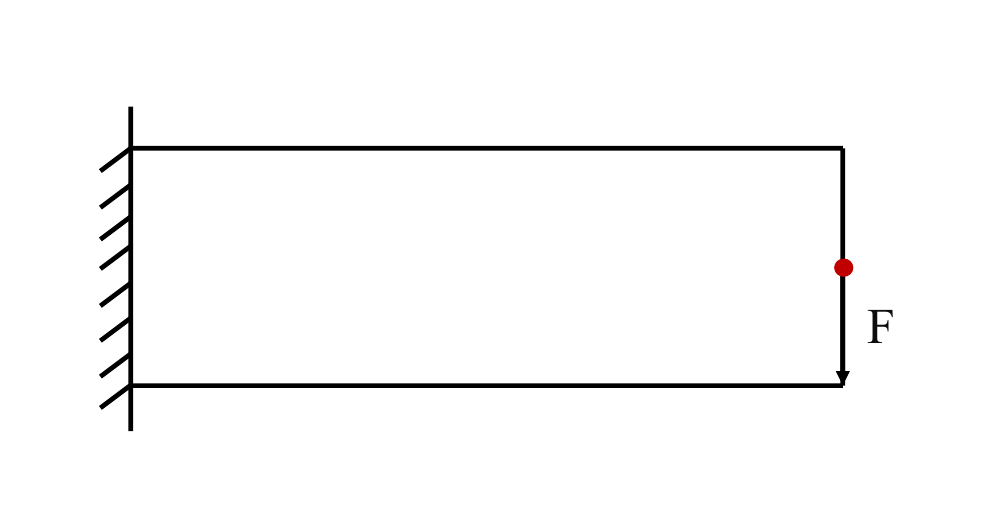}
	}
	\subfigure[]{
		\label{fig:mid_beam_density}
		\includegraphics[width=0.33\linewidth]{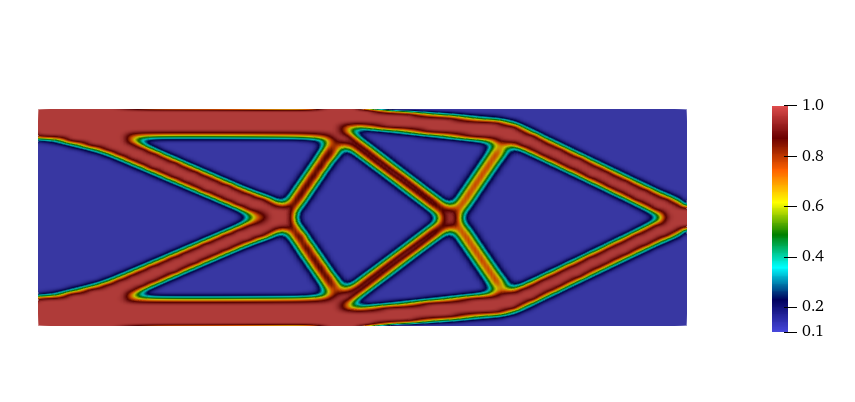}
	}
	\subfigure[]{
		\label{fig:mid_beam_no_opt_binar}
		\includegraphics[width=0.29\linewidth]{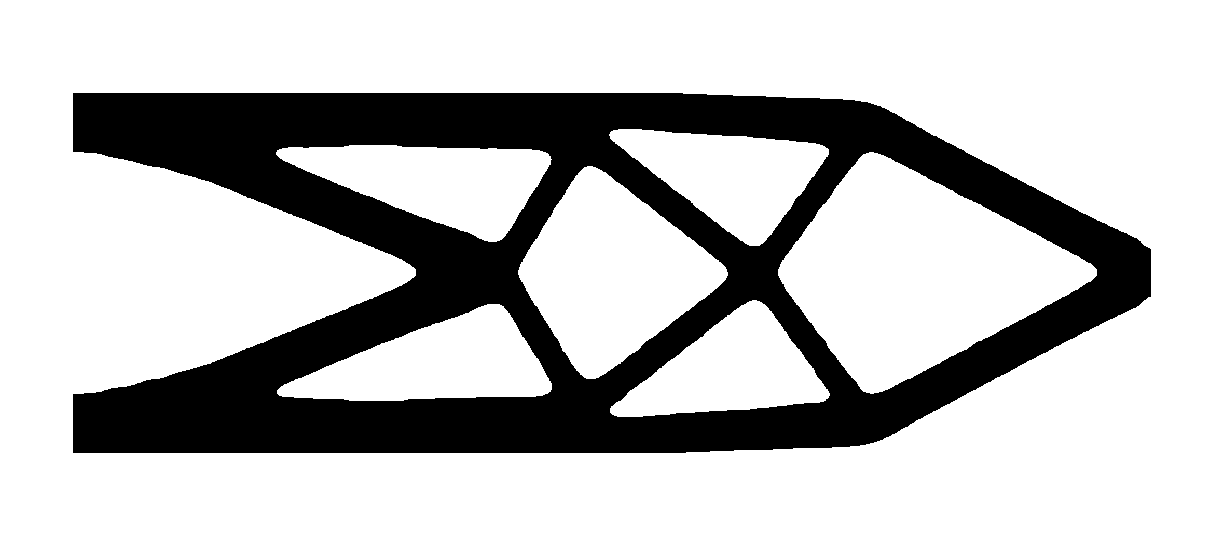}
	}
	\caption{(a) Cantilever beam model. (b) Optimization result without topology control. 
	(c) Binary structure corresponding to the density field in (b).}
	\label{fig:mid_beam_no_opt}
\end{figure}

In this subsection, the classic cantilever beam is considered. 
The beam length is three times its width. 
The boundary conditions are shown in Fig. \ref{fig:mid_beam}: 
	the left boundary is fixed, 
	while a downward force of $10^{6}$ is applied 
	at the midpoint of the right boundary.
The prescribed volume fraction is set to $0.4$.
In the IGA model, a bi-quadratic NURBS surface is employed to represent 
	both the design domain and density field. 
The number of control points in this model is $101 \times 51$, with $4704$ elements.

As shown in Figs. \ref{fig:mid_beam_density} and \ref{fig:mid_beam_no_opt_binar}, 
	there exist $6$ holes in the optimized structure without topological control.

\begin{figure}[!htb]
	\centering
	\setcounter{subfigure}{0}
	\subfigure[$\bar{N}_{1}=1$]{
	    \includegraphics[width=0.31\linewidth]{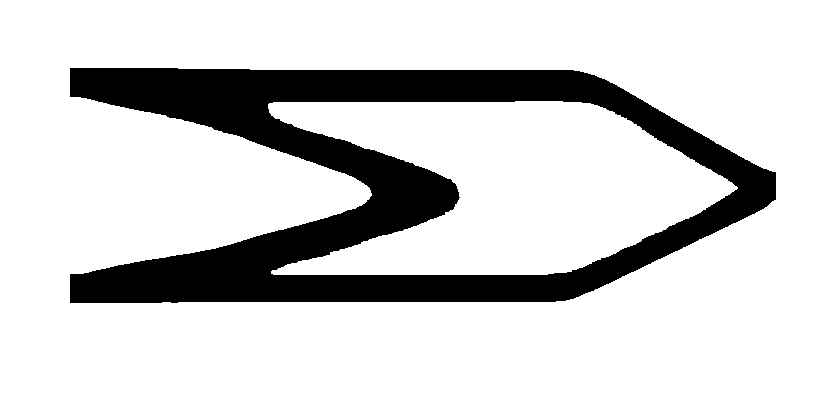}
	}
	\subfigure[$\bar{N}_{1}=2$]{
	     \includegraphics[width=0.31\linewidth]{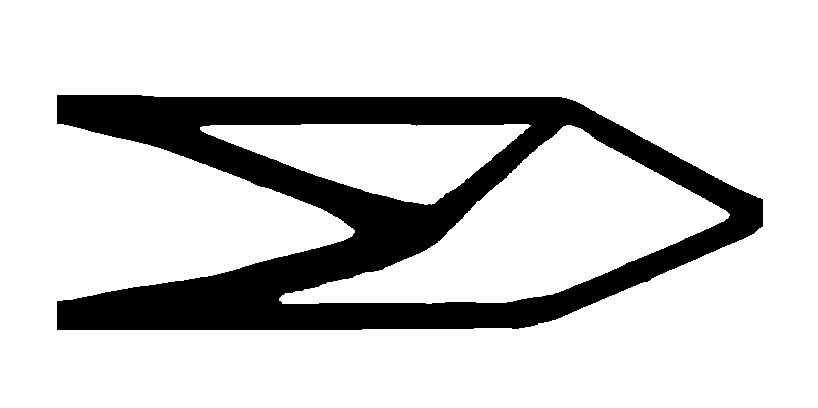}
	}
	\subfigure[$\bar{N}_{1}=3$]{
	     \includegraphics[width=0.31\linewidth]{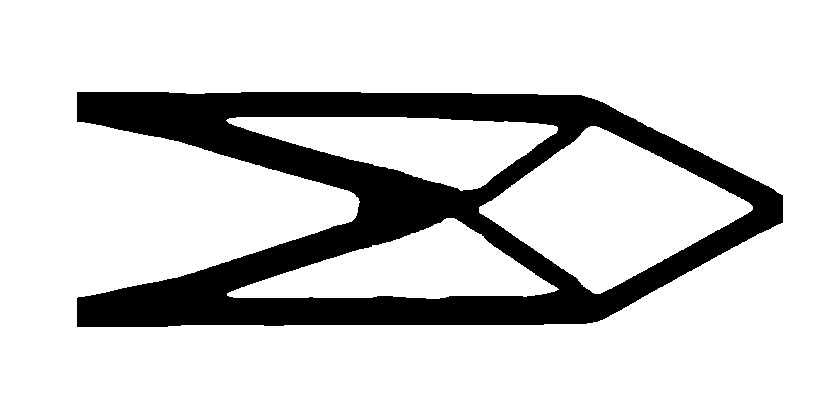}
	}\\
	\subfigure[$\bar{N}_{1}=4$]{
		\includegraphics[width=0.31\linewidth]{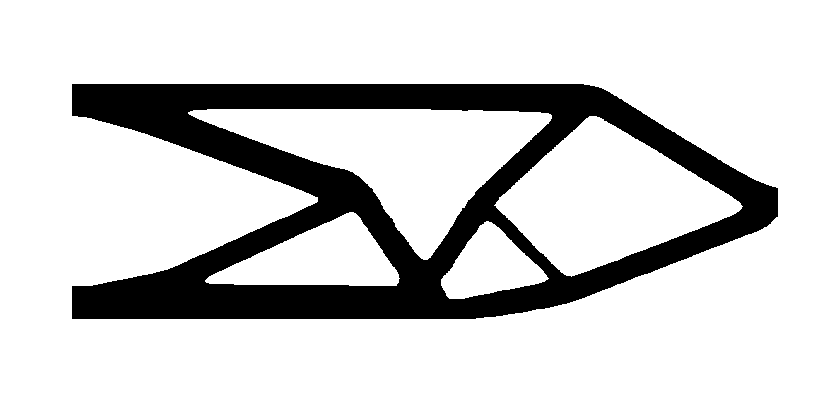}
	}
	\subfigure[$\bar{N}_{1}=5$]{
		\includegraphics[width=0.31\linewidth]{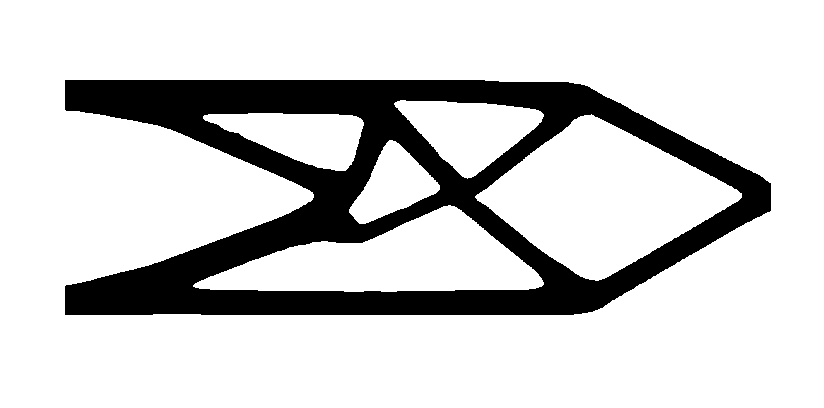}
	}
	\subfigure[$\bar{N}_{1}=6$]{
		\includegraphics[width=0.31\linewidth]{example_fig/ex3/ex3_vol00.png}
	}\\
	\subfigure[$\bar{N}_{1}=7$]{
		\includegraphics[width=0.31\linewidth]{example_fig/ex3/ex3_vol00.png}
	}
	\caption{Optimization results of the cantilever beam for $\bar{N}_{1}=1,2,3,4,5,6,7$, respectively.}
	\label{fig:m_beam_opt}
\end{figure}

Fig. \ref{fig:m_beam_opt} and Table \ref{tbl:obj_mid} summarize the topology-controlled results for $\bar{N}_{1}=1,2,3,4,5,6,7$. 
Overall, reducing the allowable number of holes leads to a moderate compliance increase, 
	consistent with the reduced admissible design space 
	under stronger topological constraints.
Notably, the unconstrained solution exhibits a symmetric layout (Fig.~\ref{fig:mid_beam_no_opt_binar}).
After topology control is imposed, 
	the optimized designs may become non-symmetric. 
This behavior is expected: topological constraints can break symmetry 
	by restricting admissible void configurations, 
	even when the loading and boundary conditions are symmetric. 
Importantly, the resulting designs remain mechanically sound and satisfy the prescribed topology.

\begin{table}[!htb]
	\caption{Comparison of compliance, number of holes, and time costs 
for the cantilever beam benchmark under different topology-control settings.}
	\begin{tabular}{llllllll}
		\toprule
		        &$\bar{N}_{1}=1$ & $\bar{N}_{1}=2$  & $\bar{N}_{1}=3$ & $\bar{N}_{1}=4$ & $\bar{N}_{1}=5$ & $\bar{N}_{1}=6$ & $\bar{N}_{1}=7$\\
		\midrule
		weights $(\mu_{0},\mu_{1})$    &(1.1,0.9)  &(1.1,0.9)  &(1.0,0.8)  &(0.9,0.7)  &(0.9,0.7)   & (0.4,0.2)   & (0.4,0.2)\\
		Compliance  &17.5  &17.3  &17.1 &16.1  & 16.5  &15.9    &15.9\\
		Time (s) &553.8  &539.7  &306.1  &231.1  & 146.5  &283.1    &285.2\\
		\bottomrule
	\end{tabular}
	\label{tbl:obj_mid}
\end{table}

\begin{figure}[!htb]
	\centering
	\setcounter{subfigure}{0}
	\subfigure[Compliance vs. iteration]{
	    \label{fig:ex3_iter}
		\includegraphics[width=0.3\linewidth]{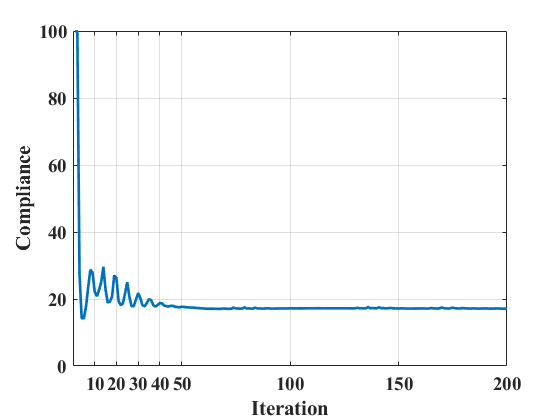}
	}
	\subfigure[iter=4: $N_{1}=0$]{
		\label{fig:ex3_iter_4}
		\includegraphics[width=0.31\linewidth]{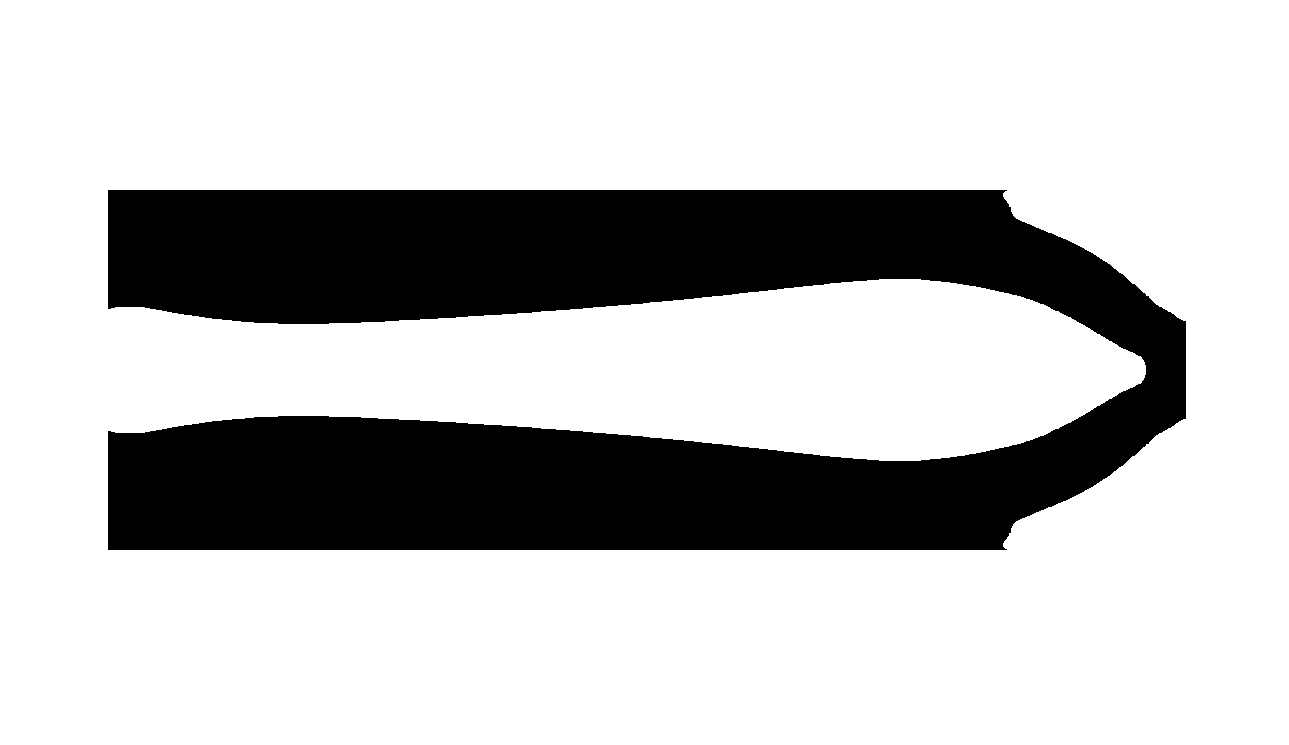}
	}
	\subfigure[iter=9: $N_{1}=2$]{
		\label{fig:ex3_iter_9}
		\includegraphics[width=0.31\linewidth]{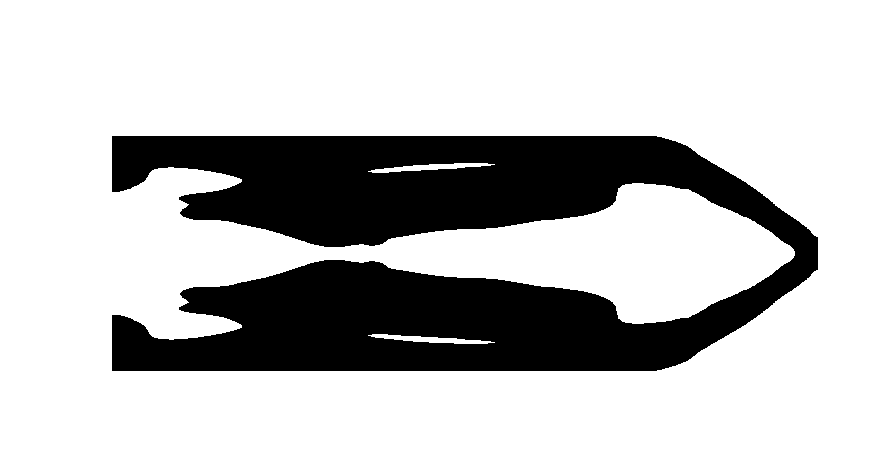}
	}\\
	\subfigure[iter=17: $N_{1}=2$]{
		\label{fig:ex3_iter_17}
		\includegraphics[width=0.31\linewidth]{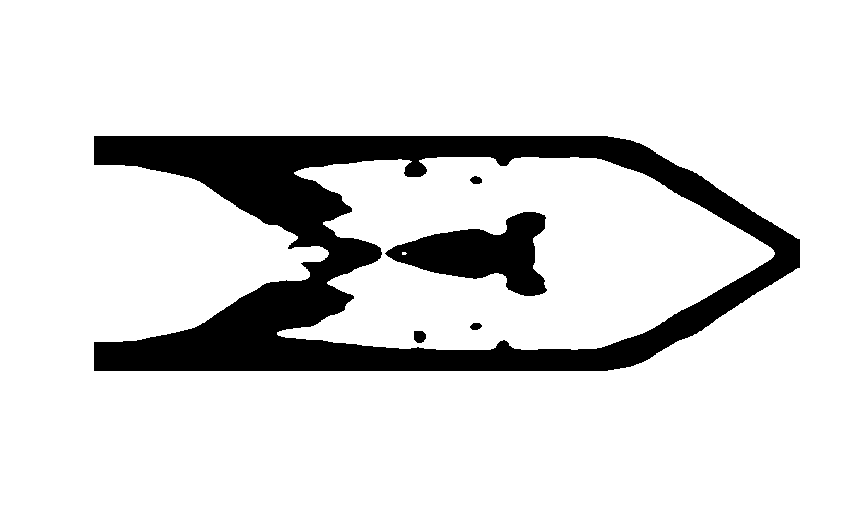}
	}
	\subfigure[iter=20: $N_{1}=12$]{
		\label{fig:ex3_iter_20}
		\includegraphics[width=0.32\linewidth]{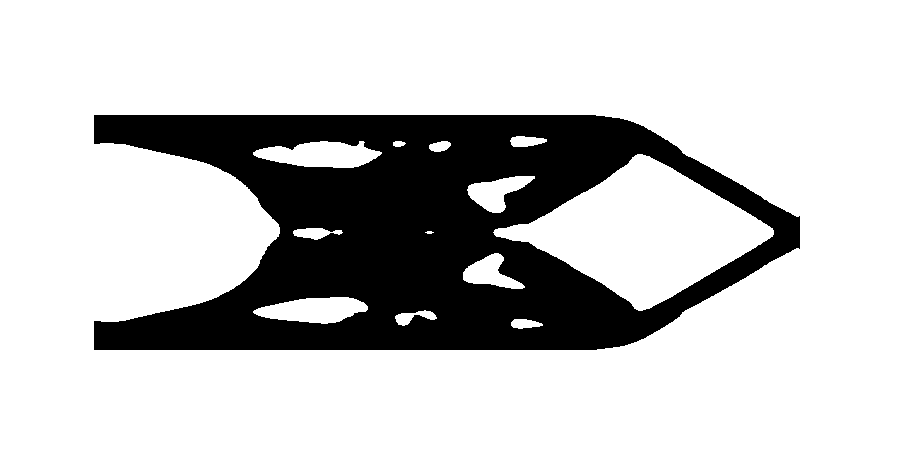}
	}
	\subfigure[iter=30: $N_{1}=8$]{
		\label{fig:ex3_iter_30}
		\includegraphics[width=0.30\linewidth]{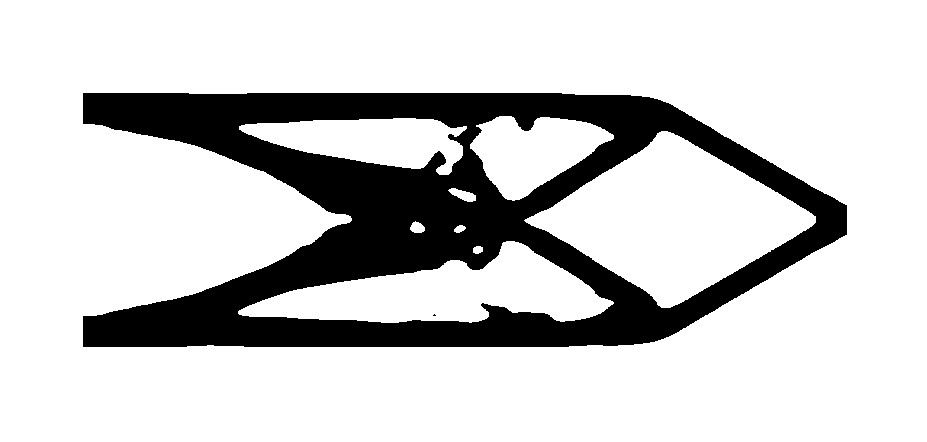}
	}\\
	\subfigure[iter=38: $N_{1}=3$]{
		\label{fig:ex3_iter_38}
		\includegraphics[width=0.31\linewidth]{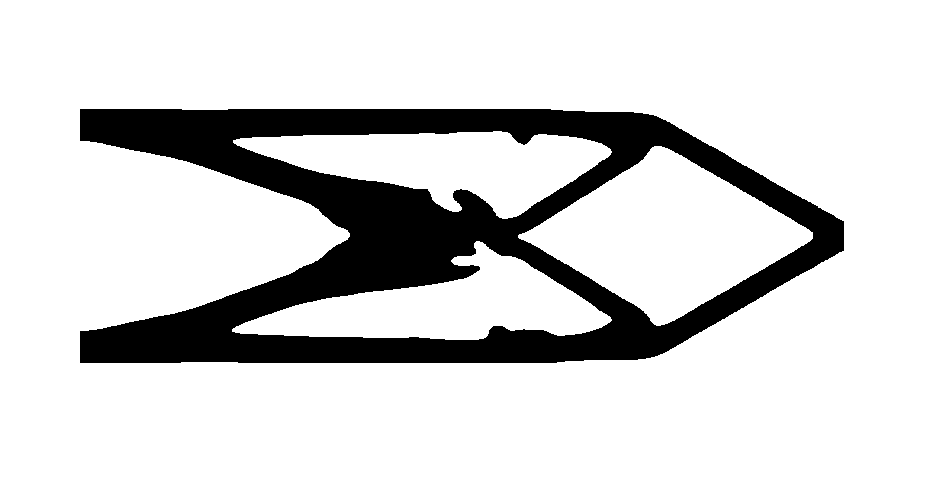}
	}
	\subfigure[iter=53: $N_{1}=3$]{
		\label{fig:ex3_iter_52}
		\includegraphics[width=0.32\linewidth]{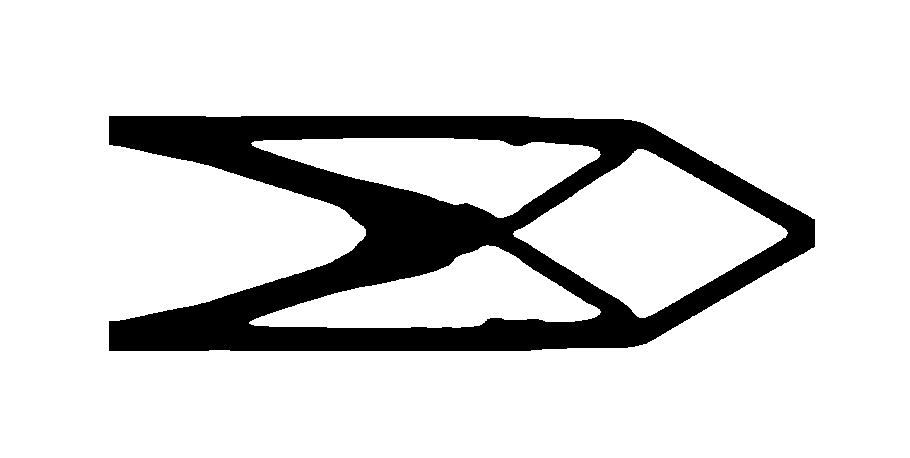}
	}
	\subfigure[iter=100: $N_{1}=3$]{
		\label{fig:ex3_iter_100}
		\includegraphics[width=0.30\linewidth]{example_fig/ex3/ex3_300.png}
	}
	\caption{Evolution of the optimized topology for the cantilever beam when $\bar{N}_{1}=3$.
Density fields at representative iterations are shown.}
	\label{fig:mid_process}
\end{figure}

Fig. \ref{fig:mid_process} illustrates the optimization process with 
	topology control $\bar{N}_{1}=3$.
In the early iterations, 
	the structure tends to generate holes and additional 
	connected components, 
	as the compliance dominates the objective function at this stage.
Fig. \ref{fig:ex3_iter} shows that the compliance exhibits 
	significant variation during the first $20$ iterations, 
	indicating that the main structure skeleton is formed during this stage.
Subsequently, the compliance becomes nearly stable 
	as iterations proceed, 
	with this phase primarily 
	focused on optimizing the number of holes.	
At the $17$th iteration, 
	the structure exhibits $5$ connected components 
	and $2$ holes; 
	at the $20$th iteration, 
	$12$ holes appear in the structure. 
After activating topology control at iteration 30, redundant small holes 
	disappear rapidly and local structural reconfigurations 
	(including the removal of ineffective thin members) 
	drive the design toward the prescribed hole constraint.
Specifically, the hole number decreases from 8 at iteration 30 to 3 at iteration 38, 
	satisfying $N_{1}\le 3$.
The compliance remains nearly constant during this topology-correction stage, 
	demonstrating that the proposed persistence-based objective regulates 
	topology without inducing numerical instability.

\begin{figure}[!htb]
	\centering
	\setcounter{subfigure}{0}
	\subfigure[]{
		\label{fig:hole_2_excess}
		\includegraphics[width=0.31\linewidth]{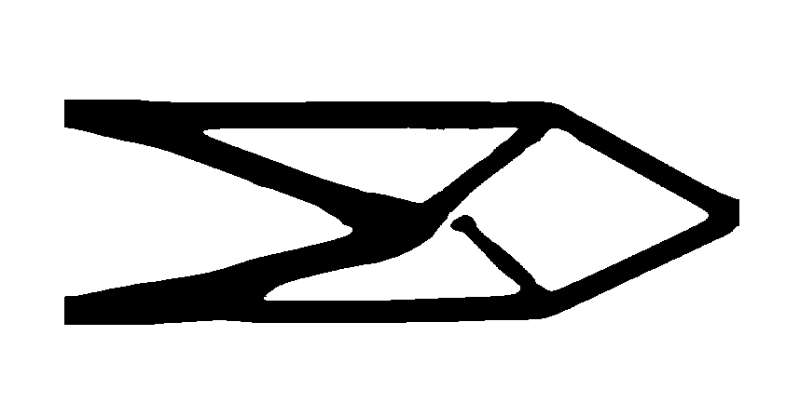}
	}
	\subfigure[]{
		\label{fig:hole_2}
		\includegraphics[width=0.31\linewidth]{example_fig/ex3/ex3_200.png}
	}
	\caption{(a) Optimized structure with excess material. 
	(b) Optimized structure after removing excess material.}
	\label{fig:hole_pre}
\end{figure}

Fig. \ref{fig:hole_2_excess} shows a typical situation 
	where the hole-number constraint is satisfied, 
	yet a small amount of excess material remains.
Although this material contributes marginally to stiffness, 
	it may persist due to the coupling between compliance minimization 
	and the topology objective.
To remove such ineffective material while preserving the achieved topology, 
	we freeze the density values within the identified holes 
	and continue optimizing compliance for the remaining design region. 
As shown in Fig. \ref{fig:hole_2}, the excess material is progressively eliminated 
	with negligible impact on structural performance. 
In Fig. \ref{fig:hole_2_excess}, the compliance of the structure is $17.4$, 
	and the compliance in Fig. \ref{fig:hole_2} is $17.3$.
This post-processing step improves manufacturability 
	while maintaining both topology constraints and mechanical efficiency.


\subsection{Quarter annulus example}
\begin{figure}[!htb]
	\centering
	\setcounter{subfigure}{0}
	\subfigure[]{
		\label{fig:quarter_annulus}
		\includegraphics[width=0.31\linewidth]{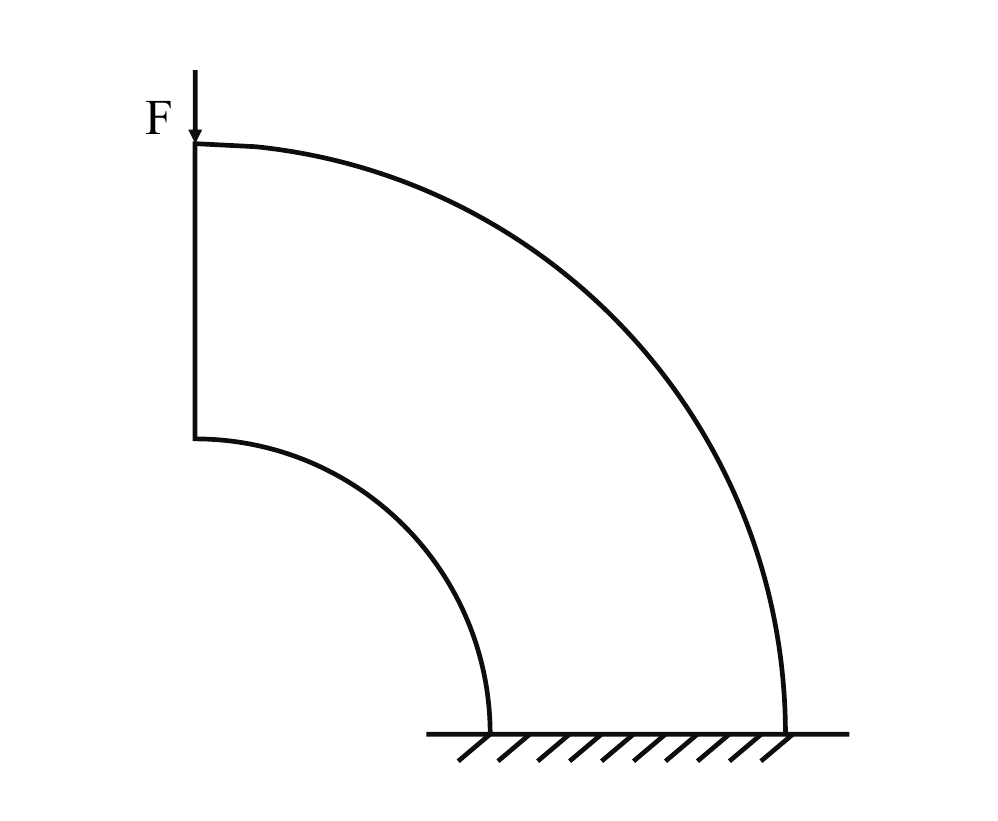}
	}
	\subfigure[]{
		\label{fig:quarter_density}
		\includegraphics[width=0.34\linewidth]{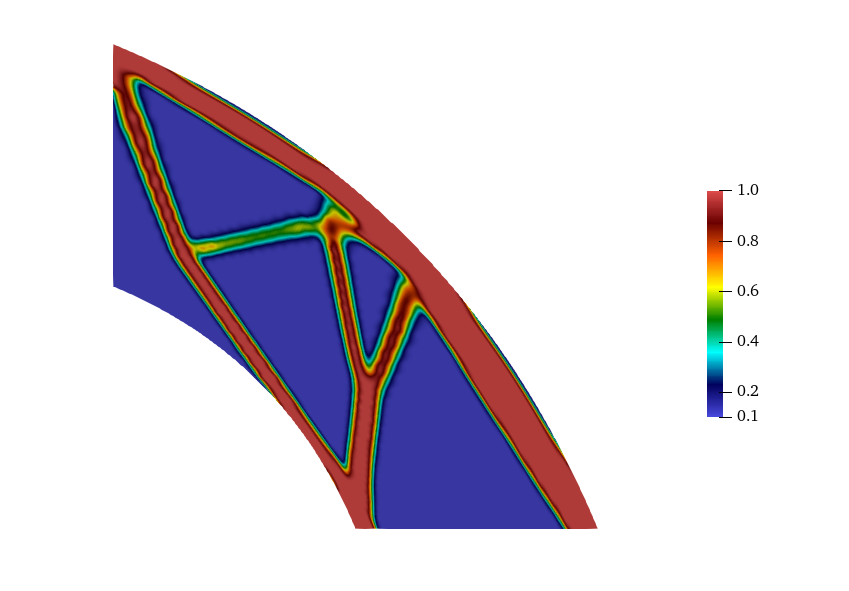}
	}
	\subfigure[]{
		\label{fig:quarter_binar}
		\includegraphics[width=0.27\linewidth]{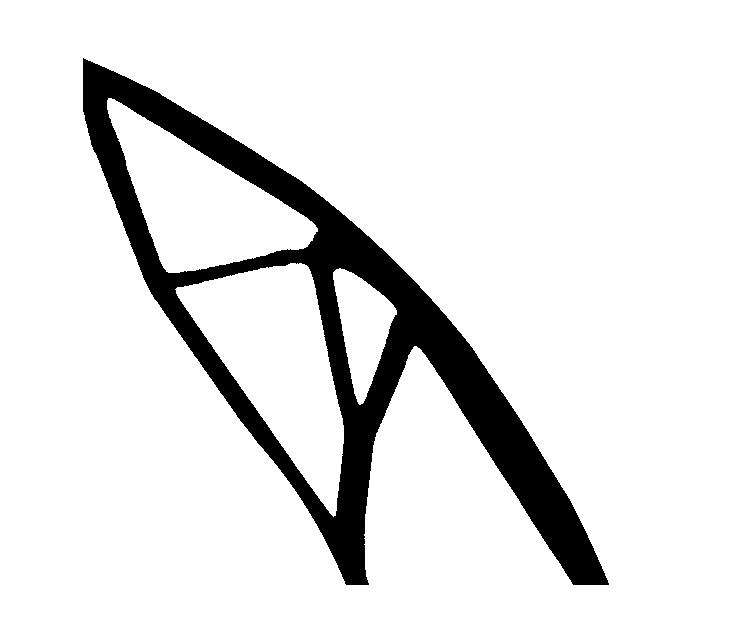}
	}
	\caption{(a) Quarter annulus model. (b) Optimization result without topology control. 
	(c) Binary structure corresponding to the density field in (b).}
	\label{fig:quarter_annulus_no_opt}
\end{figure}

In this subsection, we consider a quarter annulus with inner radius $r=5$ and 
	outer radius $R=10$, as shown in Fig. \ref{fig:quarter_annulus}. 
The bottom boundary of the design domain $\Omega$ is fully fixed, 
	and a downward load  
	of $10^{6}$ is applied at the top-right tip of $\Omega$.
The prescribed volume fraction is set to $V_{0}=0.4$. 
Both the design domain $\Omega$ and density field are 
	represented by NURBS bases of degree $2$. 
The discretization employs $62 \times 62$ control points and $3600$ elements.
With IGA, the exact geometric representation is 
	maintained throughout the analysis and optimization.

As shown in Figs. \ref{fig:quarter_density} and \ref{fig:quarter_binar}, 
	there exist $3$ holes in the optimized structure without topological control.

\begin{figure}[!htb]
	\centering
	\setcounter{subfigure}{0}
	\subfigure[$\bar{N}_{1}=1$]{
	    \includegraphics[width=0.31\linewidth]{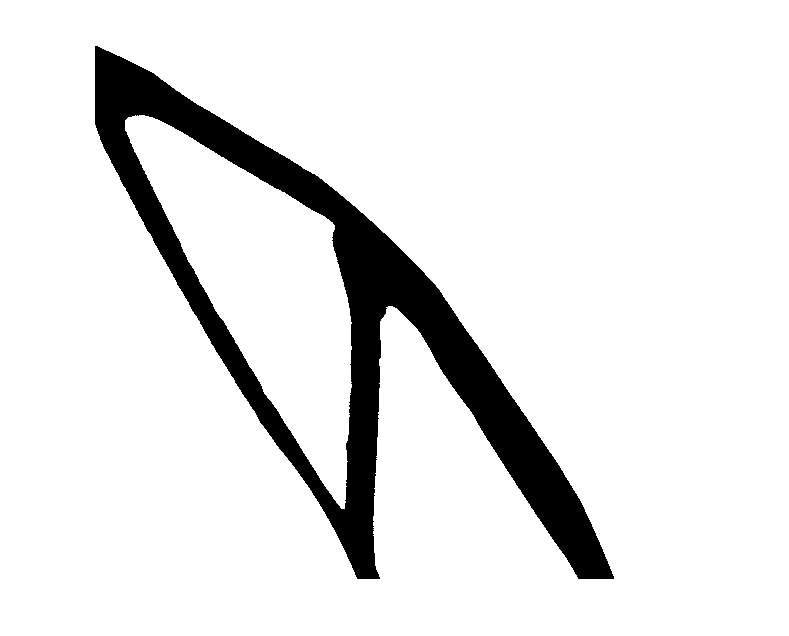}
	}
	\subfigure[$\bar{N}_{1}=2$]{
	     \includegraphics[width=0.31\linewidth]{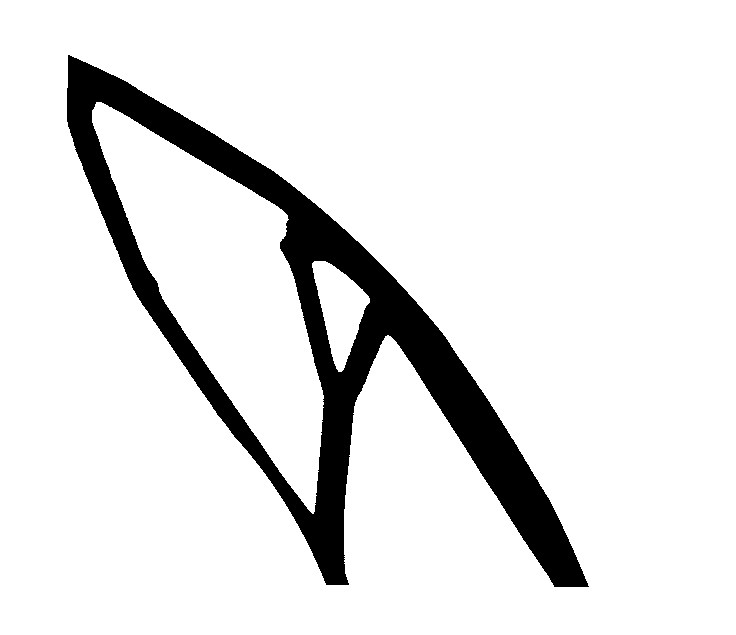}
	}\\
	\subfigure[$\bar{N}_{1}=3$]{
	     \includegraphics[width=0.31\linewidth]{example_fig/ex4/ex4_300.png}
	}
	\subfigure[$\bar{N}_{1}=4$]{
		\includegraphics[width=0.31\linewidth]{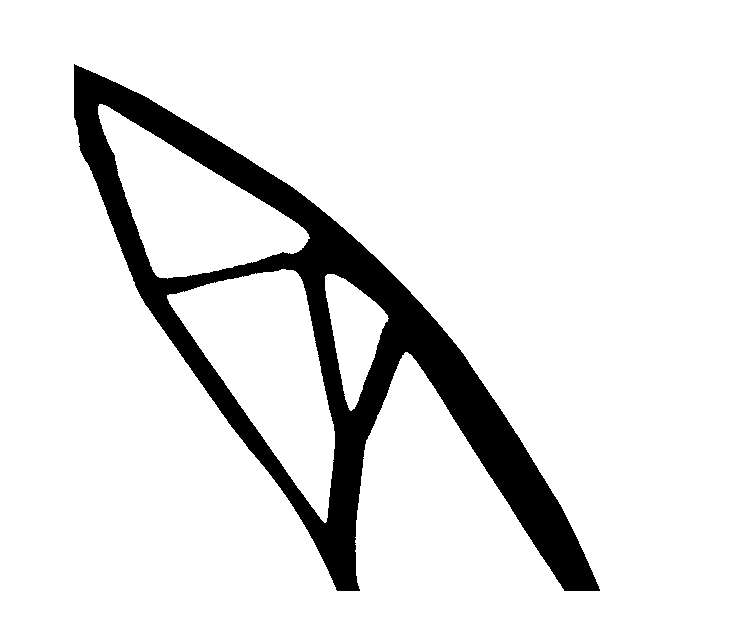}
	}
	\caption{Optimization results of the quarter annulus for $\bar{N}_{1}=1,2,3,4$, respectively.}
	\label{fig:quarter_opt}
\end{figure}

\begin{figure}[!htb]
	\centering
	\setcounter{subfigure}{0}
	\subfigure[Compliance vs. iteration]{
	\label{fig:ex4_iter}
		\includegraphics[width=0.33\linewidth]{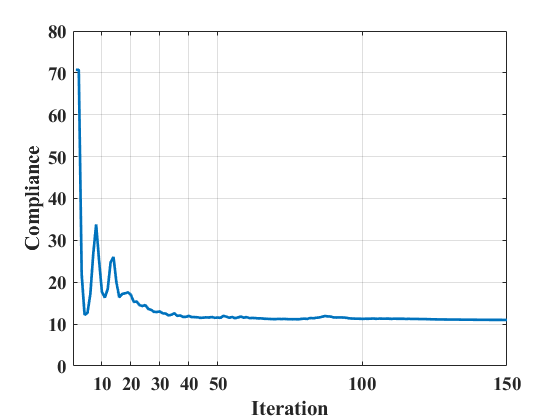}
	}
	\subfigure[iter=6: $N_{1}=0$]{
		\label{fig:ex4_iter_6}
		\includegraphics[width=0.29\linewidth]{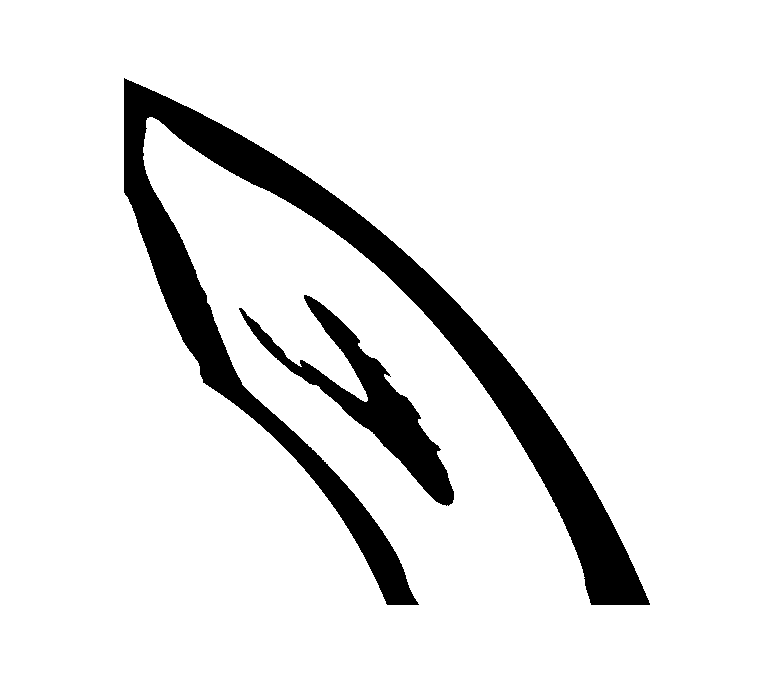}
	}
	\subfigure[iter=15: $N_{1}=2$]{
		\label{fig:ex4_iter_15}
		\includegraphics[width=0.29\linewidth]{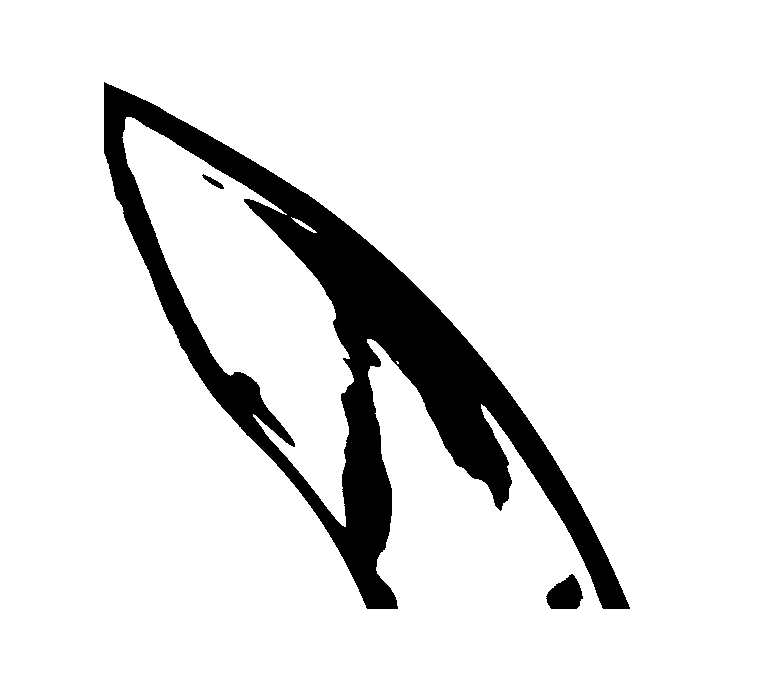}
	}\\
	\subfigure[iter=26: $N_{1}=4$]{
		\label{fig:ex4_iter_26}
		\includegraphics[width=0.29\linewidth]{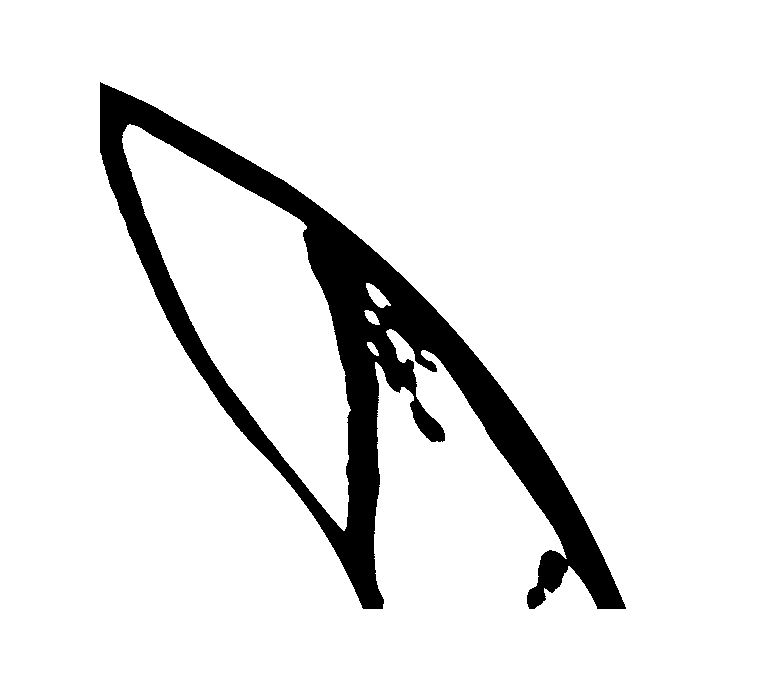}
	}
	\subfigure[iter=27: $N_{1}=1$]{
		\label{fig:ex4_iter_27}
		\includegraphics[width=0.29\linewidth]{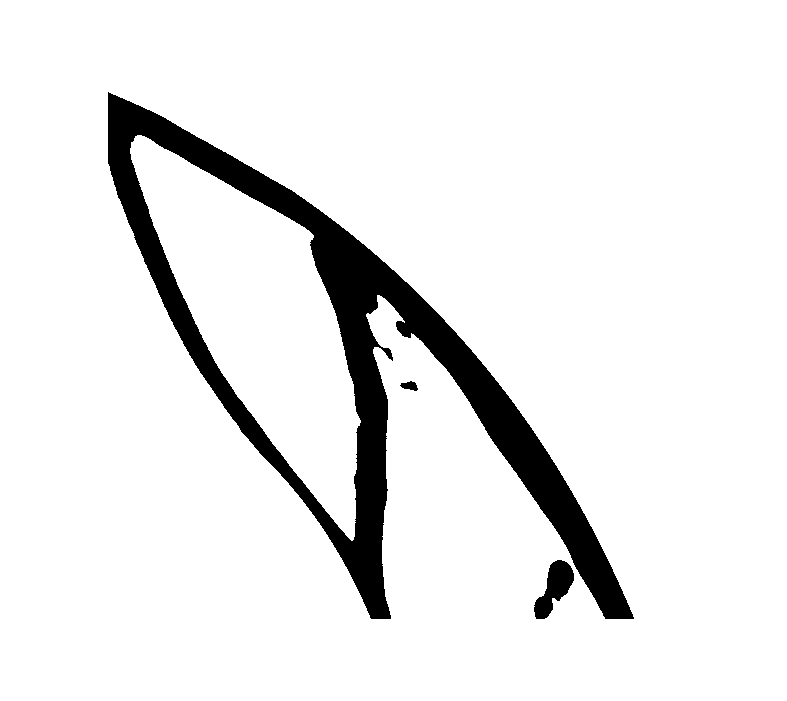}
	}
	\subfigure[iter=90: $N_{1}=1$]{
		\label{fig:ex4_iter_90}
		\includegraphics[width=0.29\linewidth]{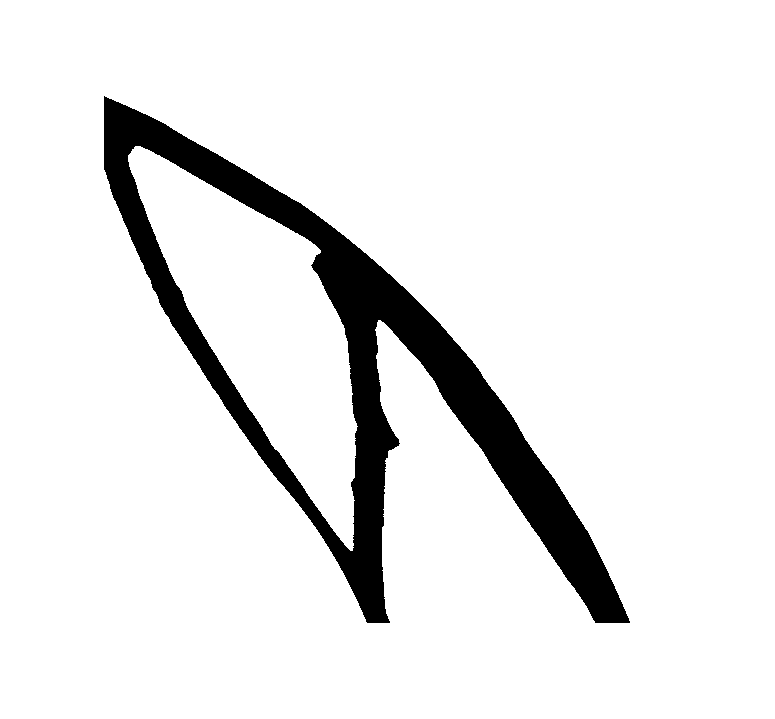}
	}
	\caption{Evolution of the optimized topology for the quarter annulus when $\bar{N}_{1}=1$.
Density fields at representative iterations are shown.}
	\label{fig:quarter_process}
\end{figure}

Fig. \ref{fig:quarter_opt} and Table \ref{tbl:obj_quarter} summarize 
	the topology-controlled results for $\bar{N}_{1}=1,2,3,4$.
A similar performance–complexity trade-off is observed. 
The compliance increases as the allowable number of holes decreases. 
The optimization process for $\bar{N}_{1}=1$ is shown in Fig.~\ref{fig:quarter_process}.
In Figs. \ref{fig:ex4_iter_26} and \ref{fig:ex4_iter_27}, 
	it can be observed that there exists $4$ holes 
	in the structure at $26$th iteration. 
After only one iteration, the hole number drops to 1 
	as small topological features are rapidly suppressed, 
	demonstrating the strong corrective capability of the proposed topology objective.

\begin{table}[!htb]
	\centering
	\caption{Comparison of compliance, number of holes, and time costs 
for the quarter annulus model benchmark under different topology-control settings.}
	\begin{tabular}{lllll}
		\toprule
		        &$\bar{N}_{1}=1$ & $\bar{N}_{1}=2$  & $\bar{N}_{1}=3$ & $\bar{N}_{1}=4$ \\
		\midrule
		weights $(\mu_{0},\mu_{1})$     &(0.4,0.3)  &(0.1,0.05)  &(0.05,0.05)  & (0.05,0.05)\\
		Compliance  &10.9  &10.7 &10.5  &10.5  \\
		Time (s) &271.3  &122.6  &189.2  &172.9  \\
		\bottomrule
	\end{tabular}
	\label{tbl:obj_quarter}
\end{table}

\section{Conclusion}
\label{sec:conclud}
In this work, we present a persistent-homology-based framework for 
	explicit topology control within density-based topology optimization. 
By incorporating persistence diagrams into the optimization process, 
	both structural connectivity and the number of holes 
	are rigorously quantified and explicitly regulated.
The resulting topology-aware objective is fully differentiable and can be seamlessly 
	integrated into gradient-based solvers, such as MMA.
Through a series of benchmark examples, 
	we demonstrate that the proposed method achieves explicit and reliable topology regulation 
	while maintaining competitive mechanical performance and stable convergence behavior.
Compared with conventional implicit strategies, 
	the present framework establishes a direct correspondence between 
	design variables and topological features, 
	thereby providing a systematic and mathematically grounded 
	mechanism for controlling structural complexity. 
Although the current study focuses on two-dimensional problems, 
	the proposed formulation is naturally extensible to three-dimensional settings. 
Future work will investigate large-scale three-dimensional applications,  
	improve computational efficiency for high-resolution discretizations, 
	and explore strategies for enforcing fully hole-free optimized structures 
	when required by manufacturing constraints.

\section*{Acknowledgements}

\bibliographystyle{elsarticle-harv} 
\bibliography{refer}

@article{bendsoe1988generating,
	title={Generating optimal topologies in structural design using a homogenization method},
	author={Bends{\o}e, Martin Philip and Kikuchi, Noboru},
	journal={Computer methods in applied mechanics and engineering},
	volume={71},
	number={2},
	pages={197--224},
	year={1988},
	publisher={Elsevier}
}

@article{li2018topology,
	title={Topology optimization design of cast parts based on virtual temperature method},
	author={Li, Quhao and Chen, Wenjiong and Liu, Shutian and Fan, Huiru},
	journal={Computer-Aided Design},
	volume={94},
	pages={28--40},
	year={2018},
	publisher={Elsevier}
}

@book{piegl2012nurbs,
	title={The NURBS book},
	author={Piegl, Les and Tiller, Wayne},
	year={2012},
	publisher={Springer Science \& Business Media}
}

@inproceedings{poulenard2018topological,
	title={Topological function optimization for continuous shape matching},
	author={Poulenard, Adrien and Skraba, Primoz and Ovsjanikov, Maks},
	booktitle={Computer Graphics Forum},
	volume={37},
	number={5},
	pages={13--25},
	year={2018},
	organization={Wiley Online Library}
}

@inproceedings{bruel2020topology,
	title={Topology-Aware Surface Reconstruction for Point Clouds},
	author={Br{\"u}el-Gabrielsson, Rickard and Ganapathi-Subramanian, Vignesh and Skraba, Primoz and Guibas, Leonidas J},
	booktitle={Computer graphics forum},
	volume={39},
	number={5},
	pages={197--207},
	year={2020},
	organization={Wiley Online Library}
}

@article{svanberg1987method,
	title={The method of moving asymptotes—a new method for structural optimization},
	author={Svanberg, Krister},
	journal={International journal for numerical methods in engineering},
	volume={24},
	number={2},
	pages={359--373},
	year={1987},
	publisher={Wiley Online Library}
}

@article{sigmund1998numerical,
	title={Numerical instabilities in topology optimization: a survey on procedures dealing with checkerboards, mesh-dependencies and local minima},
	author={Sigmund, Ole and Petersson, Joakim},
	journal={Structural optimization},
	volume={16},
	number={1},
	pages={68--75},
	year={1998},
	publisher={Springer}
}

@phdthesis{sigmund1994design,
  title={Design of material structures using topology optimization},
  author={Sigmund, Ole},
  year={1994},
  school={Technical University of Denmark Lyngby}
}

@article{bourdin2001filters,
  title={Filters in topology optimization},
  author={Bourdin, Blaise},
  journal={International journal for numerical methods in engineering},
  volume={50},
  number={9},
  pages={2143--2158},
  year={2001},
  publisher={Wiley Online Library}
}

@article{zhou2015minimum,
  title={Minimum length scale in topology optimization by geometric constraints},
  author={Zhou, Mingdong and Lazarov, Boyan S and Wang, Fengwen and Sigmund, Ole},
  journal={Computer Methods in Applied Mechanics and Engineering},
  volume={293},
  pages={266--282},
  year={2015},
  publisher={Elsevier}
}

@article{wang2003level,
  title={A level set method for structural topology optimization},
  author={Wang, Michael Yu and Wang, Xiaoming and Guo, Dongming},
  journal={Computer methods in applied mechanics and engineering},
  volume={192},
  number={1-2},
  pages={227--246},
  year={2003},
  publisher={Elsevier}
}

@article{liu2015identification,
  title={An identification method for enclosed voids restriction in manufacturability design for additive manufacturing structures},
  author={Liu, Shutian and Li, Quhao and Chen, Wenjiong and Tong, Liyong and Cheng, Gengdong},
  journal={Frontiers of Mechanical Engineering},
  volume={10},
  number={2},
  pages={126--137},
  year={2015},
  publisher={Springer}
}

@article{yamada2022topology,
  title={Topology optimization with a closed cavity exclusion constraint for additive manufacturing based on the fictitious physical model approach},
  author={Yamada, Takayuki and Noguchi, Yuki},
  journal={Additive Manufacturing},
  volume={52},
  pages={102630},
  year={2022},
  publisher={Elsevier}
}

@article{zhang2017structural,
	title={Structural complexity control in topology optimization via moving morphable component (MMC) approach},
	author={Zhang, Weisheng and Zhou, Jianhua and Zhu, Yichao and Guo, Xu},
	journal={Structural and Multidisciplinary Optimization},
	volume={56},
	number={3},
	pages={535--552},
	year={2017},
	publisher={Springer}
}

@article{zhang2017explicit,
	title={Explicit control of structural complexity in topology optimization},
	author={Zhang, Weisheng and Liu, Ying and Wei, Peng and Zhu, Yichao and Guo, Xu},
	journal={Computer Methods in Applied Mechanics and Engineering},
	volume={324},
	pages={149--169},
	year={2017},
	publisher={Elsevier}
}

@article{han2021topological,
	title={Topological constraints in 2D structural topology optimization},
	author={Han, Haitao and Guo, Yuchen and Chen, Shikui and Liu, Zhenyu},
	journal={Structural and Multidisciplinary Optimization},
	volume={63},
	number={1},
	pages={39--58},
	year={2021},
	publisher={Springer}
}

@article{zhao2020direct,
	title={A direct approach to controlling the topology in structural optimization},
	author={Zhao, Zi-Long and Zhou, Shiwei and Cai, Kun and Xie, Yi Min},
	journal={Computers \& Structures},
	volume={227},
	pages={106141},
	year={2020},
	publisher={Elsevier}
}

@article{zuo2022explicit,
  title={Explicit 2D topological control using SIMP and MMA in structural topology optimization},
  author={Zuo, Tongxing and Wang, Chong and Han, Haitao and Wang, Qianglong and Liu, Zhenyu},
  journal={Structural and Multidisciplinary Optimization},
  volume={65},
  number={10},
  pages={293},
  year={2022},
  publisher={Springer}
}

@inproceedings{carriere2021optimizing,
  title={Optimizing persistent homology based functions},
  author={Carriere, Mathieu and Chazal, Fr{\'e}d{\'e}ric and Glisse, Marc and Ike, Yuichi and Kannan, Hariprasad and Umeda, Yuhei},
  booktitle={International conference on machine learning},
  pages={1294--1303},
  year={2021},
  organization={PMLR}
}

@article{clough2020topological,
  title={A topological loss function for deep-learning based image segmentation using persistent homology},
  author={Clough, James R and Byrne, Nicholas and Oksuz, Ilkay and Zimmer, Veronika A and Schnabel, Julia A and King, Andrew P},
  journal={IEEE transactions on pattern analysis and machine intelligence},
  volume={44},
  number={12},
  pages={8766--8778},
  year={2020},
  publisher={IEEE}
}

@book{edelsbrunner2010computational,
  title={Computational topology: an introduction},
  author={Edelsbrunner, Herbert and Harer, John},
  year={2010},
  publisher={American Mathematical Soc.}
}

\end{document}